
\input amstex 
\documentstyle{amsppt}
\magnification1200
\NoRunningHeads
\refstyle{A}


\let\la=\langle 
\let\ra=\rangle 

\define\1{{\bold 1}}
\define\C{{\Bbb C}}
\define\Z{{\Bbb Z}}
\define\N{{\Bbb N}}

\define\F{{\Bbb F}}

\define\R{{\Cal R}}
\define\Rhat{{\hat{\Cal R}}}
\define\g{\frak{g}}
\define\ghat{\hat{\frak{g}}}
\define\h{\frak{h}}

\define\End{\operatorname{End}}
\define\Hom{\operatorname{Hom}}
\define\rank{\operatorname{rank}}
\define\Res{\operatorname{Res}}
\define\coeff{\operatorname{coeff}}
\define\Sym{\operatorname{Sym}}
\define\wt{\operatorname{wt}}
\define\sh{\operatorname{sh}}
\define\lt{\operatorname{\ell \!\text{{\it t\,}}}}
\define\lspan{\operatorname{\Bbb F\text{-span}}}
\define\ch{\operatorname{ch}}


\topmatter
\title{Annihilating fields of standard modules of
$\frak{sl}(2,\C)\sptilde$ and combinatorial
  identities}
\endtitle
\author Arne Meurman and Mirko Primc
\endauthor
\address
Univ\. of Lund,
Dept\. of Mathematics, 
S-22100 Lund,
Sweden
\endaddress
\email arnem\@maths.lth.se
\endemail
\address
Univ\. of Zagreb,
Dept\. of Mathematics,
Zagreb,
Croatia
\endaddress
\email primc\@cromath.math.hr 
\endemail
\thanks Partially supported by a grant from the G\"oran
Gustafsson
Foundation for Research in Natural Sciences and Medicine and
by the Ministry of Science 
and Technology of the Republic of Croatia
  \endthanks
\keywords
affine Lie algebras, vertex operator algebras, vertex
operator formula, standard modules, 
loop modules, Rogers-Ramanujan identities, colored
partitions, partition ideals
  \endkeywords
\subjclass Primary 17B67;
Secondary 05A19
\endsubjclass
\abstract
  We show that a set of local admissible fields generates a
vertex algebra. For an
affine Lie algebra $\tilde\goth g$ we construct the
corresponding level $k$ vertex operator algebra 
and we show that level $k$  highest weight $\tilde\goth
g$-modules are modules for this
vertex operator algebra. We determine the set of
annihilating fields of level $k$ standard modules and
we study the corresponding loop $\tilde\goth g$ module---the
set of relations that defines
standard modules. In the case when $\tilde\goth g$ is of
type $A_1^{(1)}$, we construct
bases of standard modules parameterized by colored
partitions and, as a consequence, we
obtain a series of Rogers-Ramanujan type combinatorial
  identities.  
\endabstract
\toc
\subhead  {} \ \ \ \ \ \   Introduction  \endsubhead
\subhead  1. Formal Laurent series and rational functions 
\endsubhead
\subhead  2. Generating fields \endsubhead
\subhead  3. The vertex operator algebra $N(k\Lambda_0)$
\endsubhead
\subhead  4. Modules over $N(k\Lambda_0)$ \endsubhead
\subhead  5. Relations on standard modules \endsubhead
\subhead  6. Colored partitions, leading terms and the main
results \endsubhead
\subhead  7. Colored partitions allowing at least two
embeddings \endsubhead 
\subhead  8. Relations among relations \endsubhead
\subhead  9. Relations among relations for two
embeddings\endsubhead
\subhead  10. Linear independence of bases of standard
modules\endsubhead
\subhead  11. Some combinatorial identities of
Rogers-Ramanujan type \endsubhead
\subhead   {} \ \ \ \ \ \ References \endsubhead
\endtoc

\endtopmatter


\document

\head{ Introduction}\endhead

In this paper we give an explicit construction of standard
(i.e. integrable highest
weight) representations of affine Lie algebra $\tilde\goth
g$ of the type $A_1^{(1)}$, and as 
a consequence we get a series of Rogers-Ramanujan type
identities. 

In order to describe our main result, let
$\g=\frak{sl}(2,\Bbb C )$ and let as usual
$\tilde\goth g=\g\otimes\Bbb C[t,t^{-1}]\oplus \Bbb C c
\oplus \Bbb C d$ be the corresponding 
affine Lie algebra. We fix the usual sl${}_2$-basis $\{x, h,
y\}$ in $\g$. 
For an element $X\in \g$ we write $X(n)=X\otimes t^n$. Set
$$
\bar B=\{b(n) \mid b\in\{x,h,y\}, \ n\in \Bbb Z \},\quad 
\bar B_-=\{b(n)\in \bar B \mid n\le 0, \ n<0 \text{  if 
}b\neq y\}.
$$
Then $\bar B_-$ is a basis of subalgebra $\tilde{\goth n}_-$
of affine  Lie algebra 
$\tilde\goth g$. We define a linear order $\preccurlyeq$ on
$\bar B$ by 
$\dots\prec x(j-1)\prec y(j)\prec h(j)\prec x(j)\prec
y(j+1)\prec \dots$ and we write 
$$
  u(\pi)=\prod_{a\in \bar B}\, a^{\pi(a)}= a_1\cdots a_s
\tag{1}
$$
to denote a (finite) product in $U(\tilde{\goth g})$ of
elements $a_j\in \bar B$, 
$a_1 \preccurlyeq \dots \preccurlyeq a_s$, where
each factor $a\in \bar B$ of $u(\pi)$ appears precisely
$\pi(a)$ times. We may think of
such product $u(\pi)$ in the universal enveloping algebra 
simply as a
colored partition $\pi \: \bar B \rightarrow \N$
consisting of $\pi(a)$ parts $a$, 
each part $a=b(n)$ having a degree
$n$ and color $b\in\{x,h,y\}$. We denote the set of all
colored partitions by
$\Cal P(\bar B)$, and by $\Cal P(\bar B_-)$ we denote all
colored partitions with parts in 
$\bar B_-$.  

The set $\{u(\pi) \mid \pi \in \Cal P(\bar B_-)\}$ is a 
Poincar\'e-Birkhoff-Witt basis of $U(\tilde{\goth n}_-)$,
so for a highest weight $\tilde\goth g$-module $V$ with a
highest weight vector $v_0$ the set of 
vectors
$$
   u(\pi)\cdot v_0,\qquad \pi\in\Cal P(\bar B_-) ,\tag{2}
$$
is a spanning set of $V$. In the case that $V=M(\Lambda)$ is
a Verma module,
the set (2) is a basis. In the case that $V=L(\Lambda)$ is a
standard
$\tilde\goth g$-module, the spanning set (2) is not a basis,
and our
main result is that for $\Lambda
=k_0\Lambda_0+k_1\Lambda_1$, $k_0, k_1\in \N$, the 
set of vectors of the form (2) subject to conditions
$$
\gathered
\pi (y(j-1)) + \pi (h(j-1)) + \pi (y(j))  \le k,\\
\pi (h(j-1)) + \pi (x(j-1)) + \pi (y(j))  \le k,\\
\pi (x(j-1)) + \pi (y(j)) + \pi (h(j)) \le k,\\
\pi (x(j-1)) + \pi (h(j)) + \pi (x(j))  \le k,\\ 
\pi (x(-1)) \le k_0, \quad \pi (y(0))  \le k_1
\endgathered \tag{3}
$$
is a basis of standard module $L(\Lambda)$ of level 
$k=k_0+k_1$.

As a consequence of this result we get a series of
combinatorial identities 
of Rogers-Ramanujan type. For example, if we write a
partition of a nonnegative integer $n$
in the form  $n=\sum_{j\ge 0} \,jf_j$, then the number of
partitions of
$n$ such that each part appears at most once (i.e\. $f_j\le 1$
for all $j$) equals the number of partitions of $n$ such that each
part appears at most twice (i.e\. $f_j\le 2$ for all $j$), but subject
to the conditions
$$
\gathered
f_{3j+2} + f_{3j+1} + f_{3j}  \le 2,\\
f_{3j+2} + f_{3j} + f_{3j-1}  \le 2,\\
f_{3j+1} + f_{3j} + f_{3j-2}  \le 2,\\
f_{3j} + f_{3j-1} + f_{3j-2}  \le 2,\\
f_1 \le  1, \quad f_2 \le  1 .
\endgathered \tag{4}
$$

Before we say something about the proof itself, let us say a
few words about 
some general ideas involved.

If $\tilde\goth g$ is an affine Lie algebra, then a standard
$\tilde\goth g$-module 
$L(\Lambda)$ is a
quotient of the Verma module $M(\Lambda)$ by its maximal
submodule $M^1(\Lambda)$,
and this construction works for any Kac-Moody Lie algebra.
On the other hand, for affine
Lie algebras we also have various vertex operator
constructions of basic (or level~1) 
representations: The  representation space is explicitly given,
usually the bosonic
or fermionic Fock space, and the 
action of the affine Lie algebra is given by the vertex operator
formula, roughly of the form
$X(z)=E^-(z)E^+(z)E^0(z)$, where $E^\pm(z)$ are some
exponentials of Heisenberg Lie algebra 
elements, and $E^0(z)$ is an exponential involving ``group
elements" (in principle 
defining operators only on integrable representations).

J.Lepowsky and R.Wilson gave a Lie-theoretic interpretation
and proof of the Rogers-Ramanujan
identities by explicitly (explicitly in a broader sense)
constructing level~3 standard 
representations of $\goth {sl}(2,\Bbb C)\sptilde$: The
vertex operator 
formula in the principal picture for level 1 standard modules
implies (very roughly speaking)
a formula of the form $X(z)^2=E^-(z)X(z)E^+(z)E^0(z)$
relating the action of various
Lie algebra elements on level 3 standard modules. In this
case one cannot write the representation 
space ``explicitly" nor an ``explicit" action for the Lie
algebra in terms of some formula
for $X(z)$. But one can start with a
Poincar\'e-Birkhoff-Witt basis of the Verma module
$M(\Lambda)$, see the corresponding spanning set in
$L(\Lambda)$, and by using the above
mentioned formula, reduce this spanning set to a spanning
set of $L(\Lambda)$ parameterized 
combinatorially by partitions satisfying the difference two
conditions. At this point
one could either recall the Rogers-Ramanujan identities (and
a specialization of the Weyl-Kac character 
formula) and conclude that this is a basis, or somehow prove
directly linear independence of 
this set and, as a consequence, prove the Rogers-Ramanujan 
identities. Later on came a third
idea for proving linear independence, somewhere in between
the first two: 
The correct character formula for $L(\Lambda)$
in terms of partitions satisfying difference two conditions
is replaced by finding a basis of 
$M^1(\Lambda)$  parameterized by partitions. It should be
noticed that in the principal picture
``group elements" $E^0(z)$  in the vertex operator formula
happen to be $\pm 1$, and hence are well 
defined operators on the Verma module. 

The Lepowsky-Wilson ideas of constructing representations by
using
vertex operator formulas work as well in the case of the
homogeneous picture, except that there the presence of
$E^0(z)$ prevents one to use Verma
modules: one has to choose among the two.

The starting point of our construction was the first
consequence of Frenkel-Kac vertex 
operator formula (in the homogeneous picture) which does not
contain terms $E^0(z)$: the relation 
$$
x_\theta(z)^{k+1}=0\tag {5}
$$
which holds on every standard $\tilde\goth g$-module
$L(\Lambda)$ of level $k$. Here $\theta$ is the
maximal root of the corresponding finite dimensional Lie algebra $\g$,
$x_\theta$ is a corresponding
root vector and  $x_\theta(z)=\sum_{n\in \Bbb
Z}\,x_\theta(n)z^{-n-1}$.  Of course, 
$x_\theta(z)^{k+1}$
is a formal Laurent series in $z$; the coefficients of
this series are sums of 
elements of the universal enveloping algebra $U(\tilde\goth g)$,
and relation (5) means that all these 
coefficients are zero on $L(\Lambda)$. 

Relation (5) exhibits several ``universal" properties: The
adjoint action of $\g$ on 
$x_\theta(z)^{k+1}$ produces new series which vanish on
$L(\Lambda)$, and these series
form an irreducible representation $R$ of $\g$ with highest
weight $(k+1)\theta$. Moreover,
all the coefficients of the series from $R$ form a loop
$\tilde\goth g$-module $\bar R$ for 
the adjoint action of $\tilde\goth g$. We
call the elements of $\bar R$ the relations (for the standard
module) because  a highest 
weight $\tilde\goth g$-module $V$ of level $k$ is standard
if and only if relation (5) holds
on $V$. Besides 
that,  $\bar R$ seems to be ``the smallest possible" loop
module that annihilates
$L(\Lambda)$, certainly ``the smallest possible" that can be
obtained from vertex operator
formula in any picture.

With regard to the Verma modules of level $k$ we have the property
that $\bar R M(\Lambda)=
M^1(\Lambda)$ if and only if $\Lambda$ is dominant integral,
otherwise $\bar R M(\Lambda)=
M(\Lambda)$. Hence in the case of the standard module
$L(\Lambda)$ we have
$$
   M^1(\Lambda)=U(\tilde\goth g)\bar R v_\Lambda,\tag {6}
$$
where $v_\Lambda$ denotes a highest weight vector in
$M(\Lambda)$. 
This result is quite similar to the one that holds in the
principal picture, and allows one to 
try to prove ``linear independence" by finding a basis in
$M^1(\Lambda)$ of the form 
$ur\cdot v_\Lambda$, $u\in U(\tilde\goth g)$, $r\in \bar R$.
Of course, as in the principal picture, one
should expect that many of these elements are linearly
dependent, and one should construct
relations among them in order to reduce this set to a basis.
This kind of relations we
call relations among relations.

By an analogy with the principal picture, for $\goth
{sl}(2,\Bbb C)\sptilde$  we could guess 
correctly a space of all relations among relations that are
needed: For level~1 modules 
it was some $\tilde\goth g$-module having a ``7-dimensional" 
loop submodule and a ``5-dimensional" loop quotient---it looked
like another ad hoc argument.
 
The elements $x_\theta(z)^{k+1}$ are zero on $L(\Lambda)$, so
they should be seen somewhere as 
nonzero  elements. The emerging theory of vertex operator algebras
provided a suitable setting where this could be achieved. It was known
that the irreducible module $L(k{\Lambda}_0)$ could be given the
structure of a vertex operator algebra in a natural way, but here the
vector $f_0^{{\Lambda}(h_0)+1}v_{\Lambda}$ that should parameterize 
$x_\theta(z)^{k+1}$ vanishes. We thus guessed and proved that a
generalized Verma module extension $N(k{\Lambda}_0)$ of $L(k{\Lambda}_0)$
also had the structure of vertex operator algebra. In $N(k{\Lambda}_0)$,
$f_0^{{\Lambda}(h_0)+1}v_{\Lambda}$ is nonzero and this provided the
required setting to explain the ``universal'' properties of (5). Later
it was also possible to prove ``all necessary'' relations among
relations by calculations in $N(k{\Lambda}_0)$: the relations among
relations arise by writing the vector
$x_{\theta}(-2)x_{\theta}(-1)^{k+1}\1$ in two
different ways, thus obtaining two different expressions for the same field.

The construction of the vertex operator algebra $N(k{\Lambda}_0)$ in
Section 3 is based on a general ``generators and relations''
construction of vertex operator algebras presented in Theorem 2.6. The
proof of Theorem 2.6 is inspired by the ideas of P\. Goddard on
locality. In our case the vertex operator algebra $V=N(k\Lambda_0)$ is
generated by a set of
fields $x(z)$, $x\in \goth g
\cong \goth g(-1)\1\subset V$. In the general construction we start
with a set of series $u(z)$, $u\in U\subset V$,
which are admissible,
local and generate $V$. By using 
normal order products we extend this generating set to a set
of mutually local fields $v(z)$
parameterized by all vectors $v\in V$, and then a general
duality argument implies that
$V$ is a vertex algebra. In Section 4 we prove that a large class of
level $k$ $\tilde\g$-modules (including all highest weight modules)
can be given the structure of modules over the vertex operator algebra 
$N(k\Lambda_0)$, a result that has been announced by B\. Feigin and
E\. Frenkel and also deduced by I\. B\. Frenkel and Y\. Zhu by
another method. Our proof uses an inductive construction of 
intertwining operators local with respect to the generating fields $x(z)$,
$x\in \goth g$, acting on a 
highest weight $\tilde\goth g$-module, and then again a
general duality argument applies. 

The results described above are the content of Sections
1--5, phrased and proved in terms
of vertex operator algebra $N(k\Lambda_0)$. In Section 8 we
construct (some) relations among 
relations in a manner mentioned above.

At this point we should say that by now there are several
proofs of the fact that $N(k\Lambda_0)$
is a vertex operator algebra and that highest weight
$\tilde\goth g$-modules are its modules. 
Moreover, several proofs about  properties of annihilating
field (5) have appeared as 
well (see remarks at the end of Sections 4 and 5).

In a sharp contrast to the general results on algebraic
properties of relations $\bar R$, we can 
successfully apply them only to the case of $\tilde\goth
g=\goth {sl}(2,\Bbb C)\sptilde$. The reader
will find Section~6 a sort of compromise between general
wishes and particular results
that we can prove, at least allowing one to speculate what
may be going on for higher
ranks. In order to describe these results, we return to our
set $\bar B$ defined above
and the corresponding set $\Cal P(\bar B)$ of colored
partitions. It is a monoid, and
we often use a multiplicative notation suggested by (1). We
construct a particular
order $\preccurlyeq$ on $\Cal P(\bar B)$ that allows
inductive arguments, behaves
well with respect to the multiplicative structure and allows
every nonzero element $r\in\bar R$
to be written in the form
$$
   r=c_\rho u(\rho)+\sum_{\pi\succ \rho}\,c_\pi u(\pi),\tag
{7}
$$
where $c_\rho, c_\pi \in \Bbb C$, $c_\rho \neq 0$. We call
$\rho$ the leading term of
$r$ and write $\rho=\lt(r)$. Moreover, the set $\lt(\bar R)$
of leading terms of all
nonzero elements in $\bar R$ parameterize a nice basis
$\{r(\rho)\mid \rho\in \lt(\bar R)\}$
of $\bar R$,
where $\rho=\lt(r(\rho))$. Since $r=0$ on $L(\Lambda)$, the
product $u(\rho)$ can be
expressed in terms of $u(\pi)$, $\pi \succ \rho$, and by
induction one easily gets that
$$
  u(\pi)\cdot v_{k\Lambda_0},\qquad \pi\in\Cal P(\bar
B_-),\quad
\pi\notin(\Cal P(\bar B_-)\cdot \lt(\bar R))\cup(\Cal P(\bar
B_-)\cdot y(0)),\tag {8}
$$
is a spanning set of $L(k\Lambda_0)$. Moreover, this is a
basis, and conditions (3)
are just explicitly spelled out conditions in (8).

One way of  proving that the set of vectors (8) is a basis
of $L(k\Lambda_0)$ would be
by invoking the correct character formula, and this is what
we do: 
The ``correct character formula" would be a basis of
$M^1(k\Lambda_0)$ parameterized by partitions
 $\Cal P(\bar B_-)\cap((\Cal P(\bar B_-)\cdot \lt(\bar
R))\cup(\Cal P(\bar B_-)\cdot y(0)))$.
So we start, as it was indicated earlier, with a spanning
set of $M^1(\Lambda)$, 
\ $\Lambda=k\Lambda_0$, of the
form
$$
    u(\pi)r(\rho)\cdot v_\Lambda, \qquad \pi \in \Cal P(\bar
B_-), \quad \rho \in \lt(\bar R).
\tag {9}
$$
Each element $u(\pi)r(\rho)$ is a sum of the form 
$$
u(\pi)r(\rho)=u(\pi\cdot\rho)+
\sum_{\kappa\succ\pi\cdot\rho}\,c_\kappa u(\kappa),\tag{10}
$$ 
and again we say that $\pi\cdot\rho$
is the leading term of $u(\pi)r(\rho)$. Almost half of this
paper is devoted to the proof of the
following general form of relations among relations: If
leading terms of $u(\pi)r(\rho)$ and
$u(\pi_1)r(\rho_1)$ are equal, then 
$u(\pi)r(\rho)-u(\pi_1)r(\rho_1)$ can be expressed
as a linear combination of terms  $u(\pi ')r(\rho')$,
$\pi'\cdot\rho'\succ\pi\cdot\rho$.
Having this, by induction we see that the spanning set (9)
of $M^1(\Lambda)$ 
can be reduced to a spanning set of vectors 
indexed by their leading terms. Since vectors
$u(\pi\cdot\rho)v_\Lambda$ are elements
of Poincar\'e-Birkhoff-Witt basis of $M(\Lambda)$, by using
(10) we see that vectors 
$u(\pi)r(\rho)v_\Lambda$ of our reduced spanning set are
linearly independent, and hence 
we have a basis of 
$M^1(k\Lambda_0)$ indexed in a desired way. Of course, the
example $\Lambda=k\Lambda_0$ was taken 
only to avoid mentioning some technical  difficulties in our
brief discussion.

As it was pointed out, relations among relations are the
central part of our proof of
linear independence. This was done in a few steps: we start
in Section 8 
with generators of the kernel of $\tilde\goth g$-module map
$$
U(\tilde \goth g) \otimes_{U(\tilde{\goth g}_{\ge 0})}
\bar{R} \bold 1
 \overset{\Psi_{N(k \Lambda_0)}} \to{\longrightarrow} N^1 (k
\Lambda_0)\rightarrow 0,
\tag{11}
$$
next we formulate relations among relations in terms of
colored partitions and leading terms
(Lemma 9.2),
and finally we prove the general case by induction and by
solving some odd dozen of 
exceptional cases (Sections 9 and 10).

If we are to single out one major point from Sections 9 and
10, it should be Lemmas 9.4 and 
10.7. They allow one to try to express
$u(\pi)r(\rho)-u(\pi_1)r(\rho_1)$ as a linear combination
of terms  $u(\pi ')r(\rho')$ first,  and to worry about
$\pi'\cdot\rho'\succ\pi\cdot\rho$ later.

Except for the infinite sums that pervade our theory, the work in
Sections 9 and 10 appears to be of the same type as that required for
producing a Groebner (or standard) basis of a suitable
$\tilde\g$-module with respect to a total order on the monomials in a
basis of $\tilde\g$ (reducing S-pairs etc.). We have not attempted to
investigate this connection in detail.

It is clear that combinatorics should enter sooner or later. Some
of the combinatorial arguments 
here are carried through with a brute force, but we still
hope that by using more of 
representation theory, say (11) for example, one can prove
relations among relations in some 
simpler way, having to use combinatorial arguments at some
later stage.

Of course, the combinatorial identity (4) is a consequence of
(3), in full generality
we discuss this in Section 11 by recalling Lepowsky-Wakimoto
product formulas for certain 
specializations of the Weyl-Kac character formula.

It is a pleasure to the second author M.P\. to thank the
Department of Mathematics of
University of Lund for their kind hospitality---a large portion
of this work was done during his visit there.

\head{1. Formal Laurent series and rational functions}\endhead

\subhead{1.1. Power series expansions}\endsubhead

We shall throughout use the calculus of formal Laurent series
as in \cite{FLM\rm, Chap. 2 and 8}. In particular $z_1,z_2,\dots$ denote
``commuting formal variables". We sometimes use the expansion convention that,
for $n<0$, $(z_0+z_1+\cdots+z_k)^n$ denotes the power series expansion in
nonnegative powers of $z_1,\dots,z_k$.

For any field $F$, denote by 
$$
F((z))=\{\sum_{n\in \Z}a_nz^n \mid a_n\in F, \exists N,
a_n=0 \hbox { unless } n\ge N\}
\tag1.1.1
$$
the field of fractions of the formal power series algebra
$F[[z]]$. Denote multiindices by
$$\alpha=(\alpha_1,\dots,\alpha_k),\qquad |\alpha|=
\alpha_1+\cdots+\alpha_k$$
and set
$$z^{\alpha}=z_1^{\alpha_1}\cdots z_k^{\alpha_k}.$$
Let $F\{z_1,\dots,z_k\}$ denote the vector space
$$
F\{z_1,\dots,z_k\} = \{\sum_{{\alpha}\in\Z^k}a_{\alpha}z^{\alpha}\mid
a_{\alpha}\in F  \}.
$$

Let $\F$ be a field of characteristic 0. We have
$$\eqalign{\F((z_1))\cdots((z_k))&=\{\sum_{\alpha}
a_{\alpha}
z^{\alpha}\mid a_{\alpha}\in \F, \exists N_k,\dots,
N_1({\alpha}_2,\dots,{\alpha}_k) \hbox{ s.t. }
 a_{\alpha}=0 
\cr&\text { unless }\quad{\alpha}_k\ge N_k,\dots,{\alpha}_1\ge
N_1({\alpha}_2,\dots,{\alpha}_k)\},\cr}$$
$$\eqalign{\F((z_k^{-1}))\cdots((z_1^{-1}))&=\{\sum_
{\alpha}a_{\alpha}
z^{\alpha}\mid a_{\alpha}\in \F, \exists N_1,\ldots\,,
N_k({\alpha}_1,\dots,{\alpha}_{k-1}) \hbox{ s.t. }
a_{\alpha}=0
\cr&\text { unless }\quad{\alpha}_1\le N_1,\dots,{\alpha}_k\le
N_k({\alpha}_1,\dots,{\alpha}_{k-1})\}.\cr}$$
There are canonical embeddings
$$\eqalign{\F[z_1,\dots,z_k]&=\F[z_1]\cdots[z_k]\to\F((z_1))
\cdots((z_k)),\cr
\F[z_1,\dots,z_k]&=\F[z_k]\cdots[z_1]\to\F((z_k^{-1}))
\cdots((z_1^{-1})).\cr}$$
Since the target algebras are fields, the embeddings
extend to embeddings
$$\eqalign{\iota:\F(z_1,\dots,z_k)&\to\F((z_1))
\cdots((z_k)),\cr
\iota':\F(z_1,\dots,z_k)&\to\F((z_k^{-1}))
\cdots((z_1^{-1})).\cr}$$

\proclaim {Lemma 1.1.1} Let  $(a_1,\dots,a_k)\in\F^k\setminus
\{(0,\dots,0)\}$. Then, in $\F\{z_1,\dots,z_k\}$,
$$\iota(\sum_{i=1}^ka_iz_i)^{-1}
=\iota'(\sum_{i=1}^ka_iz_i)^{-1}.$$
\endproclaim

\demo{Proof}
Let $a_i=0$ for $i<m$, $a_m\ne0$. Then
$$\eqalign{&\iota'(\sum_{i=1}^ka_iz_i)^{-1}
=\iota'(a_mz_m(1+\sum_{i=m+1}^k{a_iz_i\over
a_mz_m}))^{-1}\cr&\qquad=a_m^{-1}z_m^{-1}
\sum_{n\ge0}(-1)^n(\sum_{i=m+1}^k{a_iz_i\over
a_mz_m})^n\in\F((z_1))\cdots((z_k)).\cr}$$
Hence
$$\eqalign{&\iota'(\sum_{i=1}^ka_iz_i)^{-1}
=\iota'(\sum_{i=1}^ka_iz_i)^{-1}(\sum_{i=1}^ka_iz_i
\iota(\sum_{i=1}^ka_iz_i)^{-1})\cr
&\qquad=(\iota'(\sum_{i=1}^ka_iz_i)^{-1}\sum_{i=1}^ka_iz_i)
\iota(\sum_{i=1}^ka_iz_i)^{-1}
=\iota(\sum_{i=1}^ka_iz_i)^{-1}\cr}$$
by the associative law in $\F((z_1))\cdots((z_k))$.
\qed
\enddemo

Set 
$$
\align
S&=\{\sum_{i=1}^ka_iz_i\mid (a_1,\dots,a_k)\ne(0,\dots,0)\},
\tag1.1.2 \\
\Rhat_k&=\F[z_1,\dots,z_k]_{\langle S\rangle},
\tag1.1.3
\endalign
$$
the localization.

\proclaim {Corollary 1.1.2} The homomorphisms \quad $\iota$,
$\iota'$ coincide on $\Rhat_k$.
\endproclaim

\remark{Remark} We have
$$\iota(1-z_1)^{-1}=\sum_{n\ge0}z_1^n
\ne\iota'(1-z_1)^{-1}=-\sum_{n\ge0}z_1^{-n-1}$$
although $\iota(z_1-z_2z_3)^{-1}=\iota'(z_1-z_2z_3)^{-1}$.
It seems difficult to determine for which
\hfil \break 
$f\in\F(z_1,\dots,z_k)$, $\iota(f)=\iota'(f)$.
\endremark

Define 
$$
\iota_{1,\dots,k}:\Rhat_k\to\F((z_1))\cdots((z_k))
\tag1.1.4
$$
to be $\iota_{1,\dots,k}=\iota=\iota'$.
Also define the partial $\iota$-maps
$$
\iota_{\cdots\ell+1,\dots,k}:\Rhat_k\to\Rhat_{\ell}
((z_{\ell+1}))\cdots((z_k));
\tag1.1.5
$$
the restriction of
$$\F(z_1,\dots,z_k)\to\F(z_1,\dots,z_\ell)((z_{\ell+1}))
\cdots((z_k)).$$
By considering the image of $(\sum_{i=1}^ka_iz_i)^{-1}$ 
one verifies that the image of $\Rhat_k$ actually
lands in  
$$\Rhat_{\ell}((z_{\ell+1}))\cdots((z_k)).$$
For $m,n\in\N$, $m\le n$, denote by $\Rhat_{[m,n]}$ the
algebra analogous to $\Rhat_k$, based on the variables $z_i$,
$m\le i\le n$. Define also
$$
\iota_{1,\dots,\ell\cdots}:\Rhat_k\to
\Rhat_{[\ell+1,k]}((z_{\ell}^{-1}))\cdots((z_1^{-1})),
\tag1.1.6
$$
the restriction of
$$\F(z_1,\dots,z_k)\to\F(z_{\ell+1},\dots,z_k)
((z_{\ell}^{-1}))\cdots((z_1^{-1}));$$
here it is easier to see that the image is in $\Rhat_{[\ell+1,k]}\cdots\,$.

\proclaim {Proposition 1.1.3} We have
$$\iota_{1,\dots,k}=\iota_{1,\dots,\ell}\circ
\iota_{\cdots\ell+1,\dots,k}=\iota_{\ell+1,\dots,k}\circ
\iota_{1,\dots,\ell\cdots}:\Rhat_k\to
\F((z_1))\cdots((z_k)).$$
\endproclaim

\demo{Proof} The composition $\iota_{1,\dots,\ell}\circ
\iota_{\cdots\ell+1,\dots,k}$ is the restriction to $\Rhat_k$
of the canonical map
$$\iota:\F(z_1,\dots,z_k)\to\F((z_1))\cdots((z_k))$$
i.e. $\iota_{1,\dots,k}$. Also $\iota_{\ell+1,\dots,k}=
\iota=\iota'_{[\ell+1,k]}$ on $\Rhat_{[\ell+1,k]}$ so that
the composition $\iota_{\ell+1,\dots,k}\circ
\iota_{1,\dots,\ell\cdots}$ is the restriction to $\Rhat_k$
of the canonical map 
$$\iota':\F(z_1,\dots,z_k)\to
\F((z_k^{-1}))\cdots((z_1^{-1})).$$
By the corollary $\iota=\iota'$ on $\Rhat_k$.
\qed \enddemo

\remark{Remark} It follows from the formula
$\iota_{1,\dots,k}=
\iota_{2,\dots,k}\circ\iota_{1,\cdots}$ that the present
definition of $\iota_{1,\dots,k}$ coincides with the one
in  \cite{FLM\rm, (8.10.37--42)}.
\endremark

Let $\sigma \in S_k=\text{Sym}\{1,\dots,k\}$. Denote also by
$\sigma$ the automorphism
$$\eqalign{\sigma:\F[z_1,\dots,z_k]&\to\F[z_1,\dots,z_k]\cr
z_m\quad &\mapsto \quad z_{\sigma m}.\cr}$$
Then $\sigma (S)=S$, hence $\sigma$ extends to an
automorphism
$$\sigma:\Rhat_k\to\Rhat_k.$$
Let $\sigma$ also denote the vector space isomorphism
$$\eqalign{\sigma:\F\{z_1,\dots,z_k\}&\to
\F\{z_1,\dots,z_k\}\cr
\sum_{\alpha}a_{\alpha}z^{\alpha}&\mapsto
\sum_{\alpha}a_{\alpha}z_{\sigma1}^{\alpha_1}\cdots
z_{\sigma k}^{\alpha_k}=\sum_{\alpha}a_{\alpha}
z^{\alpha\circ\sigma^{-1}}.\cr}$$
Then $\sigma$ has the homomorphism property; if
$$x,y\in\F\{z_1,\dots,z_k\}$$
and $xy$ is defined, then $\sigma(x)\sigma(y)$ is also
defined and
$$\sigma(xy)=\sigma(x)\sigma(y).$$

Define also  
$$
\align
\iota_{\sigma1,\dots,\sigma k}&=\sigma\circ
\iota_{1,\dots,k}\circ\sigma^{-1},
\tag1.1.7 \\
\iota_{\cdots\sigma(\ell+1),\dots,\sigma(k)}&=
\sigma\circ\iota_{\cdots \ell+1,\dots,
k}\circ\sigma^{-1},
\tag1.1.8 \\
\iota_{\sigma(1),\dots,\sigma(\ell)\cdots}&=
\sigma\circ\iota_{1,\dots,\ell\cdots}\circ\sigma^{-1}.
\tag1.1.9
\endalign
$$
Then we have more generally
$$\eqalign{&\iota_{\sigma1,\dots,\sigma k}=
\iota_{\sigma1,\dots,\sigma \ell}\circ
\iota_{\cdots\sigma(\ell+1),\dots,\sigma(k)}
=\iota_{ \sigma(\ell+1),\dots,\sigma(k)}\circ
\iota_{\sigma(1),\dots,\sigma(\ell)\cdots}\cr
&\qquad:\Rhat_k\to\F((z_{\sigma1}))\cdots((z_{\sigma
k})).\cr}\tag1.1.10$$

\subhead{1.2. Cauchy's theorem}\endsubhead

Set
$$
{\delta}(z) = \sum_{n\in\Z}z^n.
\tag1.2.1
$$
For $k\in\N$ set
$$\align
\R_k &= \F[z_i,z_i^{-1},(z_m-z_n)^{-1}; 1\le i\le k, 1\le m<n\le k],
\tag1.2.2 \\
{\Cal L}_k &= \F[z_i,z_i^{-1}; 1\le i\le k].
\tag1.2.3
\endalign
$$

The following Proposition was first formulated by H.~Garland, cf.\ \cite{FLM}.

\proclaim{Proposition 1.2.1} We have
$$
{\delta}(z)p(z) = {\delta}(z)p(1),
$$
$$
{\delta}(az_1/z_2)q(z_1,z_2) = {\delta}(az_1/z_2)q(z_1,az_1)
$$
for all $ p\in\F[z,z^{-1}]$, $ a\in\F$, and $q\in{\Cal L}_2$.
\endproclaim

\proclaim {Proposition 1.2.2} Let $g(z_0, z_1, z_2) \in {\Cal L}_3$. Then
$$\Res_{z_0}\,z_0^{-1}\delta({z_1-z_2\over z_0})g(z_0, z_1, z_2) = \iota_{12}
g(z_1-z_2, z_1, z_2).$$
\endproclaim

\proclaim {Proposition 1.2.3} Let $f( z_1, z_2) \in {\Cal R}_2$. Then
$$\Res_{z_1, z_2}\,\iota_{12}f - \Res_{z_1, z_2}\,\iota_{21}f
= \Res_{z_0, z_2}\,\iota_{20}f(z_2+z_0, z_2).$$
\endproclaim

\demo
{Proof} By  partial fraction decomposition one may assume that
$$f\in \{z_1^mz_2^n, z_2^m(z_1-z_2)^n\mid m, n \in \Z\}.$$
If $m+n \ne -2$ we have $0-0=0$. If $f = {\partial\over \partial z_1}g$ then
all terms vanish as well. Let $f = z_1^{-1}z_2^{-1}
$. Then
$$\Res_{z_1, z_2}\,\iota_{12}f - \Res_{z_1, z_2}\,\iota_{21}f
= 0,$$
$$\Res_{z_0, z_2}\,\iota_{20}f(z_2+z_0, z_2)=\Res_{z_0, z_2}\,\iota_{20}
(z_2+z_0)^{-1}z_2^{-1} = 0.$$
Let $f = z_2^{-1}(z_1-z_2)^{-1}$. Then
$$\Res_{z_1, z_2}\,\iota_{12}f= \Res_{z_1, z_2}\,z_2^{-1}\sum_{i\ge 0 }
{-1 \choose i}z_1^{-1-i}(-z_2)^i = 1,$$
$$\Res_{z_1, z_2}\,\iota_{21}f= \Res_{z_1, z_2}\,z_2^{-1}\sum_{i\ge 0 }
{-1 \choose i}(-z_2)^{-1-i}z_1^i = 0,$$
$$\Res_{z_0, z_2}\,\iota_{20}f(z_2+z_0, z_2)= \Res_{z_0, z_2}\,z_2^{-1}z_0^{-1}
= 1. \quad \qed $$
\enddemo

\remark{Remark} Proposition 1.2.3 follows alternatively from  Cauchy's theorem:
Take real numbers $\epsilon$, $r$, $\rho$, $R$, such that $0 < r < \rho < R$,
$0 < \epsilon < \min\{R-\rho, \rho-r\}$. Then (if $f$ has
coefficients in $\C$)
$$
\leqalignno{\Res_{z_1, z_2}\,&\iota_{12}f - \Res_{z_1, z_2}\,\iota_{21}f\cr
&={1\over (2\pi i)^2}\int_{C_R}(\int_{C_\rho}f(z_1, z_2)dz_2)dz_1
-{1\over (2\pi i)^2}\int_{C_\rho}(\int_{C_r}f(z_1, z_2)dz_1)dz_2\cr
&={1\over (2\pi i)^2}\int_{C_\rho}(\int_{C_\epsilon(z_2)}f(z_1, z_2)dz_1)
dz_2 &(1.2.4)\cr
&= {1\over (2\pi i)^2}\int_{C_\rho}(\int_{C_\epsilon}f(z_2+z_0, z_2)dz_0)
dz_2\cr
&= \Res_{z_0, z_2}\,\iota_{20}f(z_2+z_0, z_2).\cr}
$$
\endremark

An equivalent formulation of this result is:
\proclaim {Proposition 1.2.4} Let $f( z_1, z_2) \in {\Cal R}_2$. Then
$$
\Res_{z_1, z_2}\,\iota_{12}f = \Res_{z_1, z_2}\,(\iota_{21}-\iota_{12})
f(z_1, z_1+z_2).
$$
\endproclaim

\proclaim{Lemma 1.2.5} For all $f(z_1,z_2)\in\R_2$,
$$
T{\iota}_{34}f(z_4,z_3) = \Res_{z_4}{\iota}_{34}z_2^{-1}{\delta}(\frac{z_3-z_4}
{-z_2})f(z_4,z_3) = ({\iota}_{23}-{\iota}_{32})f(z_2+z_3,z_3).
$$
\endproclaim

\demo
{Proof} It is enough to check the coefficients of
$z_2^{-m-1}z_3^{-n-1}$ on both sides. Let $g(z_2,z_3) =
z_2^mz_3^n\in{\Cal L}_2$. Then
$$\split
\Res_{z_2,z_3}&\Res_{z_4}{\iota}_{34}z_2^{-1}{\delta}(\frac{z_3-z_4}
{-z_2})f(z_4,z_3)g(z_2,z_3)\\
&= \Res_{z_3,z_4}{\iota}_{34}f(z_4,z_3)g(-z_3+z_4,z_3)\\
&= \Res_{z_3,z_4}({\iota}_{43}-{\iota}_{34})f(z_3+z_4,z_3)g(z_4,z_3)\\
&= \Res_{z_2,z_3}({\iota}_{23}-{\iota}_{32})f(z_2+z_3,z_3)g(z_2,z_3).
\quad \qed
\endsplit
$$
\enddemo

\head{2. Generating fields}\endhead

\definition{Definition 2.1} \cite{B}, \cite{FLM}, \cite{FHL} A {\it vertex operator algebra}
$(V,Y,\1,{\omega})$ consists of a vector space
$V$ over $\F$, distinguished vectors $\1,{\omega}\in V$, and a linear map
$$
\align
Y:V&\to (\End_\F V)[[z,z^{-1}]],\\
v&\mapsto Y(v,z)=\sum_{n\in\Z}v_nz^{-n-1}
\endalign
$$
satisfying the following axioms:
$$
\align
\text{For all } &u,v\in V\text{ there is an } N=N(u,v)\in\Z \text{ such that }
\tag2.1\\
&u_nv=0 \text{ for all } n > N,
\endalign
$$
and
$$
v_{-1}\1 = v,\qquad v_n\1=0
\tag2.2
$$
for $v\in V, n\ge 0$.
Define $L_n$ by
$$
Y({\omega},z)=\sum_{n\in\Z}L_nz^{-n-2}
\tag2.3
$$
(i\.e\. $L_n={\omega}_{n+1}$). Then
$$
[L_{-1}, Y(v,z)] = \frac d{dz}Y(v,z) = Y(L_{-1}v,z)
\tag2.4
$$
for all $v\in V$, and
$$
[L_m,L_n] = (m-n)L_{m+n}+\frac{m^3-m}{12}{\delta}_{m+n,0}(\rank V)1_V
\tag2.5
$$
for some scalar $\rank V\in\F$. Set, for $n\in\Z$,
$$
V_n = \{v\in V\mid L_0v=nv\}.
\tag2.6
$$
Then
$$
\align
&V=\coprod_{n\in\Z} V_n,
\tag2.7 \\
&\dim_\F V_n < \infty \text{ for all $n$,}
\tag2.8 \\
&V_n=0 \text { for $n$ sufficiently small}.
\tag2.9
\endalign
$$
Finally,
$$
\split
z_0^{-1}&{\delta}(\frac{z_1-z_2}{z_0})Y(u,z_1)Y(v,z_2)-
z_0^{-1}{\delta}(\frac{z_2-z_1}{-z_0})Y(v,z_2)Y(u,z_1)\\
&=z_2^{-1}{\delta}(\frac{z_1-z_0}{z_2})Y(Y(u,z_0)v,z_2)
\endsplit
\tag2.10
$$
for all $u,v\in V$ (the Jacobi identity).
\enddefinition

We say that $Y(v,z)$ is a {\it vertex operator} or a {\it
field} associated with a vector 
$v$. Later on we shall consider a formal Laurent series
$Y(v,z)$ as a generating function
for the coefficients $v_n$, $n\in \Bbb Z$.

In later sections we shall use the following consequences of the axioms
\cite{B}, \cite{FLM}, \cite{FHL}: For all $u,v\in V, n\in\Z$,
$$
\split
Y(u_nv,z)&=\Res_{z_0}z_0^n\big(Y(u,z_0+z)Y(v,z) \\
&-Y(v,z)\Res_{z_1}z_0^{-1}{\delta}(\frac{z-z_1}{-z_0})Y(u,z_1)\big),
\endsplit
\tag2.11
$$
with the notation
$$
Y(v,z)^-=\sum_{n<0}v_nz^{-n-1},\qquad Y(v,z)^+=\sum_{n\ge 0}v_nz^{-n-1},
\tag2.12
$$
we have
$$
\split
Y(u_nv,z)&=
\Big(\frac{(\partial/\partial z)^{-n-1}}{(-n-1)!}Y(u,z)^-\Big)Y(v,z)\\
&+Y(v,z)\Big(\frac{(\partial/\partial z)^{-n-1}}{(-n-1)!}Y(u,z)^+\Big)
\endsplit
\tag2.13
$$
for $n\le -1$, in particular
$$
Y(u_{-1}v,z)=Y(u,z)^-Y(v,z)+Y(v,z)Y(u,z)^+,
\tag2.14
$$
furthermore
$$
[Y(u,z_1),Y(v,z_2)] =
\sum_{i=0}^\infty\frac{(-1)^i}{i!}\big(\frac\partial{\partial z_1}\big)^iz_2^{-1}{\delta}(z_1/z_2)Y(u_iv,z_2),
\tag2.15
$$
and
$$
Y(u,z)v = e^{zL_{-1}}Y(v,-z)u.
\tag2.16
$$

The goal of this section is to derive Theorem 2.6 below giving a kind
of ``generators and relations" construction of vertex operator algebras.

Let $V$ be a vector space over $\F$, and let $L_{-1}, L_0 \in
\End_\F V$, $\1\in V$. Set, for $n\in \Z$,
$$
V_n = \{v\in V\mid L_0v=nv\}.
$$
Assume that
$$\align
V&=\coprod_{n\in\Z}V_n,\tag2.17\\
\1&\in V_0,\quad \text{ i.e. }L_0\1 = 0,\tag2.18\\
[L_0,L_{-1}] &= L_{-1},\quad \text{ i.e. }L_{-1}V_n\subset V_{n+1}.\tag2.19
\endalign
$$
\definition{Definition 2.2} A generating function
$$
a(z) = \sum_{n\in\Z}a_nz^{-n-1}
$$
of operators $a_n\in\End_\F V$ is said to  be {\it admissible of conformal
weight} $h$ if
$$\align
\text{ for all }v\in V &\text{ there is } N(a,v)\in\Z \text{ s.t. }a_nv=0
\text{ if }n > N(a,v),\tag2.20\\
[L_{-1},a(z)]&=\frac{d}{dz}a(z),\tag2.21\\
[L_0,a(z)]&=z\frac{d}{dz}a(z)+ha(z),\tag2.22\\
a(z)\1&\in V[[z]],\quad \text{ i.e. }a_n\1=0 \text{ if } n\ge0.\tag2.23
\endalign
$$
A generating function is {\it admissible} if it is a (finite) linear
combination of admissible generating functions of operators with (various)
conformal weights.
\enddefinition

It follows from (2.18) and (2.22) that
$$
a_{-1}\1\in V_h,\quad \text{ i.e. } L_0a_{-1}\1=ha_{-1}\1.\tag2.24
$$

\proclaim{Proposition 2.3} 
Let $a(z),b(z)$ be admissible generating functions. Define
$d_n(z)$ for $n\in\Z$ by
$$
\sum_{n\in\Z}d_n(z_3)z_2^{-n-1}= a(z_2+z_3)b(z_3)-b(z_3)Ta(z_4),
\tag2.25
$$
where $T$ denotes the operator
$$
T = \Res_{z_4}z_2^{-1}\delta(\frac{z_3-z_4}{-z_2})(\cdot)
\tag2.26
$$
(cf. \cite{FLM\rm, (8.8.31)}).
Then $d_n(z)$ is admissible for all $n\in\Z$.
\endproclaim
\demo{Proof}
(2.20): Let $v\in V$. Then
$$\split
(d_n)_kv &= \coeff_{z_2^{-n-1}z_3^{-k-1}}\{a(z_2+z_3)b(z_3)-b(z_3)
[a(z_2+z_3)-a(z_3+z_2)]\}v \\
&= \sum_{r\le n}\binom{-r-1}{n-r}a_rb_{k+n-r}v
- \sum_{0\le r\le n}\binom{-r-1}{n-r}b_{k+n-r}a_rv \\
&+ \sum_{r\ge 0}\binom{-r-1}{-n-1}b_{k+n-r}a_rv.
\endsplit
\tag2.27
$$
There are $N_1,N_2\in\Z$ such that $N_1\ge 0$, $a_mv = 0$ if $m\ge
N_1$, $b_mv = 0$ and $b_ma_rv = 0$ if $m\ge N_2$ and $0\le r \le N_1-1$.
If $N = \max\{-n,0\}+N_1+N_2$ we then have $(d_n)_kv = 0$ if $k\ge N$.

(2.21): It follows from
$$
Ta(z_4) = \sum_{j\in\Z}a_j\{(z_2+z_3)^{-j-1}-(z_3+z_2)^{-j-1} \}
$$
that
$$
[L_{-1}, Ta(z_4)] = \frac{\partial}{\partial z_3}Ta(z_4).
$$
Hence
$$\split
[L_{-1}&, \sum_{n\in\Z}d_n(z_3)z_2^{-n-1}] 
= [L_{-1},a(z_2+z_3)]b(z_3)+a(z_2+z_3)[L_{-1},b(z_3)] \\
&-[L_{-1},b(z_3)]Ta(z_4)
-b(z_3)[L_{-1},Ta(z_4)] \\
&= \frac{\partial}{\partial z_3}(\sum_{n\in\Z}d_n(z_3)z_2^{-n-1}).
\endsplit
$$

(2.22): We may assume that $a$ and $b$ have conformal weights $h_a$,
$h_b$. Then
$$\split
[L_0&, \sum_{n\in\Z}d_n(z_3)z_2^{-n-1}] 
= [L_0,a(z_2+z_3)]b(z_3)+a(z_2+z_3)[L_0,b(z_3)] \\
&-[L_0,b(z_3)]Ta(z_4)
-b(z_3)[L_0,Ta(z_4)] \\
&= (z_2\frac{\partial}{\partial z_2}+z_3\frac{\partial}{\partial z_3}+h_a+h_b)
(\sum_{n\in\Z}d_n(z_3)z_2^{-n-1}) \\
&= \sum_{n\in\Z}\{
(z_3\frac{\partial}{\partial z_3}+h_a+h_b-n-1)
d_n(z_3) \}z_2^{-n-1}
\endsplit
$$
so that $d_n$ has conformal weight $h_a+h_b-n-1$.

(2.23): This follows from (2.27) and the fact that $a_k\1=0=b_k\1$ if
$k\ge 0$.
\qed \enddemo

Set
$$
V'=\coprod_{n\in \Z}V_n^*,
$$
where $V_n^* = \Hom_\F(V_n,\F)$ and denote the pairing between $V'$
and $V$ by
$$
\langle\cdot,\cdot\rangle:V'\times V\to \F.
$$

\definition{Definition 2.4} Let $a(z), b(z)$ be admissible generating
functions. We say that $a(z), b(z)$ are mutually {\it local} or {\it
commutative} if there exists $N\in \N$ such that for all $w\in V,
w'\in V'$,
$$
(z_1-z_2)^N\langle w', a(z_1)b(z_2)w\rangle = (z_1-z_2)^N\langle
w', b(z_2)a(z_1)w\rangle
$$
i.e. $(z_1-z_2)^N[a(z_1),b(z_2)] = 0$.
\enddefinition

Equivalently, $a(z)$ and $b(z)$ are mutually local iff there is $f(z_1,z_2)\in\R_2$
such that
$$\align
\la w', a(z_1)b(z_2)w\ra &= {\iota}_{12}f,\\
\langle w', b(z_2)a(z_1)w\ra &= {\iota}_{21}f,
\endalign
$$
with $\iota$ as in Section 1.1, 
and  there is an upper bound that is independent of $w$ and $w'$ for the
order of the pole $(z_1-z_2)^{-1}$ in such matrix
coefficients. 

\proclaim{Proposition 2.5}
Let $a(z),b(z), c(z)$ be pairwise local admissible generating
functions. Let 
$d_n(z)$ for $n\in\Z$ be defined as in (2.25), (2.26).
Then $d_n(z)$ and $c(z)$ are mutually local for
each $n\in \Z$.
\endproclaim
\demo{Proof}
Set
$$
d(z_2,z_3) = a(z_2+z_3)b(z_3)-b(z_3)Ta(z_4)\in\End_\F V\{z_2,z_3\}.
$$
Let $w\in V$, $w'\in V'$. We shall calculate the matrix coefficients
$\la w',c(z_1)d(z_2,z_3)w\ra $ and $\la w',d(z_2,z_3)c(z_1)w\ra $.
It follows from (2.17), (2.20), (2.22) that there is
$$
f(z_1,z_2,z_3)\in\R_3
$$
such that
$$
\align
\la w',c(z_1)a(z_2)b(z_3)w\ra &= {\iota}_{123}f,\\
\la w',c(z_1)b(z_3)a(z_2)w\ra &= {\iota}_{132}f.
\endalign
$$
Then
$$
\la w',c(z_1)a(z_2+z_3)b(z_3)w\ra = {\iota}_{123}f(z_1,z_2+z_3,z_3).
$$
Let ${\iota}_{1\cdots}f = \sum_{n\in\Z}g_n(z_2,z_3)z_1^n $.
Using Lemma 1.2.5 and Proposition 1.1.3 we obtain
$$
\split
\la w'&,c(z_1)b(z_3)Ta(z_4)w\ra = T{\iota}_{134}f(z_1,z_4,z_3)\\
&= \sum_{n\in\Z}T{\iota}_{34}g_n(z_4,z_3)z_1^n  \\
&= \sum_{n\in\Z}({\iota}_{23}-{\iota}_{32})g_n(z_2+z_3,z_3)z_1^n \\
&= ({\iota}_{123}-{\iota}_{132})f(z_1,z_2+z_3,z_3).
\endsplit
$$
Hence
$$
\la w',c(z_1)d(z_2,z_3)w\ra = {\iota}_{132}f(z_1,z_2+z_3,z_3).
\tag2.28
$$
Similarly,
$$
\la w',a(z_2+z_3)b(z_3)c(z_1)w\ra = {\iota}_{231}f(z_1,z_2+z_3,z_3),
$$
and
$$
\split
\la w'&,b(z_3)Ta(z_4)c(z_1)w\ra = T{\iota}_{341}f(z_1,z_4,z_3)\\
&= \sum_{n\in\Z}T{\iota}_{34}h_n(z_4,z_3)z_1^n  \\
&= \sum_{n\in\Z}({\iota}_{23}-{\iota}_{32})h_n(z_2+z_3,z_3)z_1^n \\
&=({\iota}_{231}-{\iota}_{321})f(z_1,z_2+z_3,z_3)
\endsplit
$$
if ${\iota}_{\cdots 1}f = \sum_{n\in\Z}h_n(z_2,z_3)z_1^n $. This gives
$$
\la w',d(z_2,z_3)c(z_1)w\ra = {\iota}_{321}f(z_1,z_2+z_3,z_3).
$$
Since 
$$
f(z_1,z_2+z_3,z_3)\in\F[z_1,z_2,z_3,z_1^{-1},z_2^{-1},z_3^{-1},
(z_1-z_2-z_3)^{-1},(z_1-z_3)^{-1},(z_2+z_3)^{-1}],
$$
we have ${\iota}_{321}f(z_1,z_2+z_3,z_3) = {\iota}_{312}f(z_1,z_2+z_3,z_3) $
so that
$$
\la w',d(z_2,z_3)c(z_1)w\ra = {\iota}_{312}f(z_1,z_2+z_3,z_3).
\tag2.29
$$
Define $f_n(z_1,z_3)\in\R_{\{ 1,3 \}} $ by
$$
{\iota}_{\cdots 2}f(z_1,z_2+z_3,z_3) = \sum_{n\in\Z}f_n(z_1,z_3)z_2^{-n-1}.
$$
Then (2.28) and (2.29) imply
$$\align
\la w',c(z_1)d_n(z_3)w\ra &= {\iota}_{13}f_n,\\
\la w',d_n(z_3)c(z_1)w\ra &= {\iota}_{31}f_n.
\endalign
$$
Let
$$
f(z_1,z_2,z_3) = \frac{p(z_1,z_2,z_3)}{(z_1-z_2)^K(z_1-z_3)^L(z_2-z_3)^M}
$$
with $p\in{\Cal L}_3 $.
Then
$$
\split
{\iota}_{\cdots 2}&f(z_1,z_2+z_3,z_3)
= {\iota}_{\cdots 2}\frac{p(z_1,z_2+z_3,z_3)}{(z_1-z_2-z_3)^K(z_1-z_3)^L
z_2^M}  \\
&= \Bigl(\sum_{j\ge 0}\binom {-K}j(z_1-z_3)^{-K-j}(-z_2)^j\Bigr)
(z_1-z_3)^{-L}{\iota}_{\cdots 2}\frac{p(z_1,z_2+z_3,z_3)}{z_2^M}.
\endsplit
$$
Since ${\iota}_{\cdots 2}p(z_1,z_2+z_3,z_3) $ only involves nonnegative
powers of $z_2$, the $j$-th term
contributes to $f_n$ only if $j-M\le -n-1 $. Hence the order of the pole
$(z_1-z_3)^{-1} $ in
$f_n$ is at most $K+L+M-n-1 $, and since $K,L,M $ can be taken to be
independent of $w,w'$, we are done.
\qed \enddemo

We shall now assume that
$$
L_{-1}\1 = 0.\tag2.30
$$
The proof of the following theorem uses the ideas of P\. Goddard in
\cite{G}. A similar result appears in \cite{Xu, Theorem 2.4}.

\proclaim{Theorem 2.6} Let $U$ be a subspace of $V$. Assume that there
is given a linear map
$$\align
Y:U&\to (\End_\F V)[[z,z^{-1}]] \\
u&\mapsto Y(u,z)=\sum_{n\in\Z}u_nz^{-n-1}
\endalign
$$
such that
$$\align
&Y(u,z) \text{ is admissible for all } u\in U, \tag2.31 \\
&u_{-1}\1 = u \text{ for all } u\in U, \tag2.32 \\
&Y(u,z) \text{ and } Y(v,z) \text{ are local for each pair } u,v\in U,
\tag2.33 \\
&V = \F\text{-span}\{u_{n_1}^{(1)}\cdots u_{n_k}^{(k)}\1\mid k\in\N,
u^{(i)}\in U, n_i\in\Z\}. \tag2.34
\endalign
$$
Then $Y$ extends uniquely to a linear map
$$
Y: V\to (\End_\F V)[[z,z^{-1}]]
$$
that makes $V$ into a vertex operator algebra (except that we have
only $L_0$, $L_{-1}$ out of the Virasoro algebra).
\endproclaim
\demo{Proof}
Let $W$ be the space of all admissible generating functions
$$a(z)\in (\End_\F
V)[[z,z^{-1}]]
$$
such that
$$
a(z) \text{ and } Y(u,z) \text{ are local for each } u\in U.
\tag2.35
$$
Define
$$\align
{\phi}: W&\to V \\
a(z)&\mapsto a_{-1}\1.
\endalign
$$

Step 1. (uniqueness) ${\phi}$ is injective.
\demo{Proof}
Assume that $a(z)\in W$, $a_{-1}\1 = 0$. It follows from (2.21) and (2.30) that
$$
a(z)\1 = e^{zL_{-1}}a_{-1}\1 = 0.\tag2.36
$$
Let $X = \{v\in V\mid a(z)v = 0\}$ and let $v\in X$, $u\in U$. Then there
exist $N\in \N$ such that
$$
(z_1-z_2)^N\la w', a(z_1)Y(u,z_2)v\ra = (z_1-z_2)^N\la w', Y(u,z_2)a(z_1)v\ra
 = 0
$$
for all $w'\in V' $, hence $a(z_1)u_nv = 0 $ and $u_nX\subset X $. Since by
(2.36) $\1\in X $ it follows from (2.34)
that $X = V $ so that $a(z) = 0 $ and ${\phi}$ is injective.
\enddemo

By (2.31), (2.33), $Y(u,z)\in W$ for all $u\in U$, hence ${\phi}(W)\supset U $
by (2.32). Set
$$
X = {\phi}(W),\qquad U\subset X\subset V.
$$
Since ${\phi}$ is injective we can define
$$
Y: X\to(\End_\F V)[[z,z^{-1}]]
$$
by
$$
Y(v,z) = {\phi}^{-1}(v)
$$
for $v\in X $. By (2.32) the two meanings of $Y(u,z)$ for $u\in U$ denote the
same series.

Step 2. $X = {\phi}(W) = V$.
\demo{Proof}
For $u\in U $, $v\in X $, $n\in\Z $ define $Y(u_nv,z) $ by
$$
\sum_{n\in\Z}Y(u_nv,z_3)z_2^{-n-1} = Y(u,z_2+z_3)Y(v,z_3)-
Y(v,z_3)TY(u,z_4).
$$
Then $Y(u_nv,z)\in W $ by Propositions 2.3 and 2.5, and we have
$$
{\phi}Y(u_nv,z)=(u_nv)_{-1}\1 = u_nv_{-1}\1 = u_nv
$$
by (2.27).
Hence $u\in U $, $v\in X $, $n\in \Z $ implies
$$
u_nv\in X.
\tag2.37
$$
Since $\text{id}_Vz^0\in W $ and $\coeff_{z^0}(\text{id}_Vz^0)\1 = \1 $ we have
$$
\1\in X,
\tag2.38
$$
and
$$
Y(\1,z) = \text{id}_Vz^0.
$$
Now (2.34), (2.37), (2.38) imply $X = V$, i.e. we have well defined vertex
operators $Y(v,z) $ for all $v\in V $.
\enddemo

Step 3. $Y(u,z)$ and $Y(v,z) $ are mutually local for all pairs
$u,v\in V $.
\demo{Proof}
Fix $v^{(0)}\in V $ and set
$$X =\{v\in V\mid Y(v,z) \text{ and }Y(v^{(0)},z)
\text{ are mutually local}\}.
$$
Clearly $\{\1\}\cup U\subset X$. If
$u\in U $, $v\in X $, $n\in\Z $ then
$$
Y(u_nv,z_3) = \coeff_{z_2^{-n-1}}(
Y(u,z_2+z_3)Y(v,z_3)-Y(v,z_3)TY(u,z_4))
$$
together with Proposition 2.5 show that  $u_nv\in X $. Hence (2.34) again
implies
$X = V$ as required.
\enddemo

We have showed that $V$ satisfies the conditions of the axiomatic
characterization of vertex operator algebras \cite{FLM,
A.2.5, A.2.8, A.3.1; FHL} (except that we only have $L_0$, $L_{-1}$ out
of the Virasoro algebra).
\qed \enddemo

\head{3. The vertex operator algebra $N(k\Lambda_0)$}\endhead

Let $\g$ be a finite dimensional split simple Lie algebra
over $\F$, and let
$(\cdot,\cdot):\g\times\g\to\F$ be the nondegenerate
invariant
symmetric bilinear form normalized so that
$({\alpha},{\alpha}) = 2$
for each long root ${\alpha}$. We then have the associated
affine
Kac-Moody algebra
$$
\align
&\hat\g = \g\otimes\F [t,t^{-1}] \oplus \F c,
\qquad \tilde\goth g=\ghat\oplus\Bbb F L_0,\tag3.1 \\
&c \text{ central },\quad c\ne 0, \qquad [L_0,x\otimes
t^m]=-m\,x\otimes t^m,\\
&[x\otimes t^m, y\otimes t^n] = [x,y]\otimes t^{m+n} +
m{\delta}_{m+n,0}
(x,y) c.
\endalign
$$
Note that $d=-L_0$ defines the usual homogeneous grading on
$\ghat$.
We shall frequently denote $x\otimes t^m$ by $x(m)$ and
$\g\otimes t^m$ by $\g(m)$. We often identify $\g(0)$ and $\g$. Set
$$
\frak{p} = \coprod_{n\ge 0}\g\otimes t^n \oplus\F c,
\qquad \tilde{\goth g}_{\ge 0}=\goth p\oplus \Bbb F L_0.
\tag3.2
$$
For $k\in\F$ denote by $\F v_{k\Lambda_0}$ the
1-dimensional $\tilde{\goth g}_{\ge 0}$-module
such that
$$
c\cdot v_{k\Lambda_0} = kv_{k\Lambda_0},\qquad L_0\cdot
v_{k\Lambda_0} =0,
\qquad (\g\otimes t^n)\cdot v_{k\Lambda_0} = 0
\tag3.3
$$
for $n\ge 0$. Form the induced $\tilde\g$-module (a
generalized Verma module)
$$
N(k\Lambda_0) = U(\hat\g)\otimes_{U(\frak{p})}\F
v_{k\Lambda_0}=
U(\tilde\g)\otimes_{U(\tilde{\goth g}_{\ge 0})}\Bbb F
v_{k\Lambda_0}.
\tag3.4
$$
Then $N(k\Lambda_0)\cong U(\tilde\g_{<0})$ as vector spaces, where
$$
\tilde\g_{<0} = \coprod_{n < 0} \g\otimes t^n.
$$
We denote by $x(n)$ the operator on $N(k\Lambda_0)$
corresponding to $x\otimes t^n$.
Set
$$
\1 = 1\otimes v_{k\Lambda_0} \in N(k\Lambda_0).
\tag3.5
$$
Then $L_0$ satisfies
$$
L_0\cdot x_1(n_1)\cdots x_r(n_r)v_{k\Lambda_0}
= \sum_j (-n_j)x_1(n_1)\cdots x_r(n_r)v_{k\Lambda_0}.\tag3.6
$$
Let $L_{-1}$ be the operator on $N(k\Lambda_0) $
determined by
$$
L_{-1}\cdot x_1(n_1)\cdots x_r(n_r)v_{k\Lambda_0}
= \sum_j (-n_j)x_1(n_1)\cdots x_j(n_j-1)\cdots
x_r(n_r)v_{k\Lambda_0}.\tag3.7
$$
Then $N(k\Lambda_0)$, $L_{-1}$, $L_0$, $\1$ satisfy
(2.17--2.19, 2.30). For
$x\in\g$ consider the generating function of operators on
$N(k\Lambda_0)$
$$
x(z) = \sum_{n\in\Z}x(n)z^{-n-1}.
\tag3.8
$$
When we think of $x(n)$ as a coefficient of $x(z)$, we
sometimes write $x_n$ instead.
It is easy to show that each $x(z)$  is admissible of
conformal weight
1 in the sense of Definition 2.2. Set
$$
U = \g(-1)\otimes v_{k\Lambda_0} = N(k\Lambda_0)_1\subset
N(k\Lambda_0).
\tag3.9
$$
We shall often identify $U$ with $\g$.
Define
$$
Y:U\to(\End_\F N(k\Lambda_0))[[z,z^{-1}]]
\tag3.10
$$
by
$$
Y(x(-1)\otimes v_{k\Lambda_0}, z) = x(z),
\tag3.11
$$
for $x\in\g$. Since for $x,y\in\g$,
$$
[x(z_1),y(z_2)] = z_2^{-1}{\delta}(z_1/z_2)[x,y](z_2)
-k(x,y)z_2^{-1}\frac{\partial}{\partial z_1}{\delta}(z_1/z_2),
\tag3.12
$$
$Y$ satisfies the hypotheses of Theorem 2.6. Hence $Y$ extends to make
$N(k\Lambda_0) $ a vertex operator algebra except that (so far)  only $L_{-1},
L_0$ out of the Virasoro algebra act on $N(k\Lambda_0)$.

To complete the construction of a vertex operator algebra we need a
conformal vector ${\omega}$ giving rise to the Virasoro algebra. For this we
use the following so called Sugawara construction.
Let $\{x^i\}_{i\in I}$, $\{y^i\}_{i\in I}$ be dual bases in $\g$. Set
$$
{\phi} = \sum_i(x^i)_{-1}(y^i)_{-1}\1
= \sum_i(x^i)_{-1}y^i \in N(k\Lambda_0)_2.
\tag3.13
$$
There is a nondegenerate bilinear form on $\g\otimes\g$ determined by
$$
(u\otimes v, x\otimes y) = (u,x)(v,y).
$$
For $x\in\g$, $m, n \in I$,
$$
\split
&(\sum_i([x,x^i]\otimes y^i+x^i\otimes [x,y^i]), y^m\otimes x^n) \\
&=([x,x^n],y^m)+([x,y^m],x^n) = 0.
\endsplit
$$
Hence
$$
\sum_i([x,x^i]\otimes y^i+x^i\otimes[x,y^i]) = 0\in\g\otimes\g.
$$
It follows that
$$
x_0{\phi} = \sum_i\{[x,x^i]_{-1}y^i + (x^i)_{-1}[x,y^i]\} = 0
\tag3.14
$$
for all $x\in\g$.
We have
$$
\split
x_1{\phi} &= \sum_i\{[[x,x^i],y^i] + k(x,x^i)y^i+k(x,y^i)x^i\}\\
&=2kx+\sum_i[[x,x^i],y^i]\\
&=[2k+(\theta,\theta+2{\rho})]x
\endsplit
\tag3.15
$$
if $\theta$ is the highest root and ${\rho}$ half the sum of the
positive roots.
Now
$$
({\rho},\theta) = {\rho}(\nu^{-1}(\theta))={\rho}(h_\theta)
={\rho}(\sum a_j^{\vee}{\alpha}_j^{\vee}) = \sum a_j^{\vee} = g-1
\tag3.16
$$
(where $g$ is the dual Coxeter number, and $\nu$ is the isomorphism between
the Cartan subalgebra of $\goth g$ and its dual defined via the bilinear 
form $(\cdot,\cdot)$, cf. [K]) so that
$$
x_1{\phi} = 2(k+g)x.
$$
Since there is no $\g$-invariant map $\g\to\F$ we must have
$$x_2{\phi}=0
$$
for all $x\in\g$, and by degree consideration $x_n{\phi} = 0$ for $n\ge
3$. The commutator formula \cite{FLM\rm, (8.6.3), (A.3.11)} gives
$$
\align
[x(z_1),Y({\phi},z_2)] &= -2(k+g)z_2^{-1}\frac{\partial}{\partial z_1}{\delta}
(z_1/z_2)x(z_2), \\
[{\phi}_{m+1}, x_n] &= -2(k+g)nx_{m+n}.
\endalign
$$
If $k\ne -g$ we define
$$
\align
{\omega} &= \frac{1}{2(k+g)}{\phi},\tag3.17 \\
L_m &= {\omega}_{m+1}\tag3.18
\endalign
$$
and have
$$
[L_m,x_n] = -nx_{m+n}
\tag3.19
$$
so that $L_0, L_{-1}$ coincide with the operators above. Also,
$$
\align
L_0{\omega} &= 2{\omega}, \\
L_1{\omega} &= \frac{1}{2(k+g)}\sum_i[x^i,y^i] = 0, \\
L_2{\omega} &= \frac{1}{2(k+g)}\sum_i k(x^i,y^i)\1 = \frac{k\dim\g}{2(k+g)}\1
\endalign
$$
and $L_n{\omega} = 0$ for $n\ge 3$. It now follows from the commutator formula
that
$$
[L_m,L_n] = (m-n)L_{m+n} + \frac{m^3-m}{12}{\delta}_{m+n,0}\frac{k\dim\g}{k+g}
\tag3.20
$$
so that the Virasoro algebra acts on $N(k\Lambda_0)$. We can summarize
the results of this section as:

\proclaim{Theorem 3.1} For every $k\ne -g$,
$(N(k{\Lambda}_0),Y,\1,{\omega})$ is a vertex operator algebra with 
$$
\rank N(k\Lambda_0) = \frac{k\dim\g}{k+g}.
\tag3.21
$$
\endproclaim

\head{4. Modules over $N(k\Lambda_0)$}\endhead

In this section we show that when $k \ne -g$ it is possible
to extend
the operators from $\tilde\g$ on any highest weight module
of level $k$
so that the $\tilde\g$-module becomes a module over the
vertex operator
algebra $N(k\Lambda_0)$ in the sense of \cite{FLM\rm, 8.10}.

\remark{Remark} If $\g = {\frak {sl}}(2,\F)$ and
$M({\Lambda})$ is any
Verma module over $\tilde\g$ of level $-g = -2$ then there
is no
operator $L_{-1}$ on $M({\Lambda})$ such that $[L_{-1},u(n)]
=
-nu(n-1)$ for $u\in\g$ and $[L_0,L_{-1}] = L_{-1}$. It is
not even possible to
define $L_{-1}v_{\Lambda}$ so as to satisfy
$$
x(0)L_{-1}v_{\Lambda} = 0,\qquad y(1)L_{-1}v_{\Lambda} =
y(0)v_{\Lambda}.
$$
To show this, note that since the weight of
$L_{-1}v_{\Lambda}$ is determined,
$L_{-1}v_{\Lambda}$ would have to be a linear combination of
$h(-1)v_{\Lambda}$
and $x(-1)y(0)v_{\Lambda}$. It is easy to check that the
above system
of equations has no solution.
\endremark

Fix a splitting Cartan subalgebra $\h$ in $\g$, let
${\Phi}\subset \h^*$ be the
associated root system and choose a set of simple roots
$\{{\alpha}_1,\dots,{\alpha}_\ell\}\subset{\Phi}$. Let $W$
be a level
$k$ module over $\hat\g\rtimes (\F L_0\oplus\F
L_{-1})=\tilde\g\rtimes \F L_{-1}$, a Lie
algebra with the action of $L_0, L_{-1}$ on $\ghat$
determined by (3.19).
Consider $\h^e = \h(0)\oplus\F c\oplus\F L_0$ as Cartan
subalgebra
in $\ghat\rtimes (\F L_0\oplus\F L_{-1})$. Assume that
$$
W = \coprod_{{\mu}\in\h^{e*}} W_{\mu},
\tag4.1
$$
$$
\dim W_{\mu} < \infty
\tag4.2
$$
and that there is ${\tau}\in \F$ such that  for every weight ${\mu}$
of $W$
$$
{\mu}(L_0) \in {\tau}+\N.
\tag4.3
$$

\remark{Remark} An easy calculation (cf\. \cite{H}) shows that the relation
$$
[L_m, x_n] = -nx_{m+n},
$$
from the Sugawara construction in Section 3, takes place in the
completed enveloping algebra $\overline {U(\tilde\goth g)}$ \cite{MP} 
when $k\ne -g$. Hence the
Sugawara operator $L_{-1}$ can be used to show that when $k\ne -g$ any
highest weight module of level $k$ satisfies the above conditions.
\endremark

Set
$$
W' = \coprod_{{\mu}\in\h^{e*}} W_{\mu}^*
$$
and denote the pairing between $W'$ and $W$ also by
$$
\la \cdot,\cdot\ra : W'\times W\to \F.
$$
For simplicity, set $V = N(k\Lambda_0)$. Our main goal is to
construct for each $w\in W$ a so called intertwining operator
$$
Y(w,z) \in (\Hom_\F(V,W))[[z,z^{-1}]].
$$
For
$$
u = x(-1)\otimes v_{k{\Lambda}_0}\in U,
$$
we change the meaning of $Y(u,z)$ so as to denote the generating
function 
$$
Y(u,z) = x(z) = \sum_{n\in\Z} x(n)z^{-n-1}
$$
of operators on $V\oplus W$ leaving invariant $V$ and $W$. The
intertwining operator should be local with respect to this extension i\.e\.
$$
(z_1-z_2)^N[Y(u,z_1), Y(w,z_2)] = 0
$$
for some $N\in\N$, and all $u\in U$.

With $T(\ghat)$ the tensor algebra, consider the associative algebra
$$
T_k = T(\ghat)/(c-k).
$$
Then $T_k$ has a basis consisting of monomials
$$
x_1(m_1)\otimes \cdots \otimes x_r(m_r)
$$
with the $x_j$ chosen from a basis of $\g$. Let $u_2,\dots,u_n\in U$,
$w\in W$, $w'\in W'$. Set
$$
P = \prod_{2\le\ell < m\le n}(z_\ell - z_m)^2.
$$
Then
$$
\split
P&\la w', Y(u_2,z_2)\cdots Y(u_n,z_n)e^{z_1L_{-1}}w\ra \\
&=P\la w', e^{z_1L_{-1}}Y(u_2,z_2-z_1)\cdots Y(u_n,z_n-z_1)w\ra \\
&=P\la w', e^{z_1L_{-1}}Y(u_{{\pi}2},z_{{\pi}2}-z_1)\cdots
Y(u_{{\pi}n},z_{{\pi}n}-z_1)w\ra
\endsplit
$$
for all ${\pi}\in \Sym\{2,\dots,n\}$. It follows that there is $N\in
\N$ such that if
$$
Q = \prod_{\ell = 2}^n(z_\ell - z_1)^N,
$$
then
$$
PQ\la w', Y(u_2,z_2)\cdots Y(u_n,z_n)e^{z_1L_{-1}}w\ra \in \F[z_1,\dots,z_n].
\tag4.4
$$
We have also that
$$
\la w', Y(u_2,z_2)\cdots Y(u_n,z_n)e^{z_1L_{-1}}w\ra\in
\F((z_2))\cdots ((z_n))((z_1)),
$$
hence there is $f\in\R_n$ such that
$$
\la w', Y(u_2,z_2)\cdots Y(u_n,z_n)e^{z_1L_{-1}}w\ra =
{\iota}_{2\cdots n1}f.
\tag4.5
$$
Clearly the map
$$
\align
\g\times\cdots\times\g&\to\R_n \\
(u_2,\dots,u_n)&\mapsto f
\endalign
$$
is multilinear. Consideration of weights with respect to $\h^e$ now
show that we can define
$$
Y(w,z) \in (\Hom_\F(T_k,W))[[z,z^{-1}]]
$$
by the condition that
$$
\la w', Y(w,z_1)Y(u_2,z_2)\cdots Y(u_n,z_n)\ra = {\iota}_{12\cdots n}f
\tag4.6
$$
for all $w'\in W'$ with $f$ as in (4.5). Here we use $\dim W_{\mu} < \infty$
in order that $W_{\mu}^{**} = W_{\mu}$.

Set
$$
U_k(\ghat) = U(\ghat)/(c-k),\qquad  U_k(\tilde\g)=U(\tilde
\goth g)/(c-k).
$$
By the construction of the enveloping algebra, $U_k(\ghat)$ is then the
quotient of $T_k$ by the two-sided ideal generated by
$$
x(m)y(n)-y(n)x(m)-[x,y](m+n)-m{\delta}_{m+n,0}(x,y)k
$$
for all $x,y\in\g$, $m,n\in\Z$. We want to show that the components of
$Y(w,z)$ induce operators $U_k(\ghat)\to W$. We thus consider the effect of
interchanging the operators $Y(u_j,z_j)$, and $Y(u_{j+1}, z_{j+1})$ in
(4.5). For clarity we raise the indices on ${\iota}$ so that
${\iota}_{a,b,\dots,m} = {\iota}(a,b,\dots,m)$. Multiplication by
$(z_j-z_{j+1})^2$ show that
$$
\split
\la w', &Y(u_2,z_2)\cdots Y(u_{j+1},z_{j+1})Y(u_j,z_j)\cdots
Y(u_n,z_n)e^{z_1L_{-1}}w\ra = \\
&= {\iota}(2,\dots,j+1,j,\dots,n,1)f
\endsplit
$$
with $f$ as in (4.5). By (4.4) the order of the pole $(z_j-z_{j+1})^{-1}$
in $f$ is at most 2. It is easy to see that there is a unique representation
$$
f = \frac{f_1}{z_j-z_{j+1}} + \frac{f_2}{(z_j-z_{j+1})^2} + p
\tag4.7
$$
with $f_1,f_2\in\R_{\{1,\dots,\hat j,\dots,n\}}$ (i\.e\. independent of $z_j$)
and with $p$ having the poles of functions in $\R_n$ except that $p$
does not have the pole $(z_j-z_{j+1})^{-1}$. It follows that
$$
\split
\{&{\iota}(2,\dots,n,1)-{\iota}(2,\dots,j+1,j,\dots,n,1)\} f = \\
&= z_{j+1}^{-1}{\delta}(z_j/z_{j+1}){\iota}(2,\dots,\hat
j,\dots,n,1)f_1 \\
&-z_{j+1}^{-1}\frac{\partial}{\partial z_j}
{\delta}(z_j/z_{j+1}){\iota}(2,\dots,\hat j,\dots,n,1)f_2.
\endsplit
\tag4.8
$$
Let $x,y\in\g$ so that
$$
u_j = x(-1)\otimes v_{k{\Lambda}_0},\qquad
u_{j+1} = y(-1)\otimes v_{k{\Lambda}_0}
$$
and set
$$
\bar u = [x,y](-1)\otimes v_{k{\Lambda}_0}\in U.
$$
There are $g\in\R_{\{1,\dots,\hat j,\dots,n\}}$,
$h\in\R_{\{1,\dots,\hat j,\widehat{j+1},\dots,n\}}$ such that
$$
\align
\la &w', Y(u_2,z_2)\cdots Y(u_{j-1},z_{j-1})Y(\bar u,z_{j+1})\cdots
Y(u_n,z_n)e^{z_1L_{-1}}w\ra = \\
&= {\iota}(2,\dots,\hat j,\dots,n,1)g,
\endalign
$$
and
$$
\align
\la &w', Y(u_2,z_2)\cdots Y(u_{j-1},z_{j-1})Y(u_{j+2},z_{j+2})\cdots
Y(u_n,z_n)e^{z_1L_{-1}}w\ra \\
&= {\iota}(2,\dots,\hat j,\widehat{j+1},\dots,n,1)h.
\endalign
$$
Since $W$ is a $\ghat$-module we have
$$
\split
\{&{\iota}(2,\dots,n,1)-{\iota}(2,\dots,j+1,j,\dots,n,1)\}f = \\
&= z_{j+1}^{-1}{\delta}(z_j/z_{j+1}){\iota}(2,\dots,\hat j,\dots,n,1)g \\
&-k(x,y)z_{j+1}^{-1}\frac{\partial}{\partial z_j}
{\delta}(z_j/z_{j+1}){\iota}(2,\dots,\hat j,\widehat{j+1},\dots,n,1)h.
\endsplit
$$
Comparison with (4.8) yields
$$
f_1 = g,\qquad f_2 = k(x,y)h
$$
(after multiplication by $(z_j-z_{j+1})$). (Thus $f_2$ does not depend
on $z_{j+1}$.) Hence
$$
\split
\{&{\iota}(1,2,\dots,n)-{\iota}(1,2,\dots,j+1,j,\dots,n)\}f = \\
&= z_{j+1}^{-1}{\delta}(z_j/z_{j+1}){\iota}(1,2,\dots,\hat j,\dots,n)g \\
&-k(x,y)z_{j+1}^{-1}\frac{\partial}{\partial z_j}
{\delta}(z_j/z_{j+1}){\iota}(1,2,\dots,\hat j,\widehat{j+1},\dots,n)h
\endsplit
$$
so that
$$
\split
\la &w', Y(w,z_1)Y(u_2,z_2)\cdots Y(u_n,z_n)\ra -
\la w', Y(w,z_1)\cdots Y(u_{j+1},z_{j+1})Y(u_j,z_j)
\cdots\ra \\
&=z_{j+1}^{-1}{\delta}(z_j/z_{j+1})
\la w', Y(w,z_1)Y(u_2,z_2)\cdots Y(\bar u,z_{j+1})\cdots Y(u_n,z_n)\ra \\
&-k(x,y)z_{j+1}^{-1}\frac{\partial}{\partial z_j}
{\delta}(z_j/z_{j+1})
\la w', Y(w,z_1)\cdots Y(u_{j-1},z_{j-1})Y(u_{j+2},z_{j+2})
\cdots\ra.
\endsplit
$$
This means that each component of $Y(w,z)$ annihilates the two-sided
ideal which is the kernel of $T_k\to U_k(\ghat)$ and hence we have well
defined operators
$$
Y(w,z) \in (\Hom_\F(U_k(\ghat),W))[[z,z^{-1}]].
$$
Furthermore, from (4.4) there are no poles $z_n^{-1}$ in $f$
satisfying (4.5), hence
$$
{\iota}_{12\cdots n}f = \sum_{d\le -1}f_dz_n^{-d-1}
$$
with $f_d\in\F((z_1))\cdots((z_{n-1}))$ so that
$$
\la w', Y(w,z_1)Y(u_2,z_2)\cdots Y(u_{n-1},z_{n-1})(u_n)_d\ra = 0
$$
if $d\ge 0$. Hence all components of $Y(w,z)$ annihilate all vectors
in $U_k(\ghat)$ of the form
$$
x_2(m_2)\cdots x_{n-1}(m_{n-1})x_n(d)
$$
with $d\ge 0$. Since these span the left ideal in $U_k(\ghat)$ which is the
kernel of
$$
U_k(\ghat)\to N(k{\Lambda}_0),
$$
$Y(w,z)$ induces a well defined operator
$$
Y(w,z) \in (\Hom_\F(N(k{\Lambda}_0),W))[[z,z^{-1}]].
$$
This is the intertwining operator that we want. We have proved:

\proclaim{Proposition 4.1} Let $u_2,\dots,u_n\in U$, $w\in W$,
$w'\in W'$.
\item {\rm (i)} There is $f\in \R_n$ such that
$$
\la w', Y(u_2,z_2)\cdots Y(u_n,z_n)e^{z_1L_{-1}}w\ra =
{\iota}_{2\cdots n1}f.
$$
\item {\rm (ii)} Using these matrix coefficients $f$, there is a well
defined operator
$$
Y(w,z) \in (\Hom_\F(N(k{\Lambda}_0),W))[[z,z^{-1}]]
$$
determined by
$$
\la w', Y(w,z_1)Y(u_2,z_2)\cdots Y(u_n,z_n)\1\ra = {\iota}_{12\cdots n}f.
\tag4.9
$$
\endproclaim

Now it is easy to show that $Y(w,z)$ is local with respect to the
operators $Y(u,z)$, $u\in U$:

\proclaim{Proposition 4.2} Let $u_0,u_2,\dots,u_n\in U$, $w\in W$,
$w'\in W'$. Then
$$
\la w', Y(u_0,z_0)Y(w,z_1)Y(u_2,z_2)\cdots Y(u_n,z_n)\1\ra
= {\iota}_{012\cdots n}f
\tag4.10
$$
if $f\in\R_{n+1}$ is such that
$$
\la w', Y(u_0,z_0)Y(u_2,z_2)\cdots Y(u_n,z_n)e^{z_1L_{-1}}w\ra =
{\iota}_{02\cdots n1}f.
$$
If $N\in\N$ is such that $(u_0)_dw = 0$ when $d\ge N$, then the order of
the pole $(z_0-z_1)^{-1}$ in $f$ is at most $N$.
Hence $Y(w,z)$ is local with respect to $Y(u_0,z)$ for all $u_0\in U$.
\endproclaim
\demo{Proof} Let
$$
{\iota}_{0\cdots}f = \sum_{d\in\Z}f_dz_0^{-d-1}.
$$
Then
$$
\la w', (u_0)_dY(u_2,z_2)\cdots Y(u_n,z_n)e^{z_1L_{-1}}w\ra =
{\iota}_{2\cdots n1}f_d
$$
by Proposition 1.1.3. The definition of $Y(w,z_1)$ gives
$$
\la w', (u_0)_dY(w,z_1)Y(u_2,z_2)\cdots Y(u_n,z_n)\1\ra
= {\iota}_{12\cdots n}f_d,
$$
and (4.10) now follows from Proposition 1.1.3. From
$$
\split
&\prod_{\ell = 2}^n (z_0-z_\ell)^2{\iota}_{02\cdots n1}f = \\
&=\prod_{\ell = 2}^n (z_0-z_\ell)^2
\la w', Y(u_0,z_0)Y(u_2,z_2)\cdots Y(u_n,z_n)e^{z_1L_{-1}}w\ra \\
&=\prod_{\ell = 2}^n (z_0-z_\ell)^2
\la w', e^{z_1L_{-1}}Y(u_2,z_2-z_1)\cdots Y(u_n,z_n-z_1)Y(u_0,z_0-z_1)w\ra
\endsplit
$$
we can determine the order of the pole $(z_0-z_1)^{-1}$ in $f$.
\qed \enddemo

\proclaim{Theorem 4.3} Let $k\ne -g$, and let $W$ be a
$\ghat\rtimes(\F
L_0\oplus\F L_{-1})$-module of level $k$ satisfying (4.1--3).
(In
particular $W$ could be any highest weight module or a
module from 
the category $\Cal O$.) Then there is a
unique extension of the operators $Y(u,z)$, $u\in U$ on $W$
that make
$W$ into a module over the vertex operator algebra
$N(k{\Lambda}_0)$.
\endproclaim
\demo{Proof} Let $\tilde V$ be the space of all admissible
generating functions
$$
a(z) = \sum_{n\in\Z} a_nz^{-n-1}
$$
of operators $a_n$ on $V\oplus W$ stabilizing both $V$ and
$W$ such
that
$$
a(z) \text{ is local w\.r\.t\. } \{Y(u,z)\mid u\in
U\}\cup\{Y(w,z)\mid w\in W\}.
$$
Then the proof of Theorem 2.6 applies with only minor
modifications
to show that
\item {1.} the map $\tilde V\to V$, $a(z)\mapsto a_{-1}\1$
is injective,
\item {2.} the image of this map is all of $V$,
\item{3.} any two $a(z),b(z)\in \tilde V$ are mutually
local.

Given this, the  argument of \cite{FLM\rm, A.2.5} can be
employed, using the
intertwining operator $Y(w,z)$ in the r\^ole of one vertex
operator, to show
that $W$ satisfies the ``associativity condition"
\cite{FLM\rm, (8.10)} and
hence the Jacobi identity.
\qed\enddemo

We have seen that for level $k$ highest weight modules, or
modules from the 
category $\Cal O$, the action of $\tilde\goth g$ determined
by 
$x(z)=\sum x_nz^{-n-1}$ extends, by using normal order
products (2.14), to
the action of $N(k\Lambda_0)$ defined by $Y(v,z)=\sum
v_nz^{-n-1}$, $v\in N(k\Lambda_0)$.
While $x_n$ are elements in $U_k(\tilde\goth g)$, in general
the coefficients $v_n$ are 
infinite sums of elements in $U_k(\tilde\goth g)$ and should
be seen in some kind of
completion $\overline {U_k(\tilde\goth g)}$. We prefer to
think in terms of completion
defined in \cite{MP}: a net $(y_i)_{i\in I}$ converges to
$y$ if for each 
$\tilde\goth g$-module $V$ in the category $\Cal O$ and each
vector $v$ in $V$ there 
is an index $i_0\in I$ such that $i\ge i_0$ implies
$y_i\cdot v=y_{i_0}\cdot v=y\cdot v$.

\remark{Remark}
Constructions and results analogous to Theorem 3.1 and Theorem 4.3
have appeared in \cite{FF} and in relation to the
Wess-Novikov-Zumino-Witten model in the physics literature. Theorems
3.1 and 4.3 contain the precise results we require in the following sections.

It is easy to see that normal order product (2.14) coincides
with the normal order product
defined in \cite{H} if one identifies $\g\cong \g(-1)\1$
and
$S(\tilde\goth g_{<0})\cong U(\tilde\goth g_{<0})\cong
N(k\Lambda_0)$ via the symmetrization map.
Hence we work with the same set of fields, but the ``index
set" $N(k\Lambda_0)$ brings into
the game the action of $\tilde\goth g$.

Before our work was completed, other proofs of Theorems 3.1
and 4.3 appeared: in \cite{FZ} 
by using quite different ideas, and recently in \cite{Li} by
using the locality in a similar 
way, but from a different point of view, and proving in an
elegant way stronger and more 
general results.

\endremark

\head{5. Relations on standard modules}\endhead

We identify $\tilde\g$
with a Kac-Moody algebra in the usual way (cf\. \cite{K}):
first choose root vectors
$x_{\alpha}\in\g_{\alpha}$,
such that
$h_{\alpha} = [x_{\alpha},x_{-{\alpha}}]$ satisfy
${\alpha}(h_{\alpha}) = 2$
for ${\alpha}\in{\Phi}$, and then set
$e_i = x_{{\alpha}_i}(0)$, $f_i = x_{-{\alpha}_i}(0)$, for
$i = 1,\dots, \ell$,
$e_0 = x_{-{\theta}}(1)$, $f_0 = x_{\theta}(-1)$ (${\theta}$
the highest
root as in Section 3) and $h_i = [e_i,f_i]$
for $i = 0,1,\dots,\ell$. Then we have the usual triangular
decompositions
$\g=\goth n_-+\h+\goth n_+$ and $\tilde\goth g=\tilde\goth
n_-+\h^e+\tilde\goth n_+$, where
$\h^e=\h+\Bbb Fc+\Bbb F L_0$. Denote by
$$
P_+ =
\{{\Lambda}\in\h^{e*}\mid {\Lambda}(h_i)\in\N \text{ for }
i = 0,1,\dots,\ell\},
\tag5.1
$$
the set of dominant integral weights, and let for $i = 0,1,\dots,\ell$,
${\Lambda}_i\in P_+$ be the fundamental weight determined by
${\Lambda}_i(h_j) = {\delta}_{ij}$, ${\Lambda}_i(L_0) = 0$.

For ${\Lambda}\in\h^{e*}$ we shall use the notation
$M({\Lambda})$ for
the Verma module with highest weight ${\Lambda}$ and
$L({\Lambda})$
for the irreducible quotient of $M({\Lambda})$.
In this section we show that the standard $\tilde\g$-modules
$L({\Lambda})$ of level $k$, i\.e\. ${\Lambda}\in P_+$,
${\Lambda}(c) =
k$, satisfy the relation $x_{\theta}(z)^{k+1} = 0$, and that
in some
sense this relation determines the structure of
$L({\Lambda})$, see
Theorem 5.14. The relation $x_{\theta}(z)^{k+1} = 0$ and its
consequences will be studied further in the following
sections.

For a given nontrivial finite dimensional irreducible $\goth g$-module
$W$ we have
the associated loop module 
$\bar W=W\otimes \Bbb C [t, t^{-1}]$
with the action of $\tilde\goth g$ defined by (cf. [CP]):
$$
x\otimes t^n\cdot w\otimes t^m=x\cdot w\otimes t^{n+m},\quad
c\cdot w\otimes t^m=0,\quad
L(0)\cdot w\otimes t^m=-m\,w\otimes t^m.
$$
The results in this section are mostly about the loop module
$\bar R$ constructed from the coefficients of the
annihilating field $x_{\theta}(z)^{k+1}$.

\proclaim{Theorem 5.1} Let $k\in\F$ and let $R$ be a subspace of
$N(k{\Lambda}_0)$ invariant under $\g = \g(0)$ and $L_0$. Set
$$
\bar R = \F\text{-span}\{r_n\mid r\in R, n\in\Z \}\subset \End_\F N(k{\Lambda}_0).
\tag 5.2
$$
\item{\rm (i)} The following are equivalent:
\itemitem{\rm (a)} For all $x\in\g$, $r\in R$, $m,n\in\Z$,
$$
[x(m),r_n] = (x(0)r)_{m+n}
\tag 5.3
$$
so that $\bar R$ is a loop module under the ``adjoint" action of $\tilde\g$.
\itemitem{\rm (b)} We have $\g(n)R = 0$ for all positive integers $n$.
\item{\rm (ii)} If {\rm (a)} and {\rm (b)} are satisfied then $\bar R$ is irreducible
iff $R$ is a nontrivial irreducible $\g$-module.
\endproclaim

\proclaim{Proposition 5.2} There is a bijection between
\roster
\item"{(a)}" the
irreducible loop modules  $\bar R$ associated to subspaces $R$ of 
$N(k{\Lambda}_0)$ invariant under $\g$ and $L_0$ as in Theorem 5.1,
\item"{(b)}" the singular vectors (up to nonzero scalar multiples)
in $N(k{\Lambda}_0)$ that are not invariant
under $\g$.
\endroster
\endproclaim
\demo{Proof} If $\bar R$ is an irreducible loop module then $R$ is a nontrivial
irreducible $\g$-module by Theorem 5.1(ii) and is hence by the theorem
of the highest weight determined by a highest weight vector with
respect to $({\frak h}, {\Phi}^+)$. This vector is by condition (b) in
Theorem 5.1 also a singular vector. We must also show that $R$ is determined by
$\bar R$. Since $R$ is irreducible under $\g$, $R$ must be
homogeneous. It then follows from $Y(r,z)\1 = e^{zL_{-1}}r$ that
$$
\bar R\1 = \coprod_{n\ge 0}L_{-1}^nR
$$
and hence $R$ is determined as the homogeneous component of highest
degree that is nonzero in $\bar R\1$.
\qed \enddemo

\demo{Proof of Theorem 5.1}
(i) For $x\in \g$, $r\in R$,
$$
[x(z_1), Y(r,z_2)] = \sum_{i\ge 0}\frac{(-1)^i}{i!}(\frac
\partial{\partial z_1})^iz_2^{-1}
{\delta}(z_1/z_2)Y(x(i)r,z_2).
$$
Hence (b) implies (a). Furthermore, if (a) holds and $N\in\N$ such
that $x(i)r = 0$ for $i > N$, then
$$
\sum_{i = 1}^N\frac{(-1)^i}{i!}(\frac
\partial{\partial z_1})^i
{\delta}(z_1/z_2)Y(x(i)r,z_2) = 0,
$$
and
$$
\split
\coeff_{z_1^0}&\sum_{i = 1}^N\frac{(-1)^i}{i!}(\frac
\partial{\partial z_1})^i
{\delta}(z_1/z_2)Y(x(i)r,z_2)\1 = \\
&= \sum_{i = 1}^N(-1)^iz_2^{-i}(x(i)r + z_2(\cdots)) = 0,
\endsplit
$$
so that $x(N)r = 0$. Repeating this argument shows successively that
$x(N-1)r = 0$,\dots, $x(1)r = 0$, and hence (b).

(ii) If $S$ is a nonzero proper submodule of $R$ then clearly $\bar
S$ is a nonzero proper submodule of $\bar R$.
If $R$ is a 1-dimensional trivial $\g$-module, then $\F r_n$ is a
$\tilde\g$-submodule of $\bar R$ for each $r\in R, n\in\Z$.
Conversely, assume that
$R$ is nontrivial and irreducible and let 
$Z$ be a nonzero $\tilde\g$-submodule of $\bar R$. From $[L_0,
Z]\subset Z$ it follows that $Z = \oplus_{n\in\Z}Z(n)$ with $Z(n) =
\{p\in Z\mid [L_0,p] = np \}$. For some $n$, $Z(n)\ne 0$ and then $[\g,
Z(n)]\subset Z(n)$ implies that
$$
Z(n) = \{r_m\mid r\in R \}
$$
for an appropriate $m$. It now follows from (5.3) that $Z = \bar R$.
\qed \enddemo

\proclaim{Proposition 5.3} \cite{K\rm, Corollary 10.4} Let
${\Lambda}\in P_+$. Then the annihilator in $U(\tilde\g)$ of a highest
weight vector of  $L({\Lambda})$ is the left ideal
generated by $e_i$, $f_i^{{\Lambda}(h_i)+1}$, and $h-{\Lambda}(h)$
where $i = 0,1,\dots,\ell$ and $h\in\h^e$.
\endproclaim

For $k\in\F$ let $N^1(k{\Lambda}_0)$ denote the maximal $\tilde\g$-submodule of $N(k{\Lambda}_0)$.

{\it
In the remainder of this section we shall assume that $k$ is a
positive integer.
}

\proclaim{Corollary 5.4} The maximal submodule $N^1(k{\Lambda}_0)$ of
$N(k{\Lambda}_0)$ is generated by the vector  $x_{\theta}(-1)^{k+1}\1$.
\endproclaim
\demo{Proof} This follows immediately from Proposition 5.3 since for $i
= 1, \dots,\ell$,
$$
f_i^{k{\Lambda}_0(h_i)+1}\1 = f_i\1=x_{-{\alpha}_i}(0)\1=0,
$$
and we have
$f_0^{k{\Lambda}_0(h_0)+1}\1 = f_0^{k+1}\1 = x_{\theta}(-1)^{k+1}\1$.
\qed \enddemo

Set
$$
R = U(\g)x_{\theta}(-1)^{k+1}\1\subset N(k{\Lambda}_0).
\tag 5.4
$$
Since $\frak n_+(0)$ annihilates $x_{\theta}(-1)^{k+1}\1$, and
$N(k{\Lambda}_0)$ is a sum of finite dimensional $\g(0)$-modules,
$R$ is an
irreducible $\g$-module.

\proclaim{Lemma 5.5} The vector $x_{\theta}(-1)^{k+1}\1\in N(k{\Lambda}_0)$
is a singular vector. Hence $U(\ghat)x_{\theta}(-1)^{k+1}\1 =
U(\tilde{\frak n}_-)x_{\theta}(-1)^{k+1}\1$. 
\endproclaim
\demo{Proof} It remains to show that
$x_{-{\theta}}(1)x_{\theta}(-1)^{k+1}\1 = 0$. From
$$
[x_{-{\theta}}(1),x_{\theta}(-1)] = -h_{\theta}+c,
$$
and $(-h_{\theta}+c)\1 = k\1$
it follows by a standard calculation in the representation theory of
$\frak{sl}(2,\F)$ that $x_{-{\theta}}(1)x_{\theta}(-1)^{k+1}\1 = 0$.
\qed \enddemo

For $r\in R$, define $r(n)$ by
$$
Y(r,z) = \sum_{n\in\Z}r(n)z^{-n-k-1}.
\tag 5.5
$$
Then $r(n)$ has degree $n$ and $r(n) = r_{n+k}$. By Theorem 5.1
$$
[x(m),r(n)] = (x(0)r)(m+n)
\tag 5.6
$$
for $x\in\g$, $r\in R$, $m,n\in\Z$, so that
$$
\bar R = \F\text{-span}\{r(n)\mid r\in R, n\in \Z \}
\tag 5.7
$$
is a loop module. Let us record this as:

\proclaim{Corollary 5.6} The coefficients of\/ $Y(ux_{\theta}(-1)^{k+1}\1,z)$
with $u$ ranging through $U(\g)$ span a loop module.
\endproclaim

\remark{Remark} The above statement can also be verified by direct
calculation of commutators.
\endremark

By the result of V\. Chari and A\. Pressley \cite{CP\rm, Theorem 4.5} we
have that
$
\bar R\otimes L({\Lambda}) 
$
is irreducible
for each standard module $L({\Lambda})$ of level $k$, and hence the map
$$
\align
\bar R\otimes L({\Lambda})&\to L({\Lambda}) \\
u\otimes v&\mapsto uv
\endalign
$$
must be zero, i.e.  $\bar R$ annihilates $L({\Lambda})$. We shall
give another proof of this fact:

\proclaim{Lemma 5.7}
Let
$$
x = \left(\matrix 0& 1 \\
0& 0 \endmatrix\right),
\quad y = \left(\matrix 0& 0 \\
1& 0 \endmatrix\right),
\quad h = \left(\matrix 1& 0 \\
0& -1 \endmatrix\right),
\tag 5.8
$$
be the usual basis of $\g = \frak{sl}(2,\F)$, and set
${\frak n}_+ = \F x$.
Then
$$
x^{k+1}y^{k+1}\in ch(h-1)\cdots(h-k) + U(\g){\frak n}_+,
$$
for some nonzero integer $c$.
\endproclaim
\demo{Proof}
Clearly $u = x^{k+1}y^{k+1} \in p(h) + U(\g){\frak n}_+$ for some
polynomial $p$ of degree $k+1$. By applying $u$ to a highest weight
vector of a $(j+1)$-dimensional irreducible $\g$-module,
$j\in\{0,\dots,k,k+1\}$ we get $p(j) = 0$ for $j = 0,\dots,k$,
$p(k+1)\ne 0$, and hence the lemma. 
\qed \enddemo

\proclaim{Lemma 5.8} Let $\bar W = \coprod_{n\in\Z}W(n)$ be a loop module
for $\tilde\g$ which acts on a Verma module $M({\Lambda})$. Then
$$
\bar W M({\Lambda}) = M({\Lambda})\quad \text{ iff }\quad
W(0)_0v_{\Lambda}\ne 0,
$$
where $W(0)_0$ denotes the zero-weight subspace of $\bar W$.
\endproclaim
\demo{Proof} Note first that $\bar W M({\Lambda}) = M({\Lambda})$ if
and only if $v_{\Lambda}\in\bar W M({\Lambda})$. Since
$$
\bar W M({\Lambda}) = \bar W U(\tilde{\frak n}_-)v_{\Lambda} =
U(\tilde{\frak n}_-)\bar Wv_{\Lambda},
$$
the ${\Lambda}$ weight subspace of $\bar W M({\Lambda})$ equals
$W(0)_0v_{\Lambda}$, and hence the lemma holds. 
\qed \enddemo

\proclaim{Theorem 5.9} Let $L({\Lambda})$ be a standard $\tilde\g$-module of
level $k$. Then
$$
\bar RL({\Lambda}) = 0.
$$
\endproclaim
\demo{Proof} Consider first the case $\g = \frak{sl}(2,\F)$. Then
$$
R(0) = U(\g)\cdot(\sum_{j_1+\cdots +j_{k+1} = 0}y(j_1)\cdots y(j_{k+1}))
$$
(with $\cdot$ the ``adjoint" action as above and $y$ as in (5.8)) and
for $V = M({\Lambda})$,
$$
\split
R(0)_0v_{\Lambda} &=\F\Bigl(x^{k+1}\cdot\bigl(\sum_{j_1+\cdots+j_{k+1} = 0}
y(j_1)\cdots y(j_{k+1})\bigr)\Bigr)v_{\Lambda} \\
&=\F x^{k+1}y(0)^{k+1}v_{\Lambda} \\
&=\F p({\Lambda}(h))v_{\Lambda},
\endsplit
\tag 5.9
$$
where by Lemma 5.7 $p(h) = h(h-1)\cdots(h-k)$.
If ${\Lambda}(h)\in\{0,\dots,k \}$, then $R(0)_0v_{\Lambda} = 0$ and by
Lemma 5.8 $\bar RM({\Lambda})$ is a proper submodule, in particular
$\bar RM({\Lambda})\subset M^1({\Lambda})$. Since we have an exact sequence
$$
0\to M^1({\Lambda})\to M({\Lambda})\to L({\Lambda})\to 0,
$$
$\bar R$ acts trivially on $L({\Lambda})$.

In the general case we let $\frak a$ be the subalgebra isomorphic to
$\frak {sl}(2,\F)$ spanned by $x = x_{\theta}$, $y = x_{-\theta}$, $h
= h_{\theta}$, and consider $L({\Lambda})$ as an $\tilde{\frak a}$-module by
restriction. By \cite{K\rm, Theorem 10.7} $L({\Lambda})$ is a
direct sum of standard level $k$ $\tilde{\frak a}$-modules. By the
$\frak {sl}(2,\F)$-case above, $x_{\theta}(z)^{k+1}$ annihilates $L({\Lambda})$
and bracketing with $\tilde\g$ now shows that $\bar R$ annihilates
$L({\Lambda})$. 
\qed \enddemo

\definition{Definition 5.10} A subspace $I$ of a vertex operator
algebra $V$ is said to be an {\it ideal} in $V$ if and only if
$$
Y(v,z)I\subset I[[z,z^{-1}]]
$$
for all $v\in V$, i\.e\. $I$ is invariant under the action of all
components of all vertex operators.
\enddefinition

Equivalently, $I$ is an ideal if and only if $I$ is a submodule of $V$
considered as a module over itself.

\proclaim{Proposition 5.11} Let $I$ be an ideal of a vertex operator
algebra $V$. Then the vertex operator algebra structure on $V$ induces
a structure of vertex operator algebra on $V/I$.
\endproclaim
\demo{Proof} Since there is a conformal vector ${\omega}\in V$ such that
$$
Y({\omega},z) = \sum_{n\in\Z}L_nz^{-n-2},
$$
$I$ is invariant under the Virasoro algebra. From the relation
$$
Y(u,z)v = e^{zL_{-1}}Y(v,-z)u
\tag5.10
$$
valid in any vertex operator algebra \cite{FLM, (8.8.7), A.3.2; FHL}
it now follows that $Y(u,z)V\subset I[[z,z^{-1}]]$ if $u\in I$, so
that the induced operator $Y(v,z)$ on $V/I$ only depends on $v\mod I$. 
\qed \enddemo

\proclaim{Proposition 5.12} Any $\tilde\g$-submodule of $N(k{\Lambda}_0)$ is
also a vertex operator algebra ideal of $N(k{\Lambda}_0)$. There is in
particular an induced vertex operator algebra structure on
$L(k{\Lambda}_0)$. 
\endproclaim
\demo{Proof} Let $I$ be a $\tilde\g$-submodule of $N(k{\Lambda}_0)$. It then
follows from
$$
\split
&Y(u_nv,z)= \\
&=\Res_{z_0}z_0^n\{Y(u,z_0+z)Y(v,z)-\Res_{z_1}z_0^{-1}
{\delta}(\frac{z-z_1}{-z_0})Y(v,z)Y(u,z_1)\}
\endsplit
\tag5.11
$$
for $u,v\in N(k{\Lambda}_0)$
(cf. (2.25), \cite{FLM, (8.8.31)}) and induction that
$$
Y(v,z)I\subset I[[z,z^{-1}]]
$$
for all $v\in N(k{\Lambda}_0)$.
\qed \enddemo

\proclaim{Proposition 5.13} Let $\g$ be a finite dimensional split simple
Lie algebra over $\F$ and let $L({\Lambda})$ be a standard
$\tilde\g$-module of level $k$.
\roster
\item"{(i)}" We have $x_{{\theta}}(z)^{k+1} = 0$ on $L({\Lambda})$.
\item"{(ii)}" The $\tilde\g$-module structure on $L({\Lambda})$
extends to make $L({\Lambda})$ a module over the vertex operator
algebra $L(k{\Lambda}_0)$.
\item"{(iii)}" If ${\beta}\in{\Phi}$ is a short root such that
$({\beta},{\beta}) = \frac 12 ({\theta},{\theta})$ then
$$
x_{{\beta}}(z)^{2k+1} = 0 \text{ on } L({\Lambda}).
\tag5.12
$$
\item"{(iv)}" If ${\gamma}\in{\Phi}$ is a short root such that
$({\gamma},{\gamma}) = \frac 13 ({\theta},{\theta})$ then
$$
x_{{\gamma}}(z)^{3k+1} = 0 \text{ on } L({\Lambda}).
\tag5.13
$$
\endroster
\endproclaim
\demo{Proof} Part (i) is contained in Theorem 5.9.
Part (ii) follows from Theorem 4.3, Corollary 5.4, (5.11) and part (i).

(iii) By (ii)
$$
Y(u,z) = 0 \text{ on } L({\Lambda})
$$
for each $u\in N^1(k{\Lambda}_0)$, so it is enough to show that
$x_{\beta}(-1)^{2k+1}\1\in N^1(k{\Lambda}_0)$. Since $x_{\beta}(-1)$
is a root vector of a real root we have
$x_{\beta}(-1)^nv_{k{\Lambda}_0} = 0$ in $L(k{\Lambda}_0)$ for $n$
sufficiently large \cite{K} (here $ v_{k{\Lambda}_0}$ denotes a highest
weight vector in $L(k{\Lambda}_0)$). From
$x_{-{\beta}}(1)v_{k{\Lambda}_0} = 0$,
$[x_{-{\beta}}(1),x_{\beta}(-1)] = h_{\beta}+2c$,
$(h_{\beta}+2c)v_{k{\Lambda}_0} = 2kv_{k{\Lambda}_0}$ it follows from
the representation theory of $\frak{sl}(2,\F)$ that we must have
$x_{\beta}(-1)^{2k+1} v_{k{\Lambda}_0} = 0$ in $L(k{\Lambda}_0)$ so
that (in $N(k{\Lambda}_0)$)
$$
x_{\beta}(-1)^{2k+1}\1\in N^1(k{\Lambda}_0).
$$

Part (iv) is proved similarly.
\qed \enddemo

\proclaim{Theorem 5.14} Let $\g$ be a finite dimensional split simple
Lie algebra over $\F$ and let ${\Lambda}$ be dominant integral
with ${\Lambda}(c) = k$.
Then
$$
\bar RM({\Lambda}) = M^1({\Lambda})
\tag5.14
$$
where $M^1({\Lambda})$ denotes the maximal submodule of the Verma
module $M({\Lambda})$.
\endproclaim
\demo{Proof} Set $V = \bar RM({\Lambda})$. By
Theorem 5.9 $V\subset M^1({\Lambda})$. By Proposition 5.3 it
is enough to show that
$$
\{f_i^{{\Lambda}(h_i)+1}v_{\Lambda}\mid i = 0,\dots,\ell\} \subset V,
$$
with $v_{\Lambda}$ a highest weight vector in $M({\Lambda})$.
Let ${\alpha}_i$ be a long simple root.
Then
$x_{-{\alpha}_i}(-1)^{k+1}\1\in R$ (e.g. by using the action of the
Weyl group of $\g$) and it follows that
$$
x_{-{\alpha}_i}(z)^{k+1}v_{\Lambda}\in V[[z,z^{-1}]].
$$
In particular
$$
f_i^{k+1}v_{\Lambda} = x_{-{\alpha}_i}(0)^{k+1}v_{\Lambda} =
\coeff_{z^{-k-1}}x_{-{\alpha}_i}(z)^{k+1}v_{\Lambda}\in V,
$$
and now the representation theory of $\frak{sl}(2,\F)$ implies that
$f_i^{{\Lambda}(h_i)+1}v_{\Lambda} \in V$. If ${\alpha}_j$ is a
short simple root
$$
x_{-{\alpha}_j}(-1)^{nk+1}\1\in
N^1(k{\Lambda}_0) = U(\ghat)x_{\theta}(-1)^{k+1}\1,
$$
for $n = 2$ or $3$
by the argument of the proof of Proposition 5.13(iii). It follows from
(5.11) and induction that
$$
x_{-{\alpha}_j}(z)^{nk+1}v_{\Lambda}\in V[[z,z^{-1}]]
$$
so that (coefficient of $z^{-nk-1}$)
$f_j^{nk+1}v_{\Lambda}\in V$ and hence
$f_j^{{\Lambda}(h_j)+1}v_{\Lambda} \in V$. Finally,
$$
f_0^{k+1}v_{\Lambda} = x_{\theta}(-1)^{k+1}v_{\Lambda} = \coeff_{z^0}
x_{\theta}(z)^{k+1}v_{\Lambda}\in V
$$
and thus $f_0^{{\Lambda}(h_0)+1}v_{\Lambda} \in V$.
\qed \enddemo

\proclaim{Theorem 5.15} Let $V$ be a highest weight module of level
$k$. The following are equivalent:
\roster
\item"{(a)}" $V$ is a standard module,
\item"{(b)}" $\bar R$ annihilates $V$.
\endroster
\endproclaim
\demo{Proof} By Theorem 5.9 (a) implies (b). Conversely, assume that
$\bar RV = 0$, let ${\Lambda}$ be the highest weight of $V$ and
$v_{\Lambda}$ a highest weight vector. Then the same reasoning as in
the proof of Theorem 5.14 shows that for each $i = 0,1, \dots, \ell$,
$f_i^{nk+1}v_{\Lambda} \in\bar RV  = 0$ for $n = 1,2$ or $3$. Thus $V$ is a
standard module.
\qed \enddemo

\remark{Remark}
As was mentioned in the introduction, Theorem 5.9 may be
viewed as a consequence of the
Frenkel-Kac vertex operator formula (cf. \cite{FK}). We
omitted such a proof in order to 
avoid additional notation that would not be needed later.
Besides that, the ideas used
in the above proof (Lemma 5.8 in particular) work for some
nonintegrable representations
of affine Lie algebras, where there is no (known) vertex
operator formula analogous to the
integrable case. Similar ideas have been used in some other
situations, see for example
\cite{FNO}.

Theorem 5.14 was inspired by  \cite{MP, Theorem 9.28} in 
the principal picture, except
that here it is a consequence of the Weyl-Kac character formula.
Moreover, the Kac-Wakimoto
character formula can be used in a similar way for some
modular representations (cf. \cite{Ad}).

Theorems 5.9 and 5.14 put together give Theorem 5.15, which
in terms of the vertex operator
algebra $L(k\Lambda_0)$ states that the level $k$ standard
$\tilde\g$-modules are its only
irreducible highest weight modules. Results of this type
have been proved by different
methods in \cite{FZ}, \cite{DL} and \cite{Li}.

Theorem 5.1 and Proposition 5.2 are similar to some results in
\cite{H} at the critical level.
\endremark

\head{6. Colored partitions, leading terms and the main
results}\endhead

{\it In Sections 6--10 we specialize to the case $\goth g =\goth
s\goth\l (2,\Bbb F)$ and 
level $k \in \Bbb N$.}

\subhead{6.1. Colored partitions}\endsubhead

Let $A$ be a nonempty set and denote by $\Cal P (A)$ the set
of all maps
$\pi : A \rightarrow \Bbb N$, where $\pi(a)$ equals zero for
all but finitely
many $a \in A$. Clearly $\pi$ is determined by its values
$(\pi(a) \mid a \in
A)$ and we shall also write $\pi$ as a monomial
$$
\pi = \prod_{a\in A} a^{\pi(a)}. \tag{6.1.1}
$$

We shall say that $\pi$ is a {\it partition} and for $\pi(a)
> 0$ we shall say that
$a$ is a {\it part} of $\pi$. We define the {\it length}
$\ell(\pi)$ of $\pi$ by
$$
\ell(\pi) = \sum_{a\in A} \pi(a).
$$
For $\rho,\pi \in \Cal P (A)$ we write $\rho \subset \pi$
(or
 $\pi \supset \rho$) if $\rho (a) \le \pi (a)$ for all $a
\in A$ and
 we say that $\rho$ is {\it contained}  in $\pi$. For 
$\rho,\pi \in \Cal P (A)$
we define $\pi \rho$ in $\Cal P (A)$ by
$(\pi\rho)(a)=\pi(a)+\rho(a)$,
$a \in A$. If $\rho \subset \pi$, then we define $\pi /
\rho$ in
$\Cal P (A)$ by $(\pi / \rho)(a)=\pi(a)- \rho(a)$, $a \in
A$. Clearly
$\rho(\pi / \rho)= \pi$.

We also define $\pi \cup \nu$ and $\pi \cap \nu$ by
$$\gather
(\pi \cup \nu)(a) = \max \{\pi(a), \nu(a)\},\\
(\pi \cap \nu)(a) = \min \{\pi(a), \nu(a)\}.
\endgather$$

We shall say that $1=\prod_{a\in A} a^0$ is the partition
with no parts
 and length 0 (and we shall sometimes denote it as
$\varnothing$). Clearly 
$\Cal P (A)$ is a monoid, graded by length. We
shall usually use
the multiplicative notation (6.1.1). We consider elements of
$A$ as partitions
of length 1, i.e. $A \subset \Cal P (A)$.

For lack of a better terminology, we shall say that $\Cal I
\subset \Cal P (A)$
is an {\it ideal} in the monoid $\Cal P (A)$ if $\rho \in
\Cal I$ and $\pi \in \Cal P (A)$
implies $\rho \pi \in \Cal I$. This should not be confused
with a notion of
partition ideal which we shall use in Section 11.

In the case when $A$ is $\Bbb Z, \Bbb Z_{<0}$ or $\Bbb
Z_{\le 0}$ we shall
call elements of $\Cal P (A)$ {\it plain partitions}. For a
plain partition
$\pi \in \Cal P (\Bbb Z)$ we define the {\it degree} $\vert
\pi \vert$ of $\pi$ by
$$
\vert \pi \vert = \sum_{n \in \Bbb Z} n \,\pi(n),
$$
and we say that a part $n$ of $\pi$ has {\it degree} $n$.

Set $x =\left(\smallmatrix 0 & 1\\ 0 & 0\endsmallmatrix
\right)$,
 $h =\left(\smallmatrix 1 & 0\\ 0 & -1\endsmallmatrix
\right)$ and
$y =\left(\smallmatrix 0 & 0\\ 1 & 0\endsmallmatrix
\right)$, so that
$B = \{x,h,y\}$ is the usual basis of $\goth g$. Set
$$\align
&\bar{B} = \{b(n) \mid  b \in B ,  n \in \Bbb Z\},\\
&\bar{B}_- =\bar{B}_{<0} \cup \{y(0)\},\\
&\bar{B}_{<0} = \{b(n) \mid b \in B ,  n \in \Bbb Z_{<0}\},
\endalign$$
so that $\bar{B}, \bar{B}_-, \bar{B}_{<0}$ parameterize
bases of Lie algebras 
$\hat{\goth g}/\Bbb F c$, $\tilde{\goth n}_-$ and
$\tilde{\goth g}_{<0}$
 (respectively).  We shall say that elements of $\Cal P
(\bar{B})$,
 $\Cal P (\bar{B}_-)$ or 
 $\Cal P (\bar{B}_{<0})$ are {\it colored partitions}. For
an element $a=b(n)$ in
$\bar B$ we say that it is of {\it degree} $\vert a\vert =
n$, {\it color} $b$ and 
$\goth h$-{\it weight} $\wt (a) = \wt (b)$, where $\wt (x) =
\alpha$,  $\wt (h) = 0$ and  $\wt (y) = -\alpha$.
 For $\pi \in \Cal P (\bar B)$ we define the
{\it degree} $\vert \pi \vert$ and {\it weight} $\wt (\pi)$
by
$$
\vert \pi \vert = \sum_{a\in \bar B}  \,\pi (a)\,\vert a
\vert, \quad
\wt (\pi) = \sum_{a \in \bar B} \,\pi (a)\,\wt (a) .
$$
We also define the {\it shape} $\sh (\pi)$ of $\pi$ as the
plain partition
$$
\sh (\pi) = \prod_{a \in \bar B} \vert a \vert^{\pi(a)}.
$$

\subheading{6.2. An order on colored partitions}

We choose the order $\preccurlyeq$ on $B$ defined as $y
\prec h \prec x$,
and the order $\preccurlyeq$ on $\bar B$ defined by 
$$
b_1 (j_1)\prec b_2 (j_2) \quad \text{iff} \quad j_1 < j_2
\quad \text{or}\quad
j_1 = j_2, \quad b_1 \prec b_2. \tag{6.2.1}
$$
Then we may think of the colored partition $\pi \in \Cal P
(\bar B)$ of the
 form $\pi = \prod b_i (j_i)$ as 
$$
\pi =(b_1 (j_1),\dots, b_s (j_s)), \quad
 b_1 (j_1)\preccurlyeq\dots\preccurlyeq b_s (j_s),
\tag{6.2.2}
$$
where $b_i (j_i) \in \bar B$ are parts of $\pi$ and $s \ge
0$ is the
length of $\pi$.

We may visualize a colored partition $\pi \in \Cal P (\bar
B_{<0})$ by its
Young diagram consisting of $s$ rows of boxes: in the first
row $-j_1$ boxes of
color $b_1$, in the second row $-j_2$ boxes of color $b_2$,
etc. Similarly
for the shape of $\pi$ --- just without colors:

\def\sq{\lower.3ex\vbox{\hrule\hbox{\vrule height1.2ex
depth1.2ex\kern2.4ex 
\vrule}\hrule}\,}
$$\spreadlines{-1.15ex}\alignat4
\text{color}\quad &b_1\quad
&&\sq\sq&&\!\cdots\sq\!\cdots\sq\sq\quad&&-j_1\quad
\text{boxes}\\
&b_2\quad&&\sq\sq&&\!\cdots\sq\!\cdots\sq\quad&&-j_2\\
&\quad&&&&\!\cdots &&\\
&b_s&&\sq\sq&&\!\cdots\sq &&-j_s \ .
\endalignat
$$

Let $\pi=(a_1,\dots, a_s)$, $\pi'=(a' _1,\dots,a' _{s'})$,
 $\pi,\pi' \in \Cal P (\bar B)$ (i.e. $\pi$ and $\pi'$ are
written in the form (6.2.2)).
We extend the order $\preccurlyeq$ on $\bar B \subset \Cal P
(\bar B)$ defined
by (6.2.1) to the order on $\Cal P (\bar B)$ defined by 
$$
\pi \prec \pi'
$$
if $\pi \neq \pi'$ and one of the following statements hold:
\roster 
\item "{(i)}" $\ell(\pi) > \ell (\pi')$,
\item "{(ii)}" $\ell (\pi) = \ell (\pi')$,
  $\vert \pi \vert < \vert \pi' \vert$,
\item "{(iii)}"  $\ell (\pi) = \ell (\pi')$, 
 $\vert \pi \vert = \vert \pi' \vert$
and there is $i$, $\ell(\pi) \ge i \ge 1$, such that 
$\vert a_j \vert=\vert a' _j\vert$ for $\ell(\pi) \ge j > i$
and
$\vert a_i \vert < \vert a' _i \vert$,
\item"{(iv)}" $\sh \pi=\sh \pi'$ and there is $i$,
$\ell(\pi) \ge i \ge 1$,
such that
$a_j = a' _j$ for $\ell(\pi) \ge j > i$ and $a_i \prec a'
_i$.
\endroster
We also order plain partitions $\Cal P (\Bbb Z)$ by
requirements (i)--(iii).
So if (i), (ii) or (iii) holds, we have that $\sh \pi \prec
\sh \pi'$ and 
(hence)  $\pi \prec \pi'$.

As an example we may take
$$
y(-2)y(-2)x(-2) \prec y(-3)x(-2)y(-1)\prec y(-3)y(-2)x(-1)
\prec x(-4)y(-1)y(-1),
$$
all partitions being of degree $-6$ and weight $-\alpha$.
 
\remark{Remark} In what follows all notions will depend on
the order
chosen. The same results may hold for a different choice of
order as well. 
For example, if we define  $\pi < \pi'$ by  $\pi \neq \pi'$
and
\roster 
\item "{(i)}" $\ell(\pi) > \ell (\pi')$,
\item "{(ii)}" $\ell (\pi) = \ell (\pi')$,
$\vert \pi \vert < \vert \pi' \vert$,
\item "{(iii)}"  $\ell (\pi) = \ell (\pi')$,
 $\vert \pi \vert = \vert \pi' \vert$,
and there is $i$, $\ell(\pi) \ge i \ge 1$, such that 
$a_j=a' _j$ for $\ell(\pi) \ge j > i $
and  $a_i \prec a' _i$,
\endroster
then all the statements in this section hold. However, some
proofs in later
sections cannot be applied as they stand. Note that with
this order we have
$$
y(-2)y(-2)x(-2) < y(-3)x(-2)y(-1)
< x(-4)y(-1)y(-1) < y(-3)y(-2)x(-1).
$$
\endremark
\proclaim{Lemma 6.2.1} (i) The relation $\preccurlyeq$ is a
linear order on
$\Cal P (\bar B)$. The element 1 is the largest element of 
$\Cal P (\bar B)$.

(ii) Let $\ell \ge 0$, $n \in \Bbb Z$ and let 
$S \subset \Cal P (\bar B)$ be a nonempty subset such that
all $\pi$ in $S$
have length $\ell(\pi) \le \ell$ and degree $\vert \pi \vert
\ge n$. Then $S$
has a minimal element.

 (iii) The relation  $\preccurlyeq$ is a (reverse) well
order on
$\Cal P (\bar B_-)$.
\endproclaim
\proclaim {Lemma 6.2.2} Let $\mu, \nu, \pi \in \Cal P (\bar
B)$. Let
$\mu  \preccurlyeq  \nu$. Then $\pi\mu  \preccurlyeq 
\pi\nu$.
\endproclaim

\subheading{6.3. Filtrations on highest weight modules and
leading terms}

For $\pi \in \Cal P (\bar B)$ of the form (6.2.2) set

$$
u(\pi)=b_1 (j_1)\dots b_s (j_s) \in U_k (\tilde{\goth g}).
$$
This defines a map $u : \Cal P (\bar B) \rightarrow U_k
(\tilde{\goth g})$.
For a highest weight $\tilde{\goth g}$-module $V$ with a
highest weight
vector $v_0$ and for $\pi \in \Cal P (\bar B_-)$ define
$$\align
V_{[\pi]} & = \lspan\{ u(\pi')v_0 \mid
\pi' \in \Cal P (\bar B_-), \pi' \succcurlyeq \pi\},\\
V_{(\pi)} & = \lspan\{ u(\pi')v_0 \mid
\pi' \in \Cal P (\bar B_-), \pi' \succ \pi\}\\
& = \bigcup_{\pi' \succ \pi} V_{[\pi']}.
\endalign
$$
Clearly we have
$$\gather
V_{[\pi]} \supset V_{(\pi)} \supset V_{[\pi']} \supset
V_{(\pi')} \quad
\text{for} \quad \pi \prec \pi', \tag{6.3.1}\\
\bigcup_{\pi \in \Cal P(\bar B_-)} V_{[\pi]}=V,\quad
        \bigcap_{\pi \in \Cal P (\bar B_-)}  V_{(\pi)} = 0,
        \quad \dim V_{[\pi]} / V_{(\pi)} \le 1,
\tag{6.3.2}\\
u(\pi') V_{[\pi]} \subset V_{[\pi'\pi]}, \quad
u(\pi') V_{(\pi)} \subset V_{(\pi'\pi)}.
\endgather$$
Moreover, if $\varphi \: V \rightarrow V'$ is a surjective
homomorphism of 
highest weight $\tilde{\goth g}$-modules, then
$$
\varphi(V_{[\pi]})=V' _{[\pi]} \quad , \quad
\varphi(V_{(\pi)})=V' _{(\pi)}.
$$

For $v \in V \backslash \{0\}$ define the {\it leading term}
\ $\lt(v)=\lt_V (v)
\in \Cal P (\bar B_-)$ by
$$
\lt(v)=\pi \qquad \text{iff} \qquad v \in V_{[\pi]}
\backslash V_{(\pi)}.
$$
Because of (6.3.1), (6.3.2) and Lemma 6.2.1(iii) the leading
term $\lt(v)$ is
 well defined for every $v \ne 0$.
For a subset $S \subset V$ set
$$
\lt(S)=\{\lt(v)\mid v \in S \backslash \{0\}\}.
$$
Note that $\lt(M(\Lambda))=\Cal P(\bar B_-)$,
$\lt(N(k\Lambda_0))=
\Cal P(\bar B_{<0})$. Clearly we have

\proclaim{Proposition 6.3.1} (i) For $\pi \in \Cal P(\bar
B_-)$ and 
$v \in M(\Lambda)$ we have $\lt(u(\pi)v)=\pi \lt(v)$. In
particular, for
a  $\tilde{\goth g}$-submodule $V' \subset M(\Lambda)$ we
have that
$\lt_{M(\Lambda)} (V')$ is an ideal in the monoid $\Cal P
(\bar B_-)$.

 (ii) For $\pi \in \Cal P (\bar B_{<0})$ and $v \in
N(k\Lambda_0)$ we 
have $\lt(u(\pi)v)=\pi \lt(v)$. In particular, for a 
 $\tilde{\goth g}$-submodule $V' \subset N(k \Lambda_0)$ we
have that 
$\lt_{N(k \Lambda_0)} (V')$ is an ideal in the monoid
$ \Cal P (\bar B_{<0})$.
\endproclaim

\proclaim{Proposition 6.3.2} If $V'$ is a  $\tilde{\goth
g}$-submodule of a
highest weight  $\tilde{\goth g}$-module $V$, then
$$
\lt(V/V')=\lt(V)\backslash \lt_V (V').
$$
\endproclaim

\proclaim{Proposition 6.3.3} For an $\goth h^e$-invariant
subspace $V' \subset V$ and
$\mu \in \goth h^{e*}$ we have
$$
\dim V' _\mu = \#\{\pi \in \lt_V (V') \mid
\wt(\pi)+\vert\pi\vert \delta=\mu\}.
$$
\endproclaim
\noindent As usual (cf. \cite{K}), here $\delta$ denotes the
imaginary root related to
the homogeneous grading of $\tilde \goth g$.

\remark{Remark}
For a highest weight module $V$ with highest weight
$\Lambda$ we have $V=M(\Lambda)/M$
for some submodule $M$ of the Verma module $M(\Lambda)$. By
Propositions 6.3.1--3 we
can write a character formula
$$
e^{-\Lambda}\ch V=\sum_{\pi \in\Cal P(\bar B_-)\backslash
\lt _{M(\Lambda)}(M)}
e^{\wt(\pi)+\vert\pi\vert\delta}.\tag{6.3.3}
$$
The idea to parameterize bases of standard modules by classes of
partitions originated in the work \cite{LW} on a Lie theoretic
interpretation of the Rogers-Ramanujan identities.
\endremark

\proclaim{Proposition 6.3.4} Let $V$ be a highest weight
$\tilde \goth g$
-module. Then
\roster
\item"{(i)}" $\{u(\pi)v_0 \mid \pi \in \lt(V)\}$
is a basis of $V$.
\item"{(ii)}" For a linear subspace $V' \subset V$ and a map
$u' \: \lt_V (V') \rightarrow V'$ 
such that $\lt_V (u'(\pi))=\pi$ the set

$$
\{u'(\pi) \mid \pi \in \lt_V (V')\} 
$$
is a basis of $V'$.
\endroster
\endproclaim

\subheading{6.4. A filtration on a completion of the
enveloping algebra and
 leading terms}

Recall that the category $\Cal O$ consists of $\tilde\goth g$-modules 
$V$ such that $V$ is a direct sum of finite dimensional weight spaces and such that
the set of weights of $V$ is contained in a finite union of sets of the
form $\lambda -\sum \Bbb N \alpha_i$ (cf. [K]).

Recall that we denote by $\overline{U_k (\tilde \goth g)}$
the completed
enveloping algebra such that  $x_i \rightarrow x$ if for each
$\tilde \goth g$-module 
$V$ in the category $\Cal O$ and each vector $v$ in $V$
there is $i_0$
such that $i \ge i_0$ implies $x_i v = x_{i_0} v = x v$.

\proclaim{Lemma 6.4.1} Let $\pi \in \Cal P (\bar B)$. Then
there exist a
$\tilde{\goth g}$-module $M$ of level $k$ in the category
$\Cal O$, a vector
$v_\pi \in M$ and a functional $v^* _\pi \in M^*$ such that
$$\gather
\langle v^* _\pi, u(\pi)v_\pi \rangle \ne 0,\\
\langle v^* _\pi, u(\pi')v_\pi \rangle = 0 \quad \text{if}
\quad
\pi' \in \Cal P (\bar B), \pi' \ne \pi, \ell(\pi') \le
\ell(\pi).
\endgather$$
\endproclaim
\demo{Proof} Let $\pi = \pi_- \pi_0 \pi_+$, where 
$\pi_- \in \Cal P (\bar B_-)$, $\pi_0= h^q$ for some $q \ge
0$ and the
parts of $\pi_+$ are in $\tilde \goth n_+$. For an element
$b \in \bar B_-$ denote
 by $\bar b$ the element in $\bar B \cap \tilde{\goth n}_+$
such that 
$$
[\bar b, b] \in (\Bbb F h + \Bbb F c)\backslash \{0\}.
$$
Let $p=\ell(\pi_-)$, $r=\ell(\pi_+)$, $\ell=\ell(\pi)=p+q+r$
and let $M$ be
a level $k$ module of the form
$$
M=(\otimes^p _{i=1} M(0))\otimes (\oplus^{q+1} _{i=1}
M(\lambda_i))
\otimes (\otimes^r _{i=1} M(\mu_i)),
$$
where $\lambda_i, \mu_i \in \goth h^{e*}$ are chosen 
so that $\lambda_i (h)\ne \lambda_j (h)$ for $i \ne j$, 
 $\mu_i([\bar b_i,b_i]) > 0$ (i.e\. positive rational numbers)
for all $i=1,\dots,r$, and where $b_i \in \bar B_-$
are such that
$$
\pi_+ =(\bar b_1,\dots, \bar b_r), \quad
 \bar b_1 \preccurlyeq \cdots \preccurlyeq \bar b_r.
\tag{6.4.1}
$$
Let as usual $v_\Lambda \in M(\Lambda)$ denote a highest
weight vector in a Verma
module $M(\Lambda)$. Then the vectors
$$
w_i=\lambda_1 (h)^i v_{\lambda_1}+\cdots+\lambda_{q+1} (h)^i
v_{\lambda_{q+1}},
\qquad i=0,1,\dots,q
$$
form a basis of $\lspan
\{v_{\lambda_1},\dots,v_{\lambda_{q+1}}\}$
and
$$
\{u(\rho)w_i \mid \rho \in \Cal P (\bar B_-), \ 
i=0,\dots,q\}
$$
is a basis of $\oplus^{q+1} _{i=1} M(\lambda_i)$. For other
Verma modules
$M(\Lambda)$ choose the basis $u(\rho) v_\Lambda$ $\rho \in
\Cal P(\bar B_-)$, and for
$M$ the corresponding tensor products. Note that 
$$
h^i w_0 = w_i \ \ \text{for} \ \ i=0,1,\dots, q.
$$

Let $v_\pi$ be a basis vector
$$
v_\pi=(\otimes^p _{i=1} v_0)\otimes w_0 \otimes
(\otimes^r _{i=1} b_i v_{\mu_i})
$$
where $b_i \in \bar B_-$ are defined by (6.4.1).

Let $v^* _\pi$ be the element of dual basis corresponding to
the vector
$$
v=(\otimes^p _{i=1} d_i v_0) \otimes h^q w_0 \otimes
(\otimes^r _{i=1} v_{\mu_i}),
$$
where $d_i \in \bar B_-$ are parts of 
$\pi_- = (d_1,\dots, d_p)$, $d_1 \preccurlyeq \cdots
\preccurlyeq d_p$.
It is clear from our choice that 
$$
v^* _\pi(u(\pi)v_\pi) > 0.
$$

Now let $\pi' \in \Cal P (\bar B)$, $\ell (\pi') \le \ell$,
and
assume that $u(\pi')v_\pi$ has a nontrivial component along
the basis vector 
$v$. Write $\pi'=\pi' _- \pi' _0 \pi' _+$, notation being as
before. In
order to have a component $c v = c(\otimes^p _{i=1} d_i v_0)
\otimes h^q w_0 \otimes
(\otimes^r _{i=1} v_{\mu_i})$, $c \ne 0$, $\pi' _-$ must
contain all parts $d_i$ of $\pi$ and
$\ell(\pi' _0)$ must be at least $q$. Hence $\ell(\pi' _+)
\le r$. But to
have  $u(\pi' _+)(\otimes^r _{i=1} b_i v_{\mu_i})=
c_1(\otimes^r _{i=1} v_{\mu_i})+ \cdots$, $c_1 \ne 0$, we
need 
$\ell(\pi' _+) \ge r$. Hence $\ell(\pi' _-)= r$, $\pi' _0 =
\pi_0$, 
$\ell(\pi' _+) = r$. But then $\pi' _- = \pi_-$ and  $\pi'
_+ = \pi_+$,
i.e.  $\pi' = \pi$. \qed
\enddemo

\remark{Remark} The above proof is a slight refinement of
the proof of 
Lemma 4.1 in \cite{MP}.
\endremark
For $\pi \in \Cal P (\bar B)$ set
$$
U_{[\pi]} =\overline{\lspan\{u(\pi')\mid \pi' \succcurlyeq
\pi\}},
$$
$$
U_{(\pi)} =\overline{\lspan\{u(\pi')\mid \pi' \succ \pi\}},
$$
the closure taken in $\overline{U_k (\tilde \goth g)}$.

\proclaim{Lemma 6.4.2} For $\pi \in \Cal P (\bar B)$ we 
have $U_{[\pi]}=\Bbb F u(\pi)+ U_{(\pi)}$. Moreover,
$$
\dim U_{[\pi]}/ U_{(\pi)}= 1.
$$
\endproclaim
\demo{Proof} Clearly $\Bbb F u(\pi) + U_{(\pi)} \subset
U_{[\pi]}$. We shall
prove the converse inclusion as well: 
\newline Let $u = \lim_{i \in I} u_i \in U_{[\pi]}$. Choose
$v_\pi \in M$ and
 $v^* _\pi \in M^*$ as in Lemma 6.4.1. Then for some $i_0
\in I$ we have
$$\gather
u_i v_\pi = u_{i_0} v_\pi \quad \text{for} \quad i \ge
i_0,\\
v^* _\pi(u_i v_\pi)= v^* _\pi(u_{i_0} v_\pi)=c \quad
\text{for} \quad
 i \ge i_0.
\endgather$$
Since
$$
u_i = c_i u(\pi)+\sum_{\pi'\succ \pi} c^i _{\pi'} u(\pi'),
$$
by Lemma 6.4.1 we have $c_i = c$ for $i \ge i_0$. Hence
$$
u' _i= u_i - c u(\pi) \in U_{(\pi)} \quad \text{for} \quad 
i \ge i_0,
$$
the limit $\lim_{i \in I} u' _i$ exists and 
$u - c u(\pi) \in U_{(\pi)}$.
Now $U_{[\pi]}=\Bbb F u(\pi) + U_{(\pi)}$ implies
$\dim U_{[\pi]}/ U_{(\pi)} \le 1$.
By Lemma 6.4.1 \ $\langle v^* _\pi, U_{(\pi)} v_\pi \rangle
= 0$,
$\langle v^* _\pi, U_{[\pi]} v_\pi \rangle \ne 0$,
and the lemma follows.\qed
\enddemo

Set
$$
U_k(\tilde{\goth g})_{\text{loc}}=\lspan\{v_n \mid
v \in N(k\Lambda_0), n \in \Bbb Z\},
$$
where $v_n$ denotes a coefficient in $Y(v,z)$. From the
commutator 
formula (2.15) we see that
$U_k(\tilde{\goth g})_{\text{loc}}$ is a Lie algebra. From
the normal order product
formula (2.14) we see that $U_k(\tilde{\goth
g})_{\text{loc}} \subset 
\bigcup_{\pi \in \Cal P(\bar B)}\  U_{[\pi]}$.
Let us denote by $U'(U_k(\tilde{\goth g})_{\text{loc}})$ the
associative subalgebra
of $\overline{U_k (\tilde \goth g)}$ generated by
$U_k(\tilde{\goth g})_{\text{loc}}$.  
Then we have: 
$$
U'(U_k(\tilde{\goth g})_{\text{loc}}) \subset \bigcup_{\pi
\in \Cal P(\bar B)}
U_{[\pi]}.
$$

For $u \in U_{[\pi]}$, $u \notin U_{(\pi)}$ we define the
{\it leading term}
$$
\lt(u) = \pi.
$$
\remark{Remark} From Lemma 6.4.1 we see that if the leading
term $\lt(u)$ of
$u$ in $U'(U_k(\tilde{\goth g})_{\text{loc}})$ exists, it is
unique.
\endremark

\proclaim{Proposition 6.4.3} Every element
 $u \in U'(U_k(\tilde{\goth g})_{\text{loc}})$, 
$u \neq 0$, has a unique leading
 term $\lt(u)$. Moreover, for any highest weight module $V$
we have that
$$
u V_{[\pi]} \subset V_{[\pi \lt(u)]},\qquad 
u V_{(\pi)} \subset V_{(\pi \lt(u))}.
$$
\endproclaim

\demo{Proof}  We should prove
that each element in $A= U'(U_k (\tilde{\goth
g})_{\text{loc}})$ has a leading term.

First note that $A$ is graded by degree, say $A = \coprod_{n
\in \Bbb Z} A_n$.

Recall that by the definition of convergence in
 $\overline{U_k (\tilde{\goth g})}$ we have $x_i \rightarrow
x$ if for each
module $V$ in the category $\Cal O$ and each vector $v$ in
$V$ there is $i_0$
such that $i \ge i_0$ implies $x_i v = x_{i_0} v = x v$. In
this proof
we shall use the notion of {\it uniform convergence}: we say
that
 $x_i \rightarrow x$ uniformly if for each $N \in \Bbb N$
there is $i_0$
such that for each module $V$ in the category $\Cal O$ and
each $v$ in
$V$, $|\deg v| < N$, $i \ge i_0$ implies $x_i v = x_{i_0} v
= x v$.
Here by $|\deg v | < N$ we mean that for a weight vector $v
\in V_\mu$ the
degree of any maximal weight of $V$ differs from the degree
of $\mu$ 
by at most $N$.

 From the normal order product formula (2.13) we see by
induction that each $x \in A_n$
is a limit of a sequence $(x_i)_{i \in \Bbb N}$ such that
$x_i \in U_k (\tilde{\goth g})_{\ell,n}$ (where $\ell$
denotes a filtration
and $n$ a degree) and that $x_i \rightarrow x$ uniformly.

Set $\Cal P_{\ell,n} = \{\pi \mid \ell(\pi) = \ell$, $|\pi |
= n \}$. Then
$\Cal P_{\ell_1, n} \prec \Cal P_{\ell_2, n}$ for $\ell_1 >
\ell_2$, and each
$\Cal P_{\ell,n} \cong \Bbb N$ as an ordered set.

Now we show in three steps that $x \in A_n \cap U_{[\pi]}$
is the uniform
limit of a sequence $(x_i)$, $x_i \in U_k (\tilde{\goth
g})_{[\pi]} =
\lspan \{u(\pi') \mid \pi' \succcurlyeq \pi \}$, 
$\deg x_i = n$:

(1) Let $\pi_0$ be the smallest element in $\Cal
P_{\ell,n}$, 
$\ell = \ell(\pi)$, $n = |\pi|$. We (may) assume that 
$x \in U_{[\pi]} \cap A_n$ is a uniform limit of a sequence
$(x_i)$,
$x_i \in U_k (\tilde{\goth g})_{\ell,n}$. Note that the
segment 
$[\pi_0, \pi]$ is finite, and write
$$
x_i = \sum_{\pi_0 \preccurlyeq \kappa \prec \pi} C_{\kappa,
i} u(\kappa) +
\sum_{\nu \succcurlyeq \pi} C_{\nu, i} u(\nu).
$$
 
For $\kappa \prec \pi$ take $v_\kappa , v_\kappa ^*$ as in
Lemma 6.4.1. Then $x \in U_{[\pi]}$
implies $v_\kappa ^*(x v_\kappa) = 0$. Now $x_i \rightarrow
x$ implies that
there is $i_0$ such that $x_i v_\kappa = x_{i_0} v_\kappa =
x v_\kappa$
for $i \ge i_0$, and hence $C_{\kappa, i} = 0$ for $i \ge
i_0$. So
in a finite number of steps we see that for some $i_0$ the
sequence 
$(x_i)_{i \ge i_0}$, $x_i \in U_k (\tilde{\goth
g})_{[\pi]}$,
converges uniformly to $x$.

(2) Let $\ell = \ell(\pi) = \ell(\pi')$, $n = |\pi| =
|\pi'|$,
$\pi \prec \pi'$ and $x \in U_{[\pi]}$, $x \notin
U_{(\pi')}$. Since the
interval $[\pi,\pi']$ is finite, there is
 $\pi \preccurlyeq \tau \preccurlyeq \pi'$ such that
 $x \in U_{[\tau]} \backslash U_{(\tau)}$, i.e. $u$ has a
leading term.

(3) Let $x \in A_n$, $x \in  U_{[\pi]}$ for all $\pi$ such
that 
$\ell (\pi) = \ell$. For each $\pi$ we may choose a sequence
$(x_i ^\pi)$ 
in $U_k (\tilde{\goth g})_{[\pi]}$ such that
$x_i ^\pi \rightarrow x$ uniformly. Let
$$
\align
x_i ^\pi & = \sum_{\ell(\kappa)=\ell} C^\pi _{\kappa ,i}
u(\kappa) +
x^\pi _{i,\ell - 1}\\
& = x^\pi _{i,\ell} + x^\pi _{i,\ell - 1},
\endalign
$$
where $x^\pi _{i,\ell - 1} \in  U_k (\tilde{\goth g})_{\ell
- 1,n}$. For a
sequence $\pi^{(1)} \prec \pi^{(2)} \prec \cdots$ such that
$$
U_{[\pi^{(1)}]} \supset  U_{[\pi^{(2)}]} \supset \cdots
\supset
\bigcap_{\ell(\pi)=\ell} U_{[\pi]} \cap A_n = 
\bigcup_{\ell(\pi) < \ell} U_{[\pi]} \cap A_n
$$
we can construct a diagonal sequence
$$
x^{\pi^{(k)}} _{i_k} \ \ , \ \ k=1,2,\dots
$$
such that (a) $i_k$ is such that for each module $V$ and
each vector $v$
in $V$, $|\deg v | < k$, the inequality $i \ge i_k$ implies 
$x^{\pi^{(k)}} _i v = x^{\pi^{(k)}} _{i_k} v = x v$, and (b)
$k < k'$
implies $i_k < i_{k'}$. Then $x^{\pi^{(k)}} _{i_k}
\rightarrow x$ uniformly.
Since $x^{\pi^{(k)}} _{i_k} \in U_{[\pi^{(k)}]}$, we have
 $x^{\pi^{(k)}} _{i_k,\ell} \rightarrow 0$ and we may take a
subsequence 
$x_{i_{k_m},\ell} ^{\pi^{(k_m)}}$ which converges to $0$
uniformly. Hence we
have constructed a sequence
$$
x_m = x_{i_{k_m},\ell -1} ^{\pi^{(k_m)}} \in
\bigcup_{\ell(\kappa)<\ell} U_{[\kappa]} \cap
A_n 
$$
such that $x_m \rightarrow x$ uniformly.

Now in a finite number of steps we see that either $x=0$ or
$x$ has a leading
term. \qed
\enddemo

\remark{Remark} The above proof shows the following:
Let $U_k (\tilde{\goth g}) \subset A
\subset\overline{U_k(\tilde{\goth g})}$
be a graded subspace. Then every $x$ in $A$ has a leading
term if every $x$ in $A$ is a 
uniform limit of a sequence in $U_k(\tilde{\goth g})_\ell$
for some fixed $\ell$
(depending on $x$).
\endremark

\subheading{6.5. The main results}

{\it In the remainder of this section we shall assume that
$\Lambda$ is 
dominant integral and that $k=\Lambda(c)$.}

First let us recall the construction of the loop module $\bar R$: Since we 
assume that $\goth g =\goth s\goth\l (2,\Bbb F)$, we have $\theta=\alpha$ if
we write
$$
x =x_\alpha= \left(\matrix 0& 1 \\
0& 0 \endmatrix\right),
\quad y =x_{-\alpha}= \left(\matrix 0& 0 \\
1& 0 \endmatrix\right),
\quad h =\alpha^\vee= \left(\matrix 1& 0 \\
0& -1 \endmatrix\right).
$$
Then $R$ defined by (5.4), i.e.
$$
R = U(\g)x_{\alpha}(-1)^{k+1}\1\subset N(k{\Lambda}_0),
$$
is a $(2k+3)$-dimensional $\goth s\goth\l (2,\Bbb F)$-module 
and $\bar R$ is defined as a linear span of 
the coefficients of the fields
$$
Y(r,z) = \sum_{n\in\Z}r(n)z^{-n-k-1}
$$
for $r\in R$. As it has been remarked, these coefficients $r(n)$ are infinite
sums of elements in $U_k(\tilde{\goth g})$, but are well defined operators
on every highest weight $\tilde{\goth g}$-module $V$ of level $k$.

Note that every vector $v\neq 0$ in a Verma $\tilde{\goth g}$-module $M(\Lambda)$
can be written as a linear combination of terms $u(\pi)v_\Lambda$, where
$$
u(\pi)=b_1 (j_1)\dots b_s (j_s) \in U_k (\tilde{\goth g})\quad\text{for}\quad
\pi =(b_1 (j_1)\preccurlyeq\dots\preccurlyeq b_s (j_s))\in \Cal P (\bar B_-)
$$
and $v_\Lambda$ is a highest weight vector,
and then $\lt (v)$ is just the smallest $\pi$ which appears nontrivially in
this linear combination. The existence of the leading term $\rho=\lt (r(n))$ is a bit
more subtle: it is guaranteed by Proposition 6.4.3, but as we shall see
later on, $r(n)\in\bar R$ can be written in the form
$$
   r(n)=c_\rho u(\rho)+\sum_{\pi\succ \rho}\,c_\pi u(\pi),
$$
where $c_\rho, c_\pi \in \Bbb C$, $c_\rho \neq 0$, and again the leading term is
the smallest colored partition which appears nontrivially in this (infinite)
linear combination.  

For a subset $S \subset \Cal P(A)$ denote by $(S) \subset
\Cal P(A)$ the
ideal in the monoid $\Cal P (A)$ generated by $S$. 
To avoid ambiguities because of notation, let us agree with
a convention that
$$\align
& \lt (\bar R) \subset \lt (U'(U_k(\tilde {\goth
g})_{\text{loc}})) = \Cal P (\bar B), \\
&\lt _{M(\Lambda)}(\bar R v_\Lambda)\subset \lt
_{M(\Lambda)}(M(\Lambda)) = 
      \Cal P (\bar B_-),\\
&\lt _{N(k\Lambda_0)}(\bar R \bold 1)\subset \lt
_{N(k\Lambda_0)}(N(k\Lambda_0)) = 
      \Cal P (\bar B_{<0}),
\endalign$$
so that $(\lt (\bar R)) \subset \Cal P (\bar B)$, $(\lt
(\bar R v_\Lambda))\subset 
\Cal P (\bar B_-)$, $(\lt (\bar R \bold 1))\subset \Cal P
(\bar B_{<0})$.

Let $U_0 \subset U_1 \subset U_2 \subset \cdots$ be the
filtration of
$U_k(\tilde{\goth g})$. Since by Theorem
5.14  $\bar RM({\Lambda})$ is the maximal submodule $M^1({\Lambda})$ of a
Verma module $M({\Lambda})$, we have:

\proclaim{Lemma 6.5.1} Let $v_\Lambda$ be a highest weight
vector of Verma module
$M(\Lambda)$. Then 
$$\gather
\lt(\bar R v_\Lambda)\subset \lt(U_1 \bar R
v_\Lambda)\subset \cdots
\subset \lt(U_j \bar R v_\Lambda)\subset \cdots \subset 
\lt_{M(\Lambda)} (M^1(\Lambda)), \tag"{(i)}"\\
(\lt(\bar R v_\Lambda))\subset (\lt(U_1 \bar R
v_\Lambda))\subset \cdots
\subset (\lt(U_j \bar R v_\Lambda))\subset \cdots \subset 
\lt_{M(\Lambda)} (M^1(\Lambda)), \tag"{(ii)}"\\
\bigcup_{j\ge 0} (\lt(U_j \bar R v_\Lambda))= 
\lt_{M(\Lambda)} (M^1(\Lambda)). \tag"{(iii)}"
\endgather$$
\endproclaim

We think of elements of $U_j\bar R$ as relations that hold
on
$L(\Lambda)=M(\Lambda)/M^1(\Lambda)$. Note that (iii)
shows that relations $U_j \bar R$, $j \ge 0$, determine 
the character of $L(\Lambda)$ written in the form (6.3.3), i.e
$$
e^{-\Lambda}\ch L(\Lambda)=\sum_{\pi \in\Cal P(\bar B_-)\backslash
\lt _{M(\Lambda)}(M^1(\Lambda))}
e^{\wt(\pi)+\vert\pi\vert\delta}.
$$
In what follows we shall describe the ideal 
$\lt_{M(\Lambda)} (M^1(\Lambda))$ by explicitly constructing a basis
of $M^1({\Lambda})=\bar RM({\Lambda})$ parameterized by colored partitions.

We shall write $\lt$ instead of $\lt _{M(\Lambda)}$ when no confusion can arise.
Clearly there exists a map
$$
\lt(\bar R v_\Lambda)\rightarrow \bar R, \qquad \rho \mapsto
r(\rho),
$$
such that $\lt(r(\rho)v_\Lambda)=\rho$, i.e. for each colored partition
$\rho\in\lt(\bar R v_\Lambda)$ we choose a relation $r(\rho)$ in $\bar R$ which
will produce a vector $r(\rho)v_\Lambda$ in $M^1(\Lambda)$ with the 
leading term $\rho$. With given such a map we define
$$
u(\rho \subset \pi)=u(\pi / \rho)r(\rho)
$$
for $\rho \in \lt(\bar R v_\Lambda)$ and $\pi \in (\lt(\bar
R v_\Lambda))$, $\rho \subset \pi$.
Since our order $\preccurlyeq$ behaves well with respect to multiplication,
we have $\lt(u(\rho \subset \pi)v_\Lambda) =\pi$.
For each $\pi$ in the ideal $(\lt(\bar R v_\Lambda))$ there is at least
one
$\rho \in \lt(\bar R v_\Lambda)$ such that $\rho \subset
\pi$. Hence there
is a map
$$
(\lt(\bar R v_\Lambda)) \rightarrow  \lt(\bar R v_\Lambda),
\qquad
\pi \mapsto \rho(\pi),
$$
such that $\rho(\pi) \subset \pi$. 

In other words, for each colored 
partition $\pi$ in the ideal $(\lt(\bar R v_\Lambda))$ we can choose
$\rho=\rho(\pi)$ such that $\pi=(\pi/\rho)\cdot\rho$. Moreover, with this 
choice of $\rho$ we can construct a vector 
$u(\rho \subset \pi)v_\Lambda=u(\pi / \rho)r(\rho)v_\Lambda$ in 
$M^1(\Lambda)$ with the leading term $\pi$. As we shall see later on,
the main difficulty in our construction is that for a given $\pi$ this can be done in 
several ways. For example, we may have $\pi=x(-3)x(-2)x(-1)$ and vectors
$x(-3)r(x(-2)x(-1))v_\Lambda$ and $x(-1)r(x(-3)x(-2))v_\Lambda$ in
$M^1(\Lambda)$ with the same leading term $\pi$.
Our main result for $\goth g = \text{sl}(2,\Bbb F)$ is 
that for each $\pi$ we need only one such vector in order to obtain 
a basis of $M^1(\Lambda)$:

\proclaim{Theorem 6.5.2} For any fixed map 
$(\lt(\bar R v_\Lambda))\ni\pi \mapsto \rho(\pi)\in\lt(\bar
R v_\Lambda)$ such that
$\rho(\pi) \subset \pi$ the set of vectors
$$
u(\rho(\pi) \subset \pi)v_\Lambda, \qquad
 \pi \in (\lt(\bar R v_\Lambda)), \tag{6.5.1}
$$
is a basis of the maximal submodule $M^1(\Lambda)$ of the Verma
module $M(\Lambda)$.
\endproclaim

\remark{Remark}
The above theorem is analogous to [MP, Theorem 9.26] in the principal picture. 
Although in \cite{MP} a notion of the leading term was not introduced, the
general ideas are quite similar to the ones used here. At this point let
us only mention that the order used in \cite{MP} on the set of partitions
$\Cal P$ slightly differs from our present order on $\Cal P(\Bbb Z)$,
but in both cases the minimal elements used in our arguments, in
\cite{MP} denoted by (p;n), are the same. So in the present 
terminology Lemma 9.9 in \cite{MP} defines the leading term $X(p;n)$ of a 
relation $R_p(n)$, and $X(\mu)$ in Theorem 9.26 
in \cite{MP} is analogous to $u(\rho(\pi) \subset \pi)$ for some specific
choice of the map $\pi \mapsto \rho(\pi)$.
\endremark

Theorem 6.5.2 will be proved in Section 10 as a consequence
of Theorem 10.1, 
and almost all of the rest of this paper is devoted to the
proof. At this point
we state the following immediate consequences of Theorem
6.5.2 and
Proposition 6.3.2:

\proclaim{Corollary 6.5.3} $\lt_{M(\Lambda)}(M^1(\Lambda))=
(\lt(\bar R v_\Lambda))$.
\endproclaim
\proclaim{Corollary 6.5.4} $\lt(L(\Lambda))=
\Cal P (\bar B_-) \backslash (\lt(\bar R v_\Lambda))$.
\endproclaim

\proclaim{Theorem 6.5.5} The set of vectors
$$
u(\pi)v_\Lambda, \qquad \pi \in 
\Cal P (\bar B_-) \backslash (\lt(\bar R
v_\Lambda)),\tag{6.5.2}
$$
is a basis of the standard $\tilde{\goth g}$-module
$L(\Lambda)$.
\endproclaim
\demo{Proof}
Since on $L(\Lambda)$ relations $r(\rho)$, $\rho \in \lt
(\bar R)$, vanish, we
see by induction that (6.5.2) is a spanning set of
$L(\Lambda)$. 
But then by Corollary 6.5.4 and Proposition 6.3.3 it must be
a basis.
\qed \enddemo

\remark{Remark}
Corollary 6.5.3 states that all inclusions in Lemma
6.5.1(ii) are equalities.
In the case of $\goth g = \goth {sl} (3,\Bbb F)$ and
$\Lambda = \Lambda_0$
take $B$ to be the ordered basis
$$
[f_1,f_2]\prec f_1\prec f_2\prec h_2\prec h_1\prec e_2\prec
e_1\prec [e_1,e_2],
$$
where $e_i,h_i,f_i$ denote the Chevalley generators of
$\goth g$, and let
the order on $\Cal P (\bar B)$  be 
defined as above. Then one can see that 
$(\lt(\bar R v_{\Lambda_0}))\subsetneqq (\lt(U_1 \bar R
v_{\Lambda_0}))$, and some
evidence suggest that one should expect equalities
$$
 (\lt(U_1 \bar R v_{\Lambda_0})) = \cdots
= (\lt(U_j \bar R v_{\Lambda_0})) = \cdots =
\lt_{M(\Lambda_0)} (M^1(\Lambda_0)).
$$
\endremark

\subheading{6.6. Difference and initial conditions}

Now we will describe the sets $\lt (\bar R)$, $\lt (\bar R
v_\Lambda)$ 
and $\lt (\bar R \bold 1)$:
In Section 5 we denoted by
$$
R = U(\goth g) x(-1)^{k+1} \bold 1 \subset N(k\Lambda_0)
$$
and by $\bar R$ the corresponding loop module. $R$ is a
$(2k+3)$-dimensional 
$\goth g$-module with a basis $\{r_{i\alpha}\mid -k-1 \le i \le
k+1\}$, where
$$
r_{i\alpha} = \tfrac{1} {(k+1+i)!} \,x(0)^{k+1+i} \cdot
y(-1)^{k+1} \bold 1. \tag{6.6.1}
$$
In particular $r_{-(k+1)\alpha}=y(-1)^{k+1} \bold 1$, and
$$\gather
Y(r_{-(k+1)\alpha} , z) = \sum_{n\in \Bbb Z} r_{-(k+1)\alpha} (n)
z^{-n-k-1},\\
r_{-(k+1)\alpha} (n) = \sum_{j_1+\cdots +j_{k+1}=n}  y(j_1) \dots
y(j_{k+1}).
\endgather$$
Clearly $r_{i\alpha}(n) \bold 1 \ne 0$ for $n \le - k - 1$. Set
$\pi_{i\alpha} (n)=\lt(r_{i\alpha} (n))$. It is also clear that for
$n=(k+1)j - s$,
$0 \le s \le k+1$,
$$
\sh (\pi_{-(k+1)\alpha} (n)) = (j - 1)^s j^{k+1 - s}. \tag{6.6.2}
$$
By the adjoint action of $x(0)$ we obtain from $r_{-(k+1)\alpha}
(n)$ all $r_{i\alpha} (n)$,
and from (6.6.2) we see that 
$$
\sh (\pi_{i\alpha} (n)) = (j-1)^s j^{k+1-s}.
$$
It is easy to see the smallest possible colorings of the
shape
$(j-1)^s j^{k+1-s}$ which appears in $r_{i\alpha}(n)$: For example,
for $k=4$ and the
shape $(-2)^2 (-1)^3$ the elements $\pi_{i\alpha} (-7)$, $i=-5,\dots ,5$, are
$$\spreadlines{-1.15ex}\alignat 2
&\sq\sq \quad &&y\ y\ h\ h\ x\ x\ x\ x\ x\ x\ x  \\
&\sq\sq \quad &&y\ h\ h\ x\ x\ x\ x\ x\ x\ x\ x  \\
&\sq \quad    &&y\ y\ y\ y\ y\ y\ y\ h\ h\ h\ x  \\
&\sq \quad    &&y\ y\ y\ y\ y\ y\ h\ h\ h\ x\ x  \\
&\sq \quad    &&y\ y\ y\ y\ y\ h\ h\ h\ x\ x\ x  \ .
\endalignat
$$
In general we have:

\proclaim{Proposition 6.6.1} The set $\lt(\bar R )$ consists
of
elements 
$\pi_{m\alpha} (n)=\lt(r_{m\alpha} (n))$, $n \in \Bbb Z$,
of the shape $(j-1)^a j^b$, $a(j-1)+ bj = n$, $a+b=k+1$,
 $0 \le a,b \le k+1$,
$$
\pi_{m\alpha} (n) = \left\{
\aligned y(j-1)^r h(j-1)^{a-r} y(j)^b & \qquad 
a \ge r \ge 0,  \ m=-k-1+a-r,\\
 h(j-1)^r x(j-1)^{a-r} y(j)^b & \qquad 
a > r \ge 0,  \ m=-k-1+2a-r,\\
 x(j-1)^a y(j)^r h(j)^{b-r} & \qquad 
b \ge r \ge 0, \  m=k+1-b-r,\\
 x(j-1)^a h(j)^r x(j)^{b-r} & \qquad 
b > r \ge 0,  \ m=k+1-r\ .\endaligned
\right.$$
The set $\lt(\bar R \bold 1)$ consists of
the elements 
$\pi_{m\alpha} (n)$, $n \le -k - 1$, $m = -k-1, \dots, k+1$.
\endproclaim

\remark{Remark}
Later on we shall use the map
$$
   \lt(\bar R) \rightarrow \bar R,\qquad \rho \mapsto
r(\rho),
$$
defined by $r(\pi_{m\alpha}(n)) = r_{m\alpha}(n)$ for $n \in \Bbb Z$, $|m|
\le k+1$.
Clearly $\lt (r(\rho)) = \rho$.
\endremark

In the case of Verma module $M(\Lambda)$ the same argument
as above applies for leading
terms of $r_{i\alpha}(n) v_\Lambda$, $n \le -k - 1$. We may state
this as the 
following:

\proclaim{Proposition 6.6.2} For $n \le - k - 1$,
$i=-k-1,\dots, k+1$,
we have
$$ 
\lt_{M(\Lambda)}(r_{i\alpha} (n) v_\Lambda) =
\lt_{N(k\Lambda_0)}(r_{i\alpha}(n)\bold 1) = \pi_{i\alpha} (n).
$$
\endproclaim
Hence the ideals $(\lt(\bar R v_\Lambda))$ have in common 
(for different $\Lambda$'s) the generators
$\pi_{i\alpha} (n)$, $n \le - k - 1$. 
These generators (colored partitions) have ``very small"
differences between the
largest and the smallest part. Since a partition $\pi$ from 
$\Cal P (\bar B_-) \backslash (\lt(\bar R v_\Lambda))$ does
not contain any 
of these generators, we may say that $\pi$ satisfies
difference conditions
defined by $\lt(\bar R \bold 1)$. In this vein we may also
say that the set
$\lt(\bar R v_\Lambda) \backslash \lt(\bar R \bold 1)$
represents the
initial conditions which $\pi$ from $\Cal P(\bar B_-)
\backslash
 (\lt(\bar R v_\Lambda))$ satisfies (since these are the
conditions on the
``initial part" of colored partition).

Write $\Lambda = k_0 \Lambda_0+k_1 \Lambda_1$, $k_0+k_1=k$, 
$k_0, k_1 \in \Bbb N$.

\proclaim{Proposition 6.6.3} The set 
$\lt(\bar R v_\Lambda) \backslash \lt(\bar R \bold 1)$
consists of
elements denoted by \linebreak
 $\pi_{m\alpha}^\Lambda (n)=\lt_{M(\Lambda)}(r_{m\alpha} (n)v_\Lambda)$,
where $- k - 1 < n \le 0$ and

$$
\pi_{m\alpha} ^\Lambda(n) = \left\{
\aligned y(-1)^r h(-1)^{a-r} y^b & \qquad 
a+b=k+1, \ a \ge r \ge 0,\\
& \qquad m=-r-b, \ n=-a,\\
h(-1)^r x(-1)^{a-r} y^b & \qquad 
a+b=k+1, \ a > r \ge 0, \\
& \qquad m=-b+a-r, \ n=-a,\\
x(-1)^a y^r & \qquad 
a+r<k+1, \ a > k_0 \quad\text{or}\quad r > k_1,\\
& \qquad m=a-r, \ n=-a.\\
\endaligned
\right.$$
\endproclaim

Proposition 6.6.3 will be proved in Section 10, see
Proposition 10.2.

\proclaim{Example} For $\Lambda$ of level $1$ we have the
difference conditions:

$$\spreadlines{-1.15ex}\alignedat 2
\cdots\ &\sq \quad    &&y\ y\ h\ h\ x  \\
\cdots\ &\sq \quad    &&y\ h\ h\ x\ x \ , 
\endalignedat\qquad\
\alignedat 2
\cdots\ &\sq\sq \quad &&y\ h\ x\ x\ x \\
\cdots\ &\sq \quad    &&y\ y\ y\ h\ x \ . 
\endalignedat\tag{6.6.3}
$$

For $\Lambda = \Lambda_0$ the initial conditions are
 $y^2, y, y(-1)y, h(-1)y, x(-1)y$.

For $\Lambda = \Lambda_1$ the initial conditions are 
 $y^2, y(-1)y, h(-1)y, x(-1)y, x(-1)$.

In particular, Theorem 6.5.5 states that for the basic
$\tilde\goth g$-module $L(\Lambda_0)$
we have a basis consisting of vectors of the form
$$
u(\pi)v_{\Lambda_0},
$$
where $\pi \in \Cal P (\bar B_{<0})$ satisfies the
difference condition defined
by $\lt(\bar R \bold 1)$, i.e. $\pi$ does not contain any of
the partitions
listed in (6.6.3).
\endproclaim

\remark{Remark}
As we have already mentioned, a notion of the leading term was 
not introduced in \cite{MP}, but still we can consider $X(p;n)$ to be the leading term 
of the relation $R_p(n)$ in the principal picture. In the
case of level 3 modules the corresponding Young diagrams for $X(2;n)$
have only one ``color'' $X$ and can be identified with

$$\spreadlines{-1.15ex}\aligned
\cdots\ &\sq \quad    \\
\cdots\ &\sq \      \ \ \text{or} 
\endaligned\qquad\
\aligned
\cdots\ &\sq\sq \\
\cdots\ &\sq \quad \    \ . 
\endaligned
$$
Partitions which do not contain any of these leading terms are precisely
the partitions satisfying Rogers-Ramanujan's difference 2 conditions.

Note that in Theorem 6.5.5 as well as in [MP, Theorem 8.7] it is easy to 
show the spanning result for the set of vectors parameterized by
partitions which satisfy difference and initial conditions, i.e. partitions
which do not contain the leading terms of relations which vanish on a given
standard module. It is the linear independence which is difficult to prove,
and that involves, one way or the other, the combinatorial identities of
Rogers-Ramanujan type.

At this point it may be interesting to note that
the combinatorial identities for
the $(1,2)$-specialization of characters of level one
$\goth {sl}(2,\C)\sptilde$-modules given by Theorem 6.5.5 are identical with
the combinatorial identities obtained by S. Capparelli for level 3 modules in 
the principal picture for the type $A^{(2)}_2$ (cf. Section 11.3).
\endremark

\heading{7. Colored partitions allowing at least two
embeddings}
\endheading

In Subsection 6.6 we described the set $\lt (\bar R)$, and, for
example, for level 1 modules we have elements $\rho_1=x(-3)x(-2)$ and 
$\rho_2=x(-2)x(-1)$ in $\lt (\bar R)$. For colored partition
$\pi=x(-3)x(-2)x(-1)$ we can construct two vectors
$u(\rho_1 \subset \pi)v_\Lambda=x(-1)r(x(-3)x(-2))v_\Lambda$ and
$u(\rho_2 \subset \pi)v_\Lambda=x(-3)r(x(-2)x(-1))v_\Lambda$,
and by Theorem 6.5.2 it is enough to take just one of them
for building a basis of the maximal submodule $M^1(\Lambda)$ of the Verma
module $M(\Lambda)$. To prove Theorem 6.5.2 we shall need relations which
relate all possible choices $u(\rho_1 \subset \pi)v_\Lambda$ and
$u(\rho_2 \subset \pi)v_\Lambda$. In this section
we classify all possible choices for $\rho_1 \subset \pi$ and $\rho_2 \subset \pi$.

\subheading{7.1. Embeddings}

Let $\pi \in \Cal P(\bar B)$, $\rho \in \lt(\bar R)$ and let
$\rho \subset \pi$. Then we shall say
that $\rho$ is {\it embedded} in $\pi$, or that $\rho
\subset \pi$ is an {\it embedding}.

If $\rho \subset \pi$ is an embedding and $\kappa \in \Cal
P(\bar B)$, then
$\rho \subset \kappa\pi$ is an embedding as well. For $\rho,
\rho' \in \lt(\bar R)$ 
we have $\rho, \rho' \subset \rho \cup \rho'$.
Since for $\rho \in \lt(\bar R)$ we have $\ell(\rho)= k+1$, 
it is clear from Proposition 6.6.1 that:

\proclaim{Lemma 7.1.1} Let  $\pi =\rho \cup \rho'$. 
Then either $\rho \cap \rho' = \varnothing$, $\ell(\pi)=2k +
2$ and $\pi = \rho\rho'$ or 
$k+2 \le \ell(\pi) \le 2k + 1$ and the shape $\sh(\pi)$ has
a form
$$
(j - 1)^a j^b (j + 1)^c
$$
for some $j \in \Bbb Z$, $b > 0$, $a, c \ge 0$,
$a+b+c=\ell(\pi)$.
\endproclaim

 For $b \in \bar B$ define $b^* \in \bar B$ by 
$$
x(j)^*=y(-j),\quad h(j)^*=h(-j),\quad y(j)^*=x(-j).
$$
For $\pi = b_1 b_2 \dots b_s \in \Cal P (\bar B)$ we define a
{\it dual} colored 
partition $\pi^*= b_s^* b_{s-1}^* \dots b_1^*$.

For $N \in \Bbb Z$ and a partition $\pi \: \bar B \to \Bbb
Z$ we may define a 
{\it translated} partition  $\pi_N = \pi \circ\tau_N$, where
$$
   \tau_N (x(i)) = x(i-N),\quad \tau_N (h(i)) = h(i-N),\quad
\tau_N (y(i)) = y(i-N).
$$ 
It follows from Proposition 6.6.1 that 

\proclaim{Proposition 7.1.2} The set $\lt(\bar R)$ is
translation invariant and
self-dual (i.e. for $\rho \in \lt(\bar R)$ we have $\rho^*
\in \lt (\bar R)$).
\endproclaim
We shall sometimes also say that $\pi$ and $\pi'$ are dual
to each other if
$\pi'$ is dual to some translation of $\pi$. In this way we
shall also speak of
dual embeddings.

\subheading{7.2. Colored partitions of length $k+2$ allowing
at least two
embeddings}

 From Proposition 6.6.1 we get the following list of
partitions of length $k+2$
 allowing two or more than two embeddings (here $N$ denotes
the number of
embbedings):

$(1)\quad \sh \pi = j^{k+2} $,
$$\alignat 2
&h(j)^a x(j)^{k+2-a}, \quad 1 \le a \le k+1,
\qquad\qquad\qquad & N & =2,\\
&y(j) h(j)^k x(j), \quad & N & =2,\\
&y(j)^{k+2-a} h(j)^a, \quad 1 \le a \le k+1, \qquad & N &
=2,\\
\endalignat
$$

$(2)\quad \sh\pi=(j-1)^a j^b, \quad a,b > 0,  \quad 
a+b=k+2$,
$$\alignat 2
&x(j-1)^a x(j)^b, & \qquad  N & =2,\\
&x(j-1)^a h(j)^c x(j)^{b-c}, \quad 1 \le c \le b - 1, & N &
=3,\\
&x(j-1)^a h(j)^b, \quad & N & =2,\\
(*)\quad&x(j-1)^a y(j) h(j)^{b-2} x(j), \quad & N & =2,  \\
&x(j-1)^a y(j)^c h(j)^{b-c}, \quad 1 \le c \le b - 1, & N &
= 3,\\
&x(j-1)^a y(j)^b, & N& = 2, \\
&h(j-1) x(j-1)^{a-1} y(j)^{b-1} h(j),& N & = 2,\\
&h(j-1)^c x(j-1)^{a-c} y(j)^b, \quad 1 \le c \le a - 1, & N
& = 3,\\
&h(j-1)^a y(j)^b, \quad & N & =2,\\
(*)\quad&y(j-1) h(j-1)^{a-2} x(j-1) y(j)^b, \quad & N & =2, 
\\
&y(j-1)^c h(j-1)^{a-c} y(j)^b, \quad 1 \le c \le a - 1 ,& N
& =3,\\
&y(j-1)^a y(j)^b, \quad & N & = 2.
\endalignat
$$

Note that for $a \ge 2$ and $b \ge 2$ elements denoted by
(*) are
well defined. Also note that elements are listed so that
their weights
descend from $(k+2)\alpha$ to $-(k+2)\alpha$, where for each
element denoted
by (*) there is another element of the same weight.
In the case $b=1$ (and $a=k+1$) the list should be modified
by elements 
$$
\cases  x(j-1)^a h(j), \qquad & N = 2, \\
h(j-1) x(j-1)^{a-1} x(j), \qquad & N = 2, \endcases
$$ 
and in the case $a=1$ (and $b=k+1$) by
$$
\cases  h(j-1) y(j)^b, \qquad & N = 2, \\
y(j-1) y(j)^{b-1} h(j), \qquad\qquad & N = 2. \endcases
$$

$(3)\quad \sh \pi = (j -1) j^k (j + 1) $,
$$\alignat 2
&x(j-1) x(j)^k x(j+1)^a h(j+1)^b, \quad  a, b \ge 0,
a+b = 1,\quad &N&=2,\\
&x(j-1) h(j)^a x(j)^{k-a} y(j+1), \quad  0 \le a \le
k,\quad&N&=2,\\
&x(j-1) y(j)^a h(j)^{k-a} y(j+1), \quad  1 \le a \le
k,&N&=2,\\
&y(j-1)^a h(j-1)^b y(j)^k y(j+1), \quad  a, b \ge 0,
a+b = 1,&N&=2.
\endalignat
$$

Note that the above list of embeddings is self dual.

\subheading{7.3. Linked embeddings and exceptional cases}

Let us call two embeddings  $\rho, \rho' \subset \pi$ {\it
linked} 
if there is a sequence 
$$
\rho_1 =\rho \subset \pi, \  \rho_2 \subset \pi, \dots,
\  \rho_{n-1} \subset \pi,\  \rho_{n}= \rho' \subset \pi
$$
of embeddings such that $\ell(\rho_{i-1} \cup \rho_i) <
\ell(\rho \cup \rho') $ 
for each $i=2,\dots,n$.

The above definition is designed for inductive arguments.
However, there are
{\it exceptional cases} that we will have to consider
separately:

\proclaim{Lemma 7.3.1} Let $\pi$ be a colored partition, 
$\rho_1, \rho_2 \subset \pi$ leading terms of relations such
that
$\pi=\rho_1 \cup \rho_2$, $\rho_1 \cap \rho_2 \neq
\varnothing$,
 $\ell(\pi)\ge k+3$, and the two embeddings $\rho_1 \subset
\pi$,
$\rho_2 \subset \pi$ are not linked. Then $\pi$ is one of
the following colored
partitions:
$$\allowdisplaybreaks \align
&y(j)^a h(j)^{k+1-a} x(j)^a,\tag{1}\\ 
  & \qquad\qquad\qquad 2 \le a \le k, \\ 
&x(j-1)^{k+1-a-b} y(j)^a h(j)^b x(j)^a, \tag{2}\\
    &\qquad\qquad\qquad a \ge 2, a+b \le k,\\
&h(j-1)^{k+1-a-b} x(j-1)^a y(j)^b x(j)^{k+1-a},\tag{3}\\
    &\qquad\qquad\qquad 1 \le a \le k - 1, \\
&h(j-1)^{k+1-a-b} x(j-1)^a y(j)^b h(j)^{k+1-a-b}, \tag{4}\\ 
    &\qquad\qquad\qquad 1 \le a+b \le k-1,\\
&y(j-1)^{k+1-a} x(j-1)^b y(j)^a h(j)^{k+1-a-b}, \tag{5}\\ 
    &\qquad\qquad\qquad 1 \le a \le k-1,\\
&y(j-1)^a h(j-1)^b x(j-1)^a y(j)^{k+1-a-b},\tag{6}\\ 
    &\qquad\qquad\qquad a \ge 2, a+b\le k, \\
&y(j-1)^{k+1-a-c} h(j-1)^c x(j-1)^b y(j)^a h(j)^{k+1-a-b},
\tag{7}\\ 
    &\qquad\qquad\qquad  a, c \ge 1, a+b+c \le k,\\
&h(j-1)^{k+1-a-b} x(j-1)^a y(j)^b h(j)^c
x(j)^{k+1-a-c},\tag{8}\\
    &\qquad\qquad\qquad  a, c \ge 1, a+b+c \le k,\\
& x(j-2)^a y(j-1)^b h(j-1)^{k+1-a-b} y(j)^a,\tag{9}\\ 
    &\qquad\qquad\qquad 2 \le a \le k,\\ 
& x(j-2)^a h(j-1)^{k+1-a-b} x(j-1)^b y(j)^a, \tag{10}\\ 
    &\qquad\qquad\qquad 2 \le a \le k,\\
&x(j-2)^{k+1-a-b} h(j-1)^a x(j-1)^b y(j)^c
h(j)^{k+1-b-c},\tag{11}\\
    &\qquad\qquad\qquad  1 \le b \le k-1, a+b+c \le k,\\
&x(j-2)^{k+1-a-b} h(j-1)^a x(j-1)^b h(j)^c
x(j)^{k+1-b-c},\tag{12}\\
    &\qquad\qquad\qquad  1 \le b \le k-1, a+b \le k, \\
& h(j-2)^{k+1-b-c} x(j-2)^c y(j-1)^b h(j-1)^a
y(j)^{k+1-a-b}, \tag{13}\\
    &\qquad\qquad\qquad  1 \le b \le k-1, a+b+c \le k,\\
& y(j-2)^{k+1-b-c} h(j-2)^c y(j-1)^b h(j-1)^a
y(j)^{k+1-a-b}, \tag{14}\\
    &\qquad\qquad\qquad  1 \le b \le k-1, a+b \le k.\\             
\endalign
$$
\endproclaim
\noindent Here (and later on) it goes without saying that
$X^a$, $X \in \bar B$, is
defined only for integers $a \ge 0$ (in particular $a=0$ is
allowed unless
explicitly ruled out).

Note that the following cases are dual to each other:
$$\gather
(1)^*=(1), \ (2)^*=(6), \ (3)^*=(5), \ (4)^*=(4),\\
(7)^*=(8), \ (9)^*=(10), \ (11)^*=(13), \ (12)^*=(14).  
\endgather$$
We could change case (7) to include c=0 covering (5), and 
we could change case (8) to include c=0 covering (3), but
for later proofs
it will be convenient to keep them separated.
Note also that in each of the cases (1)--(14) there are only
two embeddings
$\rho_1 \subset \pi, \rho_2 \subset \pi$ of leading terms of
relation, so that 
$\rho_1, \rho_2$ are determined up to permutation.

\demo{Proof} 
We shall often use the following ``replacement of
representative":
$$
(A) \left\{ \aligned
& \text{Let} \ \pi=\pi' y(j)^d, \ \rho_i=\rho' y(j)^c \ (i=1
\ \text{or} \ 2),
            \ \text{where} \  \sh(\rho_i)=(j-1)^{k+1-c}
j^c\\
& \text{and} \ \ y(j)^d, y(j)^c \ \ \text{are the maximal
powers of}\ y(j) \
           \text{contained in} \  \pi \ \text{and} \ 
\rho_i\\
&\text{(respectively). Then we may, by replacement of} \
\rho_i \ \text{by} \ 
          \rho_3=\rho'' y(j)^d \ \text{with}  \\
& \rho'' \subset \rho' \ \text{if necessary, assume that} 
           \  c=d, \  \text{that is} \ \rho_i \supset
y(j)^d.\\
&  \qquad \text{Dually, if} \ \pi \supset x(j-1)^d,
\rho_i \supset x(j-1)^c, \sh(\rho_i)=(j-1)^c j^{k+1-c},\\
& \text{then we may assume that} \quad \rho_i \supset
x(j-1)^d .\endaligned \right.
$$

 Assume first that $\sh(\pi)=j^\ell$,
 $k+3 \le \ell \le 2k+1$. Then $\pi=y(j)^a h(j)^b x(j)^c$.
Clearly
$\rho_1, \rho_2$ are linked if $a$ or $c=0$, so $a, c \ge
1$. We may then
assume $\rho_1 \subset y(j)^a h(j)^b$, $\rho_2 \subset
h(j)^b x(j)^c$. If
$a+b$ of $b+c > k+1$, $\rho_1$ and $\rho_2$ are linked. Thus
$\pi$ must be
of type (1).

Consider now the case that $\sh(\pi)=(j-1)^\ell j^m$ for
some $\ell, m \ge 1$.
Then
$$
\pi=y(j-1)^a h(j-1)^b x(j-1)^c y(j)^d h(j)^e x(j)^f
$$
for some $a, \dots , f$. It cannot be that
$\sh(\rho_1)=(j-1)^{k+1}$ and
 $\sh(\rho_2)=j^{k+1}$ since $\rho_1 \cap \rho_2 \ne
\varnothing$. We may thus assume
that  $\sh(\rho_1)=(j-1)^p j^{k+1-p}$ for some $1 \le p \le
k$.
By duality and (A) we may assume that $\rho_1 \supset
y(j)^d$. By (A) we may
also assume that $\rho_2$ satisfies one of the following:
$$\allowdisplaybreaks\align
&\rho_2 \supset y(j)^d, \quad \sh (\rho_2)=(j-1)^{k+1-d}
j^d, \tag"{$(\alpha)$}"\\
&\rho_2 \supset x(j-1)^c, \ c \ge 1, \quad 
              \sh(\rho_2)=(j-1)^c j^{k+1-c} ,
\tag"{$(\beta)$}"\\
&\rho_2 \subset y(j)^d h(j)^e, \tag"{$(\gamma)$}"\\
&\rho_2 \subset y(j-1)^a h(j-1)^b x(j-1)^c.
\tag"{$(\delta)$}"
\endalign$$

Case ($\alpha$): If $\rho_1=y(j-1)^a h(j-1)^{b'} y(j)^d$, 
 $\rho_2=h(j-1)^{b''} x(j-1)^c y(j)^d$, then $\rho_1$ and 
$\rho_2$ are 
linked unless $b'=b=b''$ in which case $\pi$ is of type (6).
If $a=0$ or
$c=0$, $\rho_1$ and  $\rho_2$ are linked.

Case ($\beta$): If $\rho_1 \subset y(j-1)^a h(j-1)^b y(j)^d$
and
 $\rho_2 \subset x(j-1)^c h(j)^e x(j)^f$ then $\rho_1 \cap
\rho_2 =
 \varnothing$.
By duality we may thus assume that  $\rho_1=h(j-1)^b
x(j-1)^{c'} y(j)^d$
and $a=0$. Then $\rho_1$ and $\rho_2$ are linked if $c' < c$
so $c'=c$.
If  $\rho_2=x(j-1)^c h(j)^e x(j)^f$ we have $\pi$ of type
(8) or (3).
In case $\pi$ of type (8) we must have $a+b+c \le k$,
otherwise
$\rho_1$ and $\rho_2$ are linked via $\rho_3 = x(j-1)^a
y(j)^b h(j)^{k+1-a-b}$.
If  $\rho_2=x(j-1)^c y(j)^{d'} h(j)^e$ then  $\rho_1$ and
$\rho_2$ are linked
 if $d' < d$, so $d'=d$ and $\pi$ is of type (4).

Case ($\gamma$): If $c \ge 1$, then  $\rho_1$ and $\rho_2$
are linked
via $\rho_3=x(j-1) y(j)^d h(j)^{k-d}$, so $c=0$. If
 $\rho_2 = y(j)^{d'} h(j)^e$ with $d' < d$, then  $\rho_1$
and $\rho_2$
 are linked via  $\rho_3 = y(j)^d h(j)^{k+1-d}$. Thus $d'=d$
and $\pi$
is of type (7).

Case ($\delta$): If $\rho_2 = y(j-1)^{k+1-b'}h(j-1)^{b'}$
and 
$\rho_1 = h(j-1)^{b''} x(j-1)^c y(j)^d$ then it is easy to
see that $\rho_1$
and $\rho_2$ are linked unless $b' = b = b''$. In that case
$\rho_1$ and
$\rho_2$ are linked  via $y(j-1)^{k-b} h(j-1)^b y(j)$. If
$c=0$ and
$\rho_1 = y(j-1)^{a'} h(j-1)^{b'} y(j)^d$, then  $\rho_1$
and
$\rho_2$ are linked  via either
 $\rho_3 = y(j-1)^{a' +1} h(j-1)^{b'} y(j)^{d-1}$ or
 $\rho_3 = y(j-1)^{a'} h(j-1)^{b'+1} y(j)^{d-1}$. Similarly
$\rho_1$ and
$\rho_2$ are linked in all cases where $\rho_2 \subset
h(j-1)^b x(j-1)^c$.

There remains the case that $\sh(\pi) = (j-2)^p (j-1)^q j^r$
for some
$p, q, r \ge 1$. We may then assume  $\sh(\rho_1) = (j-2)^p
(j-1)^{k+1-p}$,
 $\sh(\rho_2) = (j-1)^{k+1-r} j^r$. Since $\rho_1 \cap
\rho_2 \ne \varnothing$
it cannot be that $\rho_1 \supset y(j-1)^{k+1-p}$,
 $\rho_2 \supset x(j-1)^{k+1-r}$. By duality we may thus
assume that
 $\rho_1 \supset x(j-2)^p$. Then either
$$\align
 &\rho_2 \supset x(j-1)^{k+1-r}, \
\text{or}\tag"{$(\varepsilon$)}"\\ 
&\rho_2 \supset y(j)^r. \tag"{$(\zeta$)}" 
\endalign$$

Case ($\varepsilon$): By (A) $x(j-1)^{k+1-r}$ is the maximal
power of 
$x(j-1)$ contained in $\pi$. We must have
 $\rho_1 =x (j-2)^p h(j-1)^{k+1-p-a} x(j-1)^a$ for some $a$.
If $a < k+1-r$, then $\rho_1$
and $\rho_2$ are linked via $x(j-2)^{p-1} h(j-1)^{k+1-p-a}
x(j-1)^{a+1}$.
Thus $a=k+1-r$ and $\pi$ is of type (10), (11) or (12). If
$\pi$ is of type (11)
and $a \ge 1$, then $a+b+c \le k$ since otherwise $\rho_1$
 and $\rho_2$ are linked via $\rho_3 =h(j-1)^{a'} x(j-1)^b
y(j)^{k+1-a'-b}$
for suitable $a' \ge 1$.

Case($\zeta$): If $\rho_1 =x(j-2)^p y(j-1)^a
h(j-1)^{k+1-a-p}$
and $\rho_2 = h(j-1)^{k+1-b-r} x(j-1)^b y(j)^r$ and
 $1 \le a \le b$, then $\rho_1$ and $\rho_2$ are linked via 
$\rho_3 = x(j-2)^p h(j-1)^{k+1-a-p} x(j-1)^a$.
If instead $1 \le b \le a$, then  $\rho_1$ and $\rho_2$ are
linked
 via $\rho_3 = y(j-2)^b h(j-1)^{k+1-b-r} y(j)^r$. Similarly 
$\rho_1$ and $\rho_2$
 are linked if $\rho_1 = x(j-2)^p h(j-1)^{k+1-a-p}
x(j-1)^a$,
$\rho_2 = y(j-1)^b h(j-1)^{k+1-r-b} y(j)^r$. We must
therefore have
$\pi = x(j-2)^p y(j-1)^a h(j-1)^b y(j)^r$ or
$\pi = x(j-2)^p h(j-1)^a x(j-1)^b y(j)^r$, and now it is
easy to see that
 $\rho_1$ and $\rho_2$ are linked except in case $\pi$ is of
type (9) or (10).
\qed
\enddemo

\head{8. Relations among relations}\endhead

\bigskip By Theorem 5.14 we have $M^1(\Lambda) = U_k(\tilde
\goth g) \bar R v_\Lambda$, 
so $M^1(\Lambda)$ is spanned by elements of the form $u(\rho
\subset \pi)v_\Lambda
= u(\pi /\rho)r(\rho)v_\Lambda$, where $\rho$'s are leading
terms of ``relations" 
$r(\rho)$ in $\bar R$  and $\pi$'s 
are in the ideal $(\lt(\bar R ))$. We will prove that for
two elements of the form 
$u(\rho \subset \pi)$ and $u(\rho' \subset \pi)$ there are
relations that hold on 
every  highest weight module of level $k$, 
and we shall refer to them as {\it relations among
relations}. 

\bigskip

\proclaim {Proposition 8.1}
We have
$$
\sum_{j \in \Bbb Z}\, \bigl((k+2)j-n \bigr)\, x(j)
r_{(k+1)\alpha}(n-j) = 0
\tag{8.1}
$$
in $\overline {U_k(\tilde \goth g)}$ for all $n \in \Bbb Z$.
\endproclaim

Note that $x(j)$ and $r_{(k+1)\alpha} (n-j)$ commute so that the
infinite sum is convergent. 

\demo{Proof} It follows from $[L_{-1}, x(-1)] = x(-2)$ and
$L_{-1} \bold 1 = 0$ that
$$
   L_{-1} x(-1)^{k+2} \bold 1 = (k+2) x(-2) x(-1)^{k+1}
\bold 1.
\tag {8.2}$$
By (2.4), (2.13) and (2.14)
$$\gather
   Y(L_{-1} x(-1)^{k+2} \bold 1, z) = \frac d{dz} \left (
x(z) Y(r_{(k+1)\alpha}, z) \right),\\
Y(x(-2) x(-1)^{k+1} \bold 1, z) = \Bigl (\frac d{dz} x(z)
\Bigr) Y(r_{(k+1)\alpha}, z).
\endgather$$
Thus (8.2) implies
$$
\frac d{dz} \Bigl ( x(z) Y(r_{(k+1)\alpha}, z) \Bigr)
= (k+2) \Bigl (\frac d{dz} x(z) \Bigr) Y(r_{(k+1)\alpha}, z).
\tag {8.3}$$
The coefficient of $z^{-n-k-3}$ in (8.3) gives (8.1).
\qed
\enddemo
\proclaim{Proposition 8.2}
We have
$$\gather
\sum_{j<0}
\left((k+2)j-n\right)\{x(j)r_{(i-1)\alpha}(n-j)-h(j)r_{i\alpha}(n-j)-y(j)
r_{(i+1)\alpha}(n-j)\} \tag{8.4}\\
+\sum_{j\ge 0}
\left((k+2)j-n\right)\{r_{(i-1)\alpha}(n-j)x(j)-r_{i\alpha}(n-j)h(j)
-r_{(i+1)\alpha}(n-j)y(j)\} =0
\endgather 
$$
in $\overline {U_k(\tilde \goth g)}$ for all $n \in \Bbb Z$,
$i = -k-2, \dots , k+2$.
\endproclaim

\demo{Proof}
The case $i=k+2$ is (8.1). If we denote the left hand side
in (8.4) by $s_i$ then one finds
$$
  [y(0), s_i] = (k+3-i) s_{i-1}
$$
and (8.4) follows by induction. 
\qed \enddemo

\proclaim{Proposition 8.3}
We have
$$\gathered
(k+2-i) \biggl( \sum_{j<0} x(j) r_{(i-1)\alpha}(n-j) +\sum_{j \ge 0}
r_{(i-1)\alpha}(n-j) x(j)\biggr) \\
+ i \biggl( \sum_{j<0} h(j) r_{i\alpha}(n-j) +\sum_{j \ge 0} 
r_{i\alpha}(n-j)h(j) \biggr) \\
+(k+2+i)\biggl( \sum_{j<0} y(j) r_{(i+1)\alpha}(n-j) +\sum_{j \ge 0}
r_{(i+1)\alpha}(n-j) y(j)\biggr) \\
= (k+2) (-n-k-1) r_{i\alpha}(n)
\endgathered \tag{8.5}$$
in $\overline {U_k(\tilde \goth g)}$ for all $n \in \Bbb Z$,
$i = -k-1, \dots , k+1$.
\endproclaim

\demo{Proof}
By (3.17), (3.18) and (6.6.1)
$$\gather
  L_{-1} r_{i\alpha} = \tfrac 1{2(k+2)} \bigl( 2 x(-1) y(0) + h(-1)
h(0) + 2 y(-1) x(0) \bigr) r_{i\alpha} \\
= \tfrac 1{k+2} \bigl( (k+2-i) x(-1) r_{(i-1)\alpha} + i h(-1) r_{i\alpha} +
(k+2+i) y(-1) r_{(i+1)\alpha} \bigr).
\endgather $$ From
(2.4) and (2.14) it follows that
$$\gathered
(k+2-i)\bigl( x(z)^- Y(r_{(i-1)\alpha},z) + Y(r_{(i-1)\alpha},z)x(z)^+\bigr)
\\
+ i \bigl( h(z)^- Y(r_{i\alpha},z) + Y(r_{i\alpha},z)h(z)^+\bigr) \\
+ (k+2+i) \bigl( y(z)^- Y(r_{(i+1)\alpha},z) +
Y(r_{(i+1)\alpha},z)y(z)^+\bigr) \\
= (k+2) \frac d{dz} \,Y(r_{i\alpha},z).
\endgathered $$
The coefficient of $z^{-n-k-2}$ gives (8.5). \qed
\enddemo

Let $j \in \Bbb Z$, $j < 0$, $0 \le p \le k+2$ and consider
the shape
$$
  \sigma = (j-1)^p \,j^{k+2-p}.
$$
 Set  $n = \vert \sigma \vert = (k+2)j-p$. By taking linear
combinations of (8.4)
and (8.5) we can produce relations in $\overline {U_k(\tilde
\goth g)}$ of the form
$$\gather
C^x x(j\!-\!1)r_{(i-1)\alpha}(n\!-\!j\!+\!1) + C^h h(j\!-\!1)
r_{i\alpha}(n\!-\!j\!+\!1) 
+ C^y y(j\!-\!1) r_{(i+1)\alpha}(n\!-\!j\!+\!1) \\
+ C_x\,x(j)r_{(i-1)\alpha}(n-j) 
+ C_h \,h(j) r_{i\alpha}(n-j) + C_y \, y(j)
r_{(i+1)\alpha}(n-j) \tag{8.6}\\
+ \text{ terms of higher shape } = 0.
\endgather  $$
Calculations of some $2 \times 2$ determinants show the
following:
\proclaim{Lemma 8.4}
There is a relation (8.6) with
$$\align
 \{C^x, C_x\} &= \{0, 1\}  \quad \text{unless} \quad i =
k+2,\\
 \{C^y, C_y\} &= \{0, 1\}  \quad \text{unless} \quad i =
-k-2,\\
 \{C^x, C^h\} &= \{0, 1\}  \quad \text{unless} \quad \sigma
= (j-1)^{k+2},\\
 \{C^h, C^y\} &= \{0, 1\}  \quad \text{unless} \quad \sigma
= (j-1)^{k+2},\\
 \{C_x, C_h\} &= \{0, 1\}  \quad \text{unless} \quad \sigma
= j^{k+2},\\
 \{C_h, C_y\} &= \{0, 1\}  \quad \text{unless} \quad \sigma
= j^{k+2}.
\endalign$$
\endproclaim

\remark{Remark} The identities (8.4) and (8.5) are analogous to
\cite{MP, Proposition 9.3} in the principal picture.
\endremark

\remark{Remarks} Since we have exact sequences
$$
0 \rightarrow M^1(\Lambda) \rightarrow M(\Lambda)\rightarrow
L(\Lambda) \rightarrow 0,
$$
$$
0 \rightarrow N^1(k \Lambda_0) \rightarrow N(k \Lambda_0)
\rightarrow
L(k \Lambda_0) \rightarrow 0,
$$
we may call elements of $M^1(\Lambda)=U(\tilde{\goth g})
\bar R v_\Lambda$
and/or of $N^1(k\Lambda_0)=U(\tilde{\goth g}) \bar R \bold 1 $
the relations (which
define the standard module $L(\Lambda)$ and/or $L(k
\Lambda_0)$).

Obviously $\bar R v_\Lambda \subset M(\Lambda)$ is a $(\goth h^e
+\tilde{\goth n}_+)$-module,
and $\bar R \bold 1 \subset N(\Lambda)$ is a $\tilde \goth g_{\ge
0}$-module, so we have
exact sequences
$$
0 \rightarrow \ker \Psi_{M(\Lambda)} \rightarrow 
U(\tilde \goth g) \otimes_{U(\goth h^e +\tilde{\goth n}_+)}
\bar{R} v_\Lambda
 \overset{\Psi_{M(\Lambda)}} \to{\longrightarrow} M^1
(\Lambda) \rightarrow 0,
$$
$$  
0 \rightarrow \ker \Psi_{N(k \Lambda_0)} \rightarrow 
U(\tilde \goth g) \otimes_{U(\tilde{\goth g}_{\ge 0})}
\bar{R} \bold 1
 \overset{\Psi_{N(k \Lambda_0)}} \to{\longrightarrow} N^1 (k
\Lambda_0) \rightarrow 0,
$$
where $\Psi_{M(\Lambda)} (u \otimes r v_\Lambda)= u r
v_\Lambda $,
$\Psi_{N(k \Lambda_0)} (u \otimes r \bold 1) = u r \bold 1$.
 Since $\ker \Psi$
defines the quotient (``which are relations"), we may call the
elements of
 $\ker   \Psi$ relations among relations.

One can see that the $U_k (\tilde{\goth g})$-module
$\ker \Psi_{N(k \Lambda_0)}$ is generated by the elements $q^1
(n)$, 
$n \le - k - 3$, and $q^2 (n)$, $n \le - k - 2$, where
$$\align
q^1 (n) = &\sum_{j \leq -1} ((k+2) j - n) x(j) \otimes
r_{(k+1)\alpha}
(n - j) \bold 1,\\
q^2 (n) = &- \sum_{i \leq -1} (k+1) h(i) \otimes r_{(k+1)\alpha} (n
- i) \bold 1\\
&  -\sum_{i \leq -1} x (i) \otimes r_{k\alpha} (n - i) \bold 1\\
& - (k+2) (n+k+1) 1 \otimes r_{(k+1)\alpha} (n) \bold 1.
\endalign
$$
So in some sense $q^1(n)$ and $q^2 (n)$ (or the
corresponding (8.4) and (8.5)) 
form a ``complete generating set" of relations among
relations. 

For an element of the form $u \otimes r\bold 1$ in
$U(\tilde{\goth g}_{< 0}) \otimes
\bar R \bold 1$ define a ``leading term" $\ell t(u \otimes r
\bold 1) 
= (\ell t (u), \ell t (r \bold 1))$. We
may identify a ``leading term" $(\kappa,\rho)$ with the
embedding $\rho \subset
\kappa \rho$. So each nonzero element $q$ in $U
(\tilde{\goth g})
 \otimes_{U(\tilde{\goth g}_{\ge 0})} \bar R \bold 1 \cong
U(\tilde{\goth g}_{< 0}) \otimes \bar R \bold 1$ can be
written in the form
$$
q =  \sum_{\rho' \subset \pi} C_{\rho',\pi}\, u(\pi/\rho')
\otimes
r (\rho') \bold 1
 + \sum \Sb \rho' \subset \pi'\\ \pi' \succ \pi \endSb
C_{\rho', \pi'}\,
 u(\pi'/\rho') \otimes r (\rho') \bold 1 
$$
for some $\pi \in \Cal P (\bar B_{< 0})$ and
$C_{\rho',\pi}$,
 $C_{\rho',\pi'} \in \Bbb F$, where at least one coefficient
$C_{\rho',\pi}$
is nonzero. Then $q \in \ker \Psi_{N(k \Lambda_0)}$, i.e.
$$
 \sum_{\rho' \subset \pi} C_{\rho',\pi} u(\rho' \subset \pi)
\bold 1 +
 \sum \Sb \rho' \subset \pi'\\ \pi' \succ \pi \endSb
C_{\rho', \pi'}
 u(\rho' \subset \pi') \bold 1 = 0,
$$
implies $\sum_{\rho' \subset \pi} C_{\rho',\pi} = 0$. Hence
for such $q$ there
must be at least two embeddings $\rho' \subset \pi$, $\rho''
\subset \pi$ with
the corresponding coefficients being nonzero.
\endremark

\heading 9. Relations among relations for two embeddings
\endheading

We fix a map (analogous to the one in Section 6.5)
$$
\lt (\bar{R}) \rightarrow \bar R \quad , \quad \rho \mapsto
r(\rho),
$$
such that $\lt (r(\rho)) = \rho$. For $\rho \in \lt(\bar R)$
and $\pi \in (\lt (\bar{R}))$, $\rho \subset \pi$, we define
$$
u(\rho \subset \pi) = u(\pi / \rho)r(\rho).
$$
For $\pi \in \Cal P(\bar B)$ set
$$
\Cal R_{(\pi)} = \overline {\lspan\{u(\rho \subset \pi')\mid
 \rho \in \lt (\bar{R}), \pi'\in(\lt (\bar{R})),
 \rho \subset \pi', \pi \prec \pi'\}},
$$
the closure taken in $\overline{U_k(\tilde{\goth g})}$.

The main result in this section is the following consequence
of relations among
relations constructed in Section 8:
\proclaim{Theorem 9.1} Let $\pi \in (\lt(\bar{R}))$,
$\rho_1, \rho_2 \in
\lt(\bar{R})$, $\rho_1 \subset \pi$, $\rho_2 \subset \pi$.
Then
$$
u(\rho_1 \subset \pi) \in \Bbb F^\times u(\rho_2 \subset
\pi) + \Cal R_{(\pi)}.
\tag{9.1}
$$
\endproclaim
\remark{Remark}
Since $N^1(k\Lambda_0)=U(\tilde{\goth g}_{<0})\bar{R}\1 $, 
the set of vectors
$u(\rho \subset \pi)\1$, where $\rho \in \lt(\bar{R} \1)$,
 $\pi \in (\lt(\bar{R} \1))$, $\rho \subset \pi$,
is a spanning set of $N^1(k\Lambda_0)$. 
Let $(\lt(\bar{R} \1))\ni\pi \mapsto \rho(\pi)\in\lt(\bar{R}
\1)$ be such that
$\rho(\pi)\subset \pi$. Using Theorem 9.1 and induction we
see that
$$
u(\rho(\pi) \subset \pi)\1, \qquad \pi \in (\lt(\bar{R}
\1)), 
$$
is a spanning set as well. Clearly this set is linearly
independent,
i.e. it is a basis of $N^1(k\Lambda_0)$. As a consequence we
have
$$
\lt_{N(k\Lambda_0)}(N^1(k\Lambda_0)) = (\lt(\bar R \1)),
$$
which in turn implies (see the proof of Theorem 6.5.5) that
the set of vectors
$$
u(\pi)v_{k\Lambda_0}, \qquad \pi \in 
\Cal P (\bar B_{<0}) \backslash (\lt(\bar R \1)),
$$
is a basis of the standard $\tilde{\goth g}$-module
$L(k\Lambda_0) 
= N(k\Lambda_0)/N^1(k\Lambda_0)$. 
\endremark
 The rest of this section is devoted to the proof of Theorem
9.1. 

We may and we 
will assume that the map $\rho \mapsto r(\rho)$ is such that
$r(\rho)$ is 
homogeneous of weight $\wt (\rho)$ and degree $\vert \rho
\vert$.
For example, we may take the map $\pi_{m\alpha}(n) \mapsto r_{m\alpha}(n)$,
$n \in \Bbb Z$,
$|m| \le k+1$ (cf. Section 6.6).

\proclaim{Lemma 9.2} Let $ \rho_1, \rho_2 \in \lt(\bar R)$,
$\pi = \rho_1 \cup \rho_2$, $\ell(\pi) = k+2$. Then
$$
u(\rho_1 \subset \pi) \in \Bbb F^\times u(\rho_2 \subset
\pi) + \Cal R_{(\pi)}.
$$
\endproclaim
\demo{Proof} In Section 7.2 we have listed all $\pi$ of
length $k+2$ which
allow more than one embedding. For a fixed degree and weight
only one of three
 distinct cases may occur: among all colored partitons of
lenght $k+2$ of
fixed degree and weight a) there is one $\pi$ with two
embeddings, b) there
are two $\pi_1$ and $\pi_2$ with two embeddings each, and c)
there is one $\pi$
with three embeddings. Note that we have a) for $\wt (\pi) =
\pm (k+2)\alpha$,
and b) and c) for $\wt (\pi) = j\alpha$, $-k-1 \le j \le
k+1$. Also note that cases
(1) and (3) in Section 7.2 \  go together for $\vert \pi
\vert \equiv 0 \mod{(k+2)}$.

In Section 8 we have constructed relations among relations
(which hold on every
highest weight module of level $k$) of the form
$$
\sum c_i u(\rho_i \subset \pi_i) = 0, \tag{9.2}
$$
where (for the sake of notation we assume that) $c_i \ne 0$
and $\pi_i = \pi_j$
 implies $\rho_i \ne \rho_j$.
Note that relation (8.4) holds for weights $j\alpha$, $-k-2
\le j \le k+2$, and
that (8.5) holds for weights $j\alpha$, $-k-1 \le j \le
k+1$.
So in the case a) we have  one linearly independent
relation, and in the case b)
or c) there are two linearly independent relations.

Fix a degree and weight and write a relation in the form
(9.2). Let $\pi$ be
the smallest among $\pi_i$'s which appear in the sum, 
say $\pi = \pi_1$. Then 
$$
c_1 u(\rho_1\subset\pi)+
 \sum\Sb{i\ne 1}\\{\pi\preccurlyeq \pi_i}\endSb c_i
u(\rho_i\subset \pi_i)=0.
$$
In the case that $\pi \prec \pi_i$ for all other $\pi_i$ we
would have that the 
expression on the left has a leading term $\pi$, and then
the expression on the 
left would be different from 0. Hence there are 
$\rho_2 \subset \pi_2, \dots,\rho_s\subset \pi_s$ such that
$\pi_2 = \dots =\pi_s = \pi$ and we have

$$
c_1 u(\rho_1\subset\pi)+ \dots +c_s u(\rho_s\subset\pi)+
 \sum_{\pi \prec \pi_i} c_i u(\rho_i\subset \pi_i)=0.
\tag{9.3}$$
So in particular there exists $\pi$ with at least two
embeddings 
$\rho_1 \subset \pi, \rho_2 \subset \pi, \dots$.

Assume that we have the case a). Then $\rho_1 \subset \pi$
and
 $\rho_2 \subset \pi$ are uniquely determined and (9.3) is
the statement of our
lemma.

Assume that we have the case b). Then we have two linearly
independent
relations: $R$ of the form (9.2) and $R'$ of the form
$$
\sum c' _i u(\rho' _1 \subset \pi' _i) = 0, \tag {9.4}
$$
where $c' _i \ne 0$, and $\pi' _i = \pi' _j$ implies $\rho'
_i \ne \rho' _j$.
We also have two colored partitions, say $\pi_1$ and $\pi_2$
with two
embeddings each.

Each relation (9.2) and (9.4) can be written in the form
(9.3) (with $s=2$), 
in each case $\pi', \pi \in\{\pi_1,\pi_2\}$. If $\pi' \ne
\pi$, then we have
the lemma. If for both relations (say) $\pi = \pi' = \pi_1$,
then the terms
$c_1 u(\rho_1 \subset \pi_1) + c_2 u(\rho_2 \subset \pi_1)$
and
$c'_1 u(\rho'_1 \subset \pi_1) + c'_2 u(\rho'_2 \subset
\pi_1)$ are
proportional by a factor $\lambda \ne 0$. Then the
(nontrivial) relation
$R - \lambda R'$ written in the form (9.3) should give the
only possible
$\pi = \pi_2$, and the lemma holds.

Assume that we have the case c). Then two linearly
independent relations (9.2)
and (9.4) can be written in the form (9.3)
$$\gathered
c_1 u(\rho_1 \subset \pi) + c_2 u(\rho_2 \subset \pi) + 
c_3 u(\rho_3 \subset \pi) +
 \sum\Sb{\pi \prec \pi_i}\endSb c_i u(\rho_i \subset \pi_i)
= 0,\\
c'_1 u(\rho_1 \subset \pi) + c' _2 u(\rho_2 \subset \pi) + 
c'_3 u(\rho_3 \subset \pi) +
 \sum\Sb{\pi \prec \pi'_i}\endSb c'_i u(\rho'_i \subset
{\pi'}_i) = 0.
\endgathered$$
The matrix $\left(\matrix c_1 & c_2 & c_3 \\
c'_1 & c'_2 & c'_3\endmatrix\right)$ must have rank 2 (since
otherwise the
nontrivial relation of the form $R'' = R - \lambda R'$ would
give another
$\pi'' \ne \pi$ allowing more than one embedding). Hence we
can find two linear
combinations $R_1$ and $R'_1$ of $R$ and $R'$ such that the
corresponding
matrix would be (after a suitable permutation of
 columns) of the form
$$
\left(\matrix 1 & 0 & c \\
0 & 1 & c'\endmatrix\right).
$$
Hence $R_1$ and $R'_1$ give relations stated in lemma. \qed
\enddemo

\proclaim{Lemma 9.3} Let $\rho_1,\rho_2 \in \lt(\bar R)$, 
 $\pi = \rho_1 \rho_2$. Then
$$
u(\rho_1 \subset \pi) \in \Bbb F^\times u(\rho_2 \subset
\pi)+\Cal R_{(\pi)}.
$$
\endproclaim

\demo{Proof} The element $r(\rho_1)r(\rho_2)$ can be
``expanded" in two 
ways: by writing either  $r(\rho_1)$ or $r(\rho_2)$ as an
(infinite) sum of
monomials. Hence
$$
r(\rho_1)r(\rho_2) = c_1 u(\rho_1)r(\rho_2) +
\sum\Sb{\rho_1 \prec \kappa}\endSb c^1 _\kappa\,
u(\kappa)r(\rho_2),
$$
$$
r(\rho_1)r(\rho_2) = c_2 r(\rho_1)u(\rho_2) +
\sum\Sb{\rho_2 \prec \kappa}\endSb c^2 _\kappa\, r(\rho_1)
u(\kappa),
$$
$c_1, c_2 \in \Bbb F^\times$, $c^1_\kappa, c^2_\kappa \in
\Bbb F$,
and the lemma follows. \qed
\enddemo

\remark{Remark}
In Section 8 we have constructed some relations among relations, but 
for the proof of Theorem 6.5.2 we need the relations (9.1) stated in
terms of colored partitions. The first step is accomplished by Lemma 9.2
for ``adjacent'' embeddings like the one we had in the case of $\pi=x(-3)x(-2)x(-1)$,
$\rho_1=x(-3)x(-2)$ and $\rho_2=x(-2)x(-1)$. On the other hand,
Lemma 9.3 settles the problem of ``disjoint'' embeddings like for
$\pi=x(-5)x(-4)x(-2)x(-1)$, $\rho_1=x(-5)x(-4)$ and $\rho_2=x(-2)x(-1)$.
In the case of level 1 modules this is all that can happen, and, apart for
the initial conditions discussed in Section 10, our proof of Theorem 6.5.2
would be finished.

However, in the case when level $k\geq 2$, say $k=2$, we may have embeddings 
which are neither ``adjacent'' nor ``disjoint''. For example, if
$\pi=x(-3)^2 x(-2)^2 x(-1)$,  $\rho_1=x(-3)^2 x(-2)$ and 
$\rho_2=x(-2)^2 x(-1)$, we can not apply Lemmas 9.2 and 9.3 directly.
Of course, we can still apply Lemma 9.2 in two steps by considering an 
``intermediate'' embedding $\rho_3 \subset \pi$ with $\rho_3=x(-3)x(-2)^2$,
and $\rho_1 \subset \pi$ and $\rho_2 \subset \pi$ is just the example of 
what we call linked embeddings. But this kind of argument will not work
for $\pi=x(-3)^2 x(-2)x(-1)^2$, $\rho_1=x(-3)^2 x(-2)$ and 
$\rho_2=x(-2)x(-1)^2$, an example of what we call the exceptional case,
listed in Lemma 7.3.1 within the ``series'' (12) of exceptional cases.

Similar situation appears in the principal picture for modules of level
$k\geq 4$, cf. the proof of Lemma 9.18 in \cite{MP}, where both the
linked embeddings and the exceptional cases are essentially of the
form described above. Here the presence
of three colors $x$, $h$ and $y$ leads to quite a few exceptional cases.
\endremark

We shall prove Theorem 9.1 by induction on length
$\ell(\pi)$. If
 $\ell(\pi) = k+2$, then (9.1) holds by Lemma 9.2.

It is clear that if (9.1) holds for $\rho_1, \rho_2 \subset
\pi$, then it 
holds as well for  $\rho_1, \rho_2 \subset \kappa\pi$ for
any
$\kappa \in \Cal P(\bar B)$. In particular it is enough to
prove (9.1)
for any two embeddings $\rho_1, \rho_2 \subset \rho_1 \cup
\rho_2$.
The case when  $ \rho_1 \cap \rho_2 = \varnothing$ is proved
by Lemma 9.3.

If (9.1) holds for every colored partition of length less
than $s$, then
(almost by definition) (9.1) holds for two linked embeddings
$\rho_1 \subset \pi$,  $\rho_2 \subset \pi$, $\ell(\pi) =
s$.
By the above remarks it is enough to prove (9.1) in the case
when
$\pi = \rho_1 \cup \rho_2$, $ \rho_1 \cap \rho_2 \ne
\varnothing$, 
$k+3 \le \ell(\pi) \le s$, $\rho_1 \subset \pi$ and $\rho_2
\subset \pi$
are not linked, i.e. in the exceptional cases (1)--(14) 
listed in Lemma 7.3.1. In the 
proof of (9.1) for (some of) these exceptional cases we
shall use the following
two lemmas:
\proclaim{Lemma 9.4} Let
$$
u(\rho \subset \pi) = \sum C_{\rho',\pi'}\, u(\rho' \subset
\pi'),
$$
where for each $\pi' \prec \pi$ the term $u(\rho' \subset
\pi')$ appears for
only one embedding $\rho' \subset \pi'$. 
Then $C_{\pi'} = C_{\rho',\pi'} = 0$ for $\pi' \prec \pi$.
In particular
$$
u(\rho \subset \pi) \in \sum C_{\rho',\pi}\, u(\rho' \subset
\pi) +\Cal R_{(\pi)}.
$$
\endproclaim

\demo{Proof} The statement follows from the fact that the
right hand side has
the unique leading term $\pi$ (Proposition 6.4.3). \qed
\enddemo

\noindent {\bf Example.} As an example of the proof of (9.1) in the
exceptional case let us consider $\pi=x(-3)^2 x(-2)x(-1)^2$, 
$\rho_1=x(-3)^2 x(-2)$ and $\rho_2=x(-2)x(-1)^2$. Here we have
only one color $x$ and $\sh \pi=(-3)^2 (-2) (-1)^2$. 
Note that $ (-2)^5 \prec (-3) (-2)^3 (-1) \prec (-3)^2 (-2) (-1)^2$
are all plain partitions of degree $-10$ and $\preccurlyeq \sh\pi$ and
we can visualize them by Young diagrams 
$$
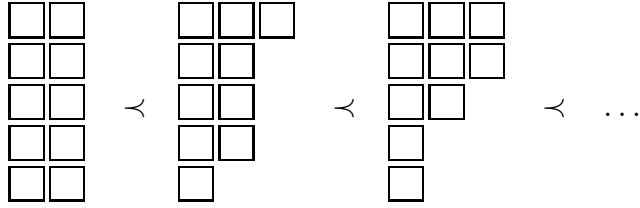
\spreadlines{-1.15ex}\alignedat 6
&\sq\sq \quad   &&\quad  &&\sq\sq\sq \quad      &&\quad  &&\sq\sq\sq  \quad && \\
&\sq\sq \quad   &&\quad  &&\sq\sq \quad         &&\quad  &&\sq\sq\sq  && \\
&\sq\sq \quad   &\prec&\quad  &&\sq\sq \quad    &\prec&\quad  &&\sq\sq &\prec &\quad\dots  \\
&\sq\sq \quad   &&\quad  &&\sq\sq \quad         &&\quad  &&\sq \quad   &&\\
&\sq\sq \quad   &&\quad  &&\sq \quad            &&\quad  &&\sq \quad   && 
\endalignedat
$$
Since we want to prove (9.1), it will be enough to consider terms 
$u(\rho'\subset\pi ')$ with $\pi'\preccurlyeq \pi$, which will amount to
$\sh\pi'\preccurlyeq \sh\pi$, and all other terms will be denoted as $\dots$.
We can also visualize the shapes of the leading terms of
$r_{3\alpha}(-4)=r(x(-2)x(-1)^2)$, $r_{3\alpha}(-5)=r(x(-2)^2x(-1))$, 
$r_{3\alpha}(-6)=r(x(-2)^3)$, $r_{3\alpha}(-7)=r(x(-3)x(-2)^2)$ and
$r_{3\alpha}(-8)=r(x(-3)^2x(-2))$ respectively as
$$
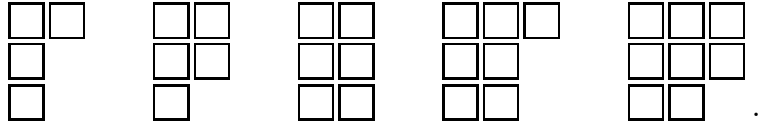
\spreadlines{-1.15ex}\alignedat 5
&\sq\sq \qquad &&\sq\sq\qquad  &&\sq\sq \qquad   &&\sq\sq\sq\qquad  &&\sq\sq\sq   \\
&\sq \quad    &&\sq\sq\quad  &&\sq\sq \quad   &&\sq\sq          &&\sq\sq\sq   \\
&\sq \quad    &&\sq\quad     &&\sq\sq \quad   &&\sq\sq \quad    &&\sq\sq \quad \ .    
\endalignedat
$$
It is clear that these are all possible $\sh\rho'$ for $\rho'\subset\pi '$ with
$\sh\pi'\preccurlyeq \sh\pi$ as above. 

Now we start with $u(\rho_2\subset\pi) =x(-3)^2 r_{3\alpha}(-4)$. First we use
the relation (8.1) with $k=2$ and $n=-7$, i.e.
$$
\sum_{j \in \Bbb Z}\, \bigl(4j+7 \bigr)\, x(j) r_{3\alpha}(-7-j) = 0,
$$
which we write as
$$
x(-3) r_{3\alpha}(-4)
=-\tfrac{1}{5} x(-2) r_{3\alpha}(-5)+\tfrac{3}{5} x(-1) r_{3\alpha}(-6) +\dots,
$$
and we get
$$
x(-3)^2 r_{3\alpha}(-4)
=-\tfrac{1}{5} x(-3)x(-2) r_{3\alpha}(-5)+\tfrac{3}{5} x(-3)x(-1) r_{3\alpha}(-6) +\dots.
$$
Now we apply the relation (8.1) again:
$$\align
&x(-3) r_{3\alpha}(-5)
=0\cdot x(-2) r_{3\alpha}(-6)+\tfrac{4}{4} x(-1) r_{3\alpha}(-7) +\dots,\\
&x(-3) r_{3\alpha}(-6)
=\tfrac{1}{3} x(-2) r_{3\alpha}(-7)+\tfrac{5}{3} x(-1) r_{3\alpha}(-8) +\dots,
\endalign$$
and we get
$$\align
x(-3)^2 r_{3\alpha}(-4)
&=-\tfrac{1}{5}\cdot 0\cdot x(-2)^2 r_{3\alpha}(-6)\\
&+(\tfrac{3}{5}\cdot\tfrac{1}{3}-\tfrac{1}{5} \cdot\tfrac{4}{4})\cdot 
x(-2)x(-1) r_{3\alpha}(-7) \\
&+\tfrac{3}{5}\cdot\tfrac{5}{3}\cdot x(-1)^2 r_{3\alpha}(-8)+\dots\\
&=x(-1)^2 r_{3\alpha}(-8) +\dots .
\endalign$$
Hence (9.1) holds. 

In the above proof we can avoid the explicit calculation of coefficients:
it is sufficient to know that, after applying the relation 
(8.1) in two steps on terms of the form $x(-3) r_{3\alpha}(n)$, we get a relation
of the form
$$\align
u( \rho_2 \subset \pi)
&=x(-3)^2 r_{3\alpha}(-4)\\
&=C_1 x(-2)^2 r_{3\alpha}(-6) \\
&+C_2 x(-2)x(-1) r_{3\alpha}(-7)\\
&+C_3 x(-1)^2 r_{3\alpha}(-8)+\dots\\
&=C_1  u( x(-2)^3 \subset x(-2)^5 )\\
&+C_2  u( x(-3)x(-2)^2 \subset x(-3)x(-2)^3 x(-1))\\
&+C_3  u( \rho_1 \subset \pi)+\dots .
\endalign$$
Then Lemma 9.4 implies $C_1=C_2=0$. After that our choice of
$r(\pi_{3\alpha}(n))= r_{3\alpha}(n)$ also implies $C_3=1$. 
A general case of this type is treated in
the proof of exceptional case (10).

As it was mentioned in the Introduction, and illustrated by the previous example, 
in some cases Lemma 9.4 enables us to prove (9.1) by 
first expressing $u(\rho_2\subset\pi)-u(\rho_1\subset\pi)$ as a linear combination
of terms  $u(\rho'\subset\pi ')$, making sure afterwards that only terms
for $\pi'\succ\pi$ appear nontrivially. In doing this we shall 
have to keep track of all
the terms that may appear, and the following lemma will be used
to check that the expression we obtain contains for each $\pi' \prec \pi$ 
at most one term of the form $u(\rho' \subset\pi')$. 

\proclaim{Lemma 9.5} Let $\pi$ be a colored partition. Let
$\pi_1$ be another
colored partition such that $\sh \pi_1 = \sh \pi$, $\wt
\pi_1 = \wt \pi$. Let 
$\rho \subset \pi_1$ be an embedding.
\roster
\item[2]
If $\pi = x(j-1)^{k+1-a-b} y(j)^a h(j)^b x(j)^a$,  $a \ge
2$,  
$a+b \le k$, and if \newline
 $\pi_1= y(j-1)^f h(j-1)^e 
y(j)^d h(j)^c \rho$, then $y(j)$ is not a part of 
$\rho $.
\item[3]
If $\pi=h(j-1)^{k+1-a-b} x(j-1)^a y(j)^b x(j)^{k+1-a}$, \
$a\ge 1$, and if \newline
$\pi_1 = y(j-1)^c h(j-1)^d x(j-1)^e y(j)^f h(j)^g \rho$,
then $y(j)$ is
not a part of $\rho$.
\item[4] 
If $\pi=h(j-1)^{k+1-a-b} x(j-1)^a y(j)^b h(j)^{k+1-a-b}$, 
$1 \le a+b \le k$, and if 
$x(j)$ is not a part of $\pi_1$, then $y(j-1)$ is
not a part of $\rho$.
\item[5] 
If $\pi=y(j-1)^{k+1-a} x(j-1)^b y(j)^a h(j)^{k+1-a-b}$, \
$a\ge 1$,  and if \newline
$\pi_1 = h(j-1)^g x(j-1)^f y(j)^e h(j)^d x(j)^c \rho$, then
$x(j-1)$ is
not a part of $\rho$.
\item[6]
If  $\pi = y(j-1)^a h(j-1)^b x(j-1)^a y(j)^{k+1-a-b}$, $a
\ge 2$, $a+b \le k$,
and if $\pi_1 = h(j-1)^c x(j-1)^d h(j)^e x(j)^f \rho$, then
$x(j-1)$ is not 
a part of $\rho$. 
\item[11]
If $\pi = x(j-2)^{k+1-a-b} h(j-1)^a x(j-1)^b y(j)^c
h(j)^{k+1-b-c}$, 
\ $a,b \ge 1$, \newline
$a+b+c \le k$, and if $\pi_1\prec \pi$, then 
$ y(j-1)$ is not a part of $\rho$.

\item[12]
If $\pi = x(j-2)^{k+1-a-b} h(j-1)^a x(j-1)^b h(j)^c
x(j)^{k+1-b-c}$, 
\ $a,b \ge 1$, \newline
$a+b \le k$, and if $\pi_1\prec \pi$, then 
$ y(j-1)$ is not a part of $\rho$.
\endroster
\endproclaim

\remark{Remark} Note that (2) and (3) have dual statements 
(6)
and (5) (respectively) and that the dual proofs hold. Also
note 
that the statement (4) is not self dual.
\endremark
\demo{Proof}
(2) The statement is clear since $\rho = \dots y(j)^p
h(j)^q$ for
$p \ge 1$, $q \ge 0$, \newline
$\sh \rho = (j-1)^{k+1-p-q} j^{p+q}$, would imply
$\wt (\pi_1) < (k+1-a-b)\alpha = \wt (\pi)$.

(3) The statement is clear since
$\rho = \dots y(j)^p h(j)^q$  for $p \ge 1$, $q \ge 0$,
\newline 
$\sh \rho = (j-1)^{k+1-p-q} j^{p+q}$, would imply
$\wt (\pi_1) < (k+1-b)\alpha = \wt (\pi)$.

(4) By assumptions $\pi_1$ must be of the form
$y(j-1)^c h(j-1)^e x(j-1)^f y(j)^d h(j)^g$,
where $c+e+f = k+1-b$, $\wt (\pi) = (a - b)\alpha = (-c -d
+f)\alpha$. Let
$$
\rho =y(j-1)^{c'} h(j-1)^{k+1-c'-d'} y(j)^{d'} \subset \pi_1
$$
for some $c' \ge 1$. Then $c\ge c' \ge 1$, $e \ge
k+1-c'-d'$, 
$d \ge d'$, and in particular $c+e \ge k+1-d'$. But then
$d < d+a+c = f+b = k+1 -(c+e) \le d' \le d$,
a contradiction.

(6) The statement is clear since  $\rho = h(j-1)^q
x(j-1)^p\cdots$
for $p \ge 1$, $q \ge 0$, $\sh \rho = (j-1)^{p+q}
j^{k+1-p-q}$, would imply
$\wt (\pi_1) > -(k+1-a-b)\alpha = \wt (\pi)$.

(11) Assume that $y(j-1)^d \subset \pi_1$ for some $d \ge
1$.
Since by assumption $\pi_1 \prec \pi$, $\wt \pi_1 = \wt
\pi$, we have that
$\pi_1$ contains
$y(j)^{c+i} h(j)^{k+1-b-c-i}$
for some $i \ge 1$. Because of $\wt \pi_1 = \wt \pi$ we also
have that $\pi_1$
contains either
$$
x(j-2)^{k+1-a-b} y(j-1)^s h(j-1)^{a-i-2s} x(j-1)^{b+i+s},\ 
s \ge d, \quad\text{or}
$$
$$\gathered
y(j-2)^r h(j-2)^p x(j-2)^{k+1-a-b-r-p} \\ 
\cdot y(j-1)^{s-r} h(j-1)^{a-i-2s-p} x(j-1)^{b+i+s+r+p},\  
s-r \ge d.
\endgathered$$
Now the assumption $\rho \subset \pi_1$, $\rho$ of the form 
$\rho = y(j-1)^{d'} h(j-1)^{e'} y(j)^{k+1-e'-d'}$
would imply $k+1-e'-d' \le c+i$ and
$$
d' + e' \le s+(a-i-2s) \quad \text{or}\quad d'+e' \le
(s-r)+(a-i-2s-p),
$$
which in turn would imply
$k+1 \le c+a \le c+a+b$, contrary to our assumption.

The possible form for $\pi_1$ shows that $x(j-1) \subset
\pi_1$, so there is
no embedding of the form
$\rho \subset \pi_1, y(j-1)\subset \rho, \sh \rho =
(j-2)^{k+1-a-b}(j-1)^{a+b}$.

(12) Assume that $y(j-1)^d \subset \pi_1$ for some $d \ge
1$. Since by
assumption $\pi_1 \prec \pi$, $\wt \pi_1 = \wt \pi$, we have
that $\pi_1$ contains
$y(j)^i h(j)^{c-i+m} x(j)^{k+1-b-c-m}$
for some $i,m \ge 0$, $i+m \ge 1$. Because of $\wt
\pi_1=\wt\pi$ we also 
have that
$\pi_1$ contains either
$$
x(j-2)^{k+1-a-b} y^s h(j-1)^{a-2s-m-i} x(j-1)^{b+m+i+s}, \ s
\ge d, \quad \text{or}
$$
$$\gathered
y(j-2)^r h(j-2)^p x(j-2)^{k+1-a-b-r-p}  \\ 
\cdot y(j-1)^{s-r} h(j-1)^{a-2s-m-i-p} x(j-1)^{b+m+i+s+r+p},
\ s-r \ge d.
\endgathered$$

Now the assumption $\rho \subset \pi_1$, $\rho$ of the form
$\rho = y(j-1)^{d'} h(j-1)^{e'} y(j)^{k+1-e'-d'}$
would imply $k+1-e'-d' \le i$ and
$$
d'+e' \le s+(a-2s-m-i) \quad \text{or}
\quad d'+e' \le  (s-r)+(a-2s-m-i-p),
$$
which in turn would imply $k+1 \le a \le a+b$,
contrary to our assumption.

Since $\wt\pi_1 = \wt \pi$ and $x(j)^f\subset\pi_1$ for
$f\le k+1-b-c$,
it must be that $x(j-1)^b\subset \pi_1$. Hence there is no
embedding of the
form
$$
\rho \subset \pi_1, \quad y(j-1) \subset\rho,\quad \sh \rho
= (j-2)^{k+1-a-b}(j-1)^{a+b}.
\qed$$
\enddemo

By using relations among relations constructed in Section 8
we shall prove 
(9.1) by induction on length $\ell(\pi)$ for all exceptional
cases. Note
that in the exceptional cases we have only two embeddings 
$\rho_1\subset\pi$, $\rho_2\subset\pi$.
At some places in the proof it will be convenient to write
$$
a \underset \,{\text{sh}}\to \sim  b 
\qquad \text{or}\qquad
a - b \underset \,{\text{sh}}\to \sim  0
$$
if $\sh\lt (a)= \sh\lt (b)$ and
$$
a-b \in  \overline {\lspan\{u(\rho \subset \pi)\mid
 \rho \in \lt (\bar{R}), \pi \in (\lt (\bar{R})),
 \rho \subset \pi, \sh \pi \succ \sh\lt (a)\}}.
$$
At some places in the proof we shall use relations (8.5), $n
= (k+2)j$,
 $\vert i \vert \le k+1$, written in the form
$$\split
(k\!+\!2\!-\!i) x(j) r_{(i-1)\alpha}((k+1)j)&+ i h(j)r_{i\alpha}((k+1)j)\\
 &+(k\!+\!2\!+\!i) y(j) r_{(i+1)\alpha}((k+1)j) \underset
\,{\text{sh}}\to\sim \ 0. 
\endsplit\tag{9.5}
$$

{\bf Now we consider all exceptional cases:}

$$
\align 
&\pi = y(j)^a h(j)^{k+1-a} x(j)^a, \quad 2 \le a \le k,
\tag"{{\bf (1\,)}}"\\
&\rho_1 = y(j)^a h(j)^{k+1-a}.
\endalign
$$
By using relations (8.5), $n = (k+2)j$, \ $i \in [-(a-1),
a-1]$, of the form (9.5)
we get
$$
x(j)^a r(\rho_1) \underset \,{\text{sh}}\to\sim
 \sum\Sb{b, c\ge 0}\\{b+c=a}\endSb C_{bc}\, y(j)^b h(j)^c
r_{b\alpha} ((k+1)j).
$$
Now (9.1) follows by applying Lemma 9.4.

$$
\align &\pi =  x(j-1)^{k+1-a-b} y(j)^a h(j)^bx(j)^a, \quad
a\ge 2, 
a+b \le k,  \tag"{{\bf (2\,)}}"\\
&\rho_1 =  x(j-1)^{k+1-a-b} y(j)^a h(j)^b.
\endalign
$$

By using in $a$ steps relations (8.6) of the form $C_x=1,
C^x =0$, or (9.5), we get
$$
x(j)^a r(\rho_1) \underset \,{\text{sh}}\to\sim
  \sum C_{cdef}\,y(j-1)^f h(j-1)^e y(j)^d h(j)^c r(\rho),
\tag{9.6}
$$
where $c+d+e+f = a$, $\rho \in \lt(\bar R)$. (Note
that the leading term $\rho$ of $r(\rho)$ is completely
determined by
$c,d,e,f$.)

Note that each time when we apply one of the relations (8.6)
or (9.5)
we replace $x(j) r(\rho)$ by a combination of terms of the
form 
$X r(\rho')$, where $X$ is either $y(j), h(j), y(j-1)$ or
$h(j-1)$ and
$\wt \rho' \le \wt \rho + 2\alpha $. Hence after applying
relations
 $i$ times, $i < a$, all $r(\rho)$'s that appear have
weights bounded by
$(k+1-2a-b)\alpha + 2 i\alpha < (k+1-b)\alpha$ and we see
that there exists a relation
of the form (8.6) or (9.5) which could be applied again. (In
the subsequent 
cases we shall state the type of relations we use and we
shall omit the proof
 that this can be done.)

By Lemma 9.5(2) $y(j)$ is not a part of $\rho$, i.e.
$$
r(\rho) = r(x(j-1)^p h(j)^q x(j)^{k+1-p-q}).
$$
Now we can apply Lemma 9.4 and (9.1) follows since
$r(\rho_1)$ does not
appear on the right hand side of (9.6).

$$\align
&\pi = h(j-1)^{k+1-a-b} x(j-1)^a y(j)^b x(j)^{k+1-a},\quad 
1 \le a \le k-1, \tag"{{\bf (3\,)}}"\\
&\rho_1 = h(j-1)^{k+1-a-b} x(j-1)^a y(j)^b.
\endalign$$
By using in \ $k+1-a$ \ steps relations: a) relations (8.6)
of the form
$C_x=1, C_h = 0$ on the terms $x(j) r(\rho)$, $\sh \rho
=(j-1)^p j^{k+1-p}$
for $p \ge 1$, b) relations (9.5) on the terms $x(j)
r(\rho)$,
$\sh \rho = j^{k+1}$, $\wt\rho < k+1$, c) no relation
otherwise, we get
$$\gathered
 x(j)^{k+1-a} r(\rho_1) \underset \,{\text{sh}}\to\sim
  \sum C_{cdef}\,y(j-1)^c h(j-1)^d x(j-1)^e y(j)^f r(\rho)
\\
+ \sum C_{cdefgh}\, y(j-1)^c h(j-1)^d x(j-1)^e y(j)^f
h(j)^g x(j)^h r(\rho).
\endgathered \tag{9.7}$$
In the first sum in (9.7) we have $\sh \rho = (j-1)^p
j^{k+1-p}$, $p \ge 1$,
and Lemma 9.5(3) implies that
$$
r(\rho) = r(x(j-1)^p h(j)^q x(j)^{k+1-p-q}).
$$
In the second sum in (9.7) we have $\sh \rho = j^{k+1}$. If
$h \ge 1$ then
$\rho = x(j)^{k+1}$ since in last steps we used c), i.e. we
didn't use any
relation. If $h=0$ then Lemma 9.5(3) implies that
$$
r(\rho) = r(h(j)^q x(j)^{k+1-q}).
$$
Since the number of $y(j)$'s in terms of the first sum is
greater than in 
terms of the second
sum, we see that we can apply Lemma 9.4. Since $r(\rho_1)$
does not appear
on the right hand side of (9.7), we get (9.1).

$$\align
&\pi = h(j-1)^{k+1-a-b} x(j-1)^a y(j)^b h(j)^{k+1-a-b}, \  1
\le a+b \le k-1,
\tag"{{\bf (4\,)}}"\\
&\rho_1 = x(j-1)^a y(j)^b h(j)^{k+1-a-b}.
\endalign$$
By using in \ $k+1-a-b$ \ steps relations (8.6) of the form
$C^h = 1, C^x = 0$
 we get
$$
 h(j-1)^{k+1-a-b} r(\rho_1) \underset \,{\text{sh}}\to\sim
   \sum C_{cdef}\, y(j-1)^c y(j)^d  h(j)^e x(j)^f r(\rho).
\tag{9.8}$$
We need not consider terms with $f > 0$ (since for those the
leading terms are
$\succ \pi)$. Consider two classes of terms with a leading
term $\pi_1$: (A) 
$y(j-1)^c \subset \pi_1$ for some $c > 0$ and (B) terms of
the form
$r(\rho) y(j)^d h(j)^e$, $x(j)$ is not a part of $\rho$
(since otherwise
$\pi_1 \succ \pi)$. By Lemma 9.5(4) and $\wt (\pi) = \wt
(\pi_1)$, 
in the case (B) we have either
$$
\align
\rho& = h(j-1)^{k+1-b-p} x(j-1)^p y(j)^b , \quad 0 < p \le
k-b,\quad \text{or}\\
\rho& = x(j-1)^{k+1-b}y(j)^b, \quad \pi_1 =  x(j-1)^{k+1-b}
y(j)^{k+1-a}.
\endalign
$$
In either case for given $\pi_1$ there is only one
embedding. In the case (A)
we may apply the induction hypothesis on $\pi_1/y(j-1)^c$. 
Hence we may apply
 Lemma 9.4 and (9.1) follows.
\remark{Remark}
The reader may find a ``dual" proof of case (4) a bit
shorter, but some changes
are needed. For example, the argument dual to:
``We need not consider terms with $f > 0$ (since for those
the leading terms are
$\succ \pi)$" is false. However, in some cases (like (5) and
(6)) dual arguments
hold all the way through the proof.
\endremark

\bigskip

\noindent {\bf (5\,)} The proof is dual to the proof of (3).

\bigskip

\noindent {\bf (6\,)} The proof is dual to the proof of (2).

\bigskip

\noindent {\bf (7\,)} \ Let $\pi = y(j-1)^{k+1-a-c} h(j-1)^c
x(j-1)^b y(j)^a
 h(j)^{k+1-a-b}$, $a \ge 1$, $a+b+c \le k$.
Note that we have already considered the case $c=0$ in (5).
We keep the 
shape
$$
(j-1)^{k+1+b-a} j^{k+1-b}
$$
fixed and prove by induction on $c$, $0 \le c \le k-a$, that
embedding of the
 form (for $c \ge 1$)
$$
y(j-1)^{k+1-a-c} h(j-1)^c r(x(j-1)^b y(j)^a h(j)^{k+1-a-b})
\tag{9.9}
$$
is, modulo terms with the leading term $\succ \pi$,
proportional to
$$
x(j-1)^b h(j)^{k+1-a-b} r(y(j-1)^{k+1-a-c} h(j-1)^c y(j)^a).
 \tag{9.10}
$$
By the induction hypothesis the embedding of the form
$$
y(j-1)^{k+2-a-c} h(j-1)^{c-1} r(x(j-1)^b y(j)^a
h(j)^{k+1-a-b}) \tag{9.11}
$$
with the leading term $\pi'$ is proportional to 
$$
x(j-1)^b h(j)^{k+1-a-b} r(y(j-1)^{k+2-a-c} h(j-1)^{c-1}
y(j)^a)  \tag{9.12}
$$
modulo terms with the leading term $\succ \pi'$.

Now we apply the adjoint action of $x$ on the relation
involving (9.11) and
(9.12).

The action of $x$ on (9.11) gives three terms, one is
(proportional to)
(9.9), the other two have the leading terms bigger than
$\pi$. The action
of $x$ on (9.12) gives (9.10) and another term $\succ \pi$.

What remains to show is that the action of $x$ on any term
(in the relation 
involving (9.11) and (9.12)) with the leading term $\pi''
\succ \pi'$ will give
terms $\succ \pi$.

If $\pi''$ contains $x(j)$, then all terms obtained by the
action of $x$
contain $x(j)$, the result being $\succ \pi$.

If $\pi''$ contains $y(j)^{a-d} h(j)^{k+1-a-b+d}$, \ $d \ge
1$, then the
action of $x$ leaves this term, gives one $x(j)$ more and
one $h(j)$ less, or
gives one $h(j)$ more and one $y(j)$ less. In either case
the result is
$\succ \pi$.
 
If $\pi''$ contains $x(j-1)^{b+d} y(j)^a h(j)^{k+1-a-b}$, \ 
$d \ge 0$, then the action of $x$ leaves this term or gives
one $x(j)$, 
$h(j)$ or $x(j-1)$ more. In either case the result is $\succ
\pi$ when $d \ge 1$.
In the case $d=0$ we have $\pi'' = \pi'$.

\bigskip

\noindent {\bf (8\,)}
\ Let $\pi=h(j-1)^{k+1-a-b} x(j-1)^a y(j)^b h(j)^c
x(j)^{k+1-a-c}$,
$a \ge 1$, $a +b + c \le k$.
Note that we have already considered the case $c=0$ in (3).
We keep the 
shape 
$$
(j-1)^{k+1-b} j^{k+1+b-a}
$$
fixed and prove by induction on $c$, \ $0 \le c \le k - a -
b$, that 
embedding of the form (for $c \ge 1$)
$$
h(j)^c x(j)^{k+1-a-c} r(h(j-1)^{k+1-a-b} x(j-1)^a y(j)^b)
\tag{9.13}
$$
is, modulo terms with the leading term $\succ \pi$,
proportional to 
$$
h(j-1)^{k+1-a-b} y(j)^b r(x(j-1)^a h(j)^c x(j)^{k+1-a-c}).
\tag{9.14}
$$
By induction hypothesis the embedding of the form
$$
h(j)^{c-1} x(j)^{k+2-a-c} r(h(j-1)^{k+1-a-b} x(j-1)^a
y(j)^b) \tag{9.15}
$$
with the leading term $\pi'$ is proportional to 
$$
h(j-1)^{k+1-a-b} y(j)^b r(x(j-1)^a h(j)^{c-1}
x(j)^{k+2-a-c}) \tag{9.16}
$$
modulo terms with the leading term $\succ \pi'$. Now we
apply the adjoint action of $y$
on the relation involving (9.15) and (9.16).

The action of $y$ on (9.15) gives three terms, one is
(proportional to) (9.13), 
the other two (with the leading terms) $\succ \pi$. The
action of $y$ on (9.16)
gives (9.14) and another term $\succ \pi$.

What remains to show is that the action of $y$ on any term
(in the relation 
involving (9.15) and (9.16)) with the leading term $\pi''
\succ \pi'$ will give
terms $\succ \pi$.

First let $\pi'' \supset y(j)^b h(j)^{c-1} x(j)^{k+2-a-c}$,
i.e.
$$
\pi'' = y(j-1)^d h(j-1)^{k+1-a-b-2d} x(j-1)^{a+d}
y(j)^b h(j)^{c-1} x(j)^{k+2-a-c},
$$
where $d \ge 1$ and $\rho \subset \pi''$. If $x(j) \notin
\rho$, the partition
(whith the obvious notation) $(y \cdot \rho) (\pi''/\rho)$
contains
$x(j)^{k+2-a-c}$, so is $\succ \pi$. If $y$ does not act on
any
$x(j) \in \pi''/\rho$, we again obtain a term $\succ \pi$.
If $y$ acts on $x(j)$, we
obtain the term (with the leading term):
$$
y(j-1)^d h(j-1)^{k+1-a-b-2d} x(j-1)^{a+d}
y(j)^b h(j)^c x(j)^{k+1-a-c} \succ \pi. \tag{9.17}
$$
If $x(j) \in \rho$, then we see that  $(y \cdot \rho)
(\pi''/\rho)$
is of the form (9.17). Also we see that the action of $y$ on
$\pi''/\rho$
will give terms $\succ \pi$.

Now let $\pi'' \supset y(j)^{b-d} h(j)^{c-1+d}
x(j)^{k+2-a-c}$, where
$d \ge 1$ and $\rho \subset \pi''$. If $x(j) \notin \rho$,
then
 $(y \cdot \rho) (\pi''/\rho)$ contains $x(j)^{k+2-a-c}$, so
is
$\succ \pi$. If $y$ does not act on $x(j) \in \pi''/\rho$,
we again obtain a term
$\succ \pi$. If $y$ acts on $x(j)$, we obtain the term
(which contains)
$$
\cdots  h(j)^{c-1+d} h(j) x(j)^{k+1-a-c} \succ \pi.
\tag{9.18}
$$
If $x(j) \in \rho$, then $y \cdot \rho$ has one less $x(j)$
and one more
$h(j)$, so $(y \cdot \rho) (\pi''/\rho)$ is of the form
(9.18). Also we 
see that the action of $y$ on $\pi''/\rho$ will give terms
$\succ \pi$.

Finally, if $\pi''\supset x(j)^{k+3-a-c}$, then the action
of $y$ will give
terms that contain $ x(j)^{k+2-a-c}$ and hence $\succ \pi$. 

\bigskip

\noindent {\bf (9\,) and (10\,).} \ We fix a shape $(j-2)^a
(j-1)^{k+1-a} j^a$,  
$2 \le a \le k$. Let
$$\align
&\pi = y(j-2)^a y(j-1)^{k+1-a} y(j)^a,\\
&\rho_1 = y(j-1)^{k+1-a} y(j)^a. \endalign
$$
By using relations (8.4), $i=-k-2$, of the form
$$
y(j-2) r_{-(k+1)\alpha} (m) + c_1 y(j-1) r_{-(k+1)\alpha} (m-1) +
c_2 y(j) r_{-(k+1)\alpha} (m-2) + \cdots = 0,
$$
where $m=(k+1)(j-1) \pm r$, $0 \le r < k+1$, $c_1, c_2$ are
some rational
constants and $\cdots$ denotes the sum of shorter elements
(i.e. elements of length 
$k+1$) and terms which contain a factor $y(p)$, $p > j$ or
$p < j-2$, we get
$$
y(j-2)^a r(\rho_1) = \sum C_{bc} \,y(j)^b y(j-1)^c r(\rho)+
A + B. \tag{9.19}
$$
Here $b+c=a$; the terms in $A$ contain a factor $y(p)$, $p >
j$ 
(and their leading terms are greater than $\pi$); the terms
$B$ contain a factor
$y(p)$, $p < j-2$.

If $\pi_1$ is a leading term of a summand in $B$, $\pi_1
\prec \pi$, then for any
embedding $\rho \subset \pi_1$ the term $y(p)$, $p < j-2$,
is not a part of $\rho$. 
By induction
 hypothesis (on lenght) for two embeddings in such $\pi_1$
(i.e. in $\pi_1/y(p)$)
relation (9.1) holds, so we may rewrite $B$ so that for each
$\pi_1 \prec \pi$
that appears there is only one embedding.

In (9.19) we have that $r(\rho)=r(y(j-2)^b y(j-1)^{k+1-b})$,
so each leading 
term appears with only one embedding. Hence we can apply
Lemma 9.4 and we
get
$$
y(j-2)^a r(\rho_1) \in \Bbb F^\times y(j)^a r( y(j-2)^a
y(j-1)^{k+1-a})+ \Cal R_{(\pi)}.
$$
But $\pi' \succ \pi$ implies $\sh(\pi')\succ(j-2)^a
(j-1)^{k+1-a} j^a$,
and hence
$$
y(j-2)^a r(y(j-1)^{k+1-a} y(j)^a) \underset
\,{\text{sh}}\to\sim 
 y(j)^a r(y(j-2)^a y(j-1)^{k+1-a}). \tag{9.20}
$$
Since the adjoint action of $x$ does not change the shape,
we get (9.1) in
the case (9) and (10) by adjoint action of $x$.

\bigskip

\noindent {\bf (11\,) and (12\,).}
\  Let $a\ge 0$, $b \ge 1$  and
$$
\align
&\pi  = x(j-2)^{k+1-a-b} h(j-1)^a x(j-1)^b P,\\
&\rho_1  = x(j-1)^b P,
\quad\text{where}\\
&P = y(j)^c h(j)^{k+1-b-c},\quad a+b+c \le k,\quad
\text{or}\\
&P = h(j)^c x(j)^{k+1-b-c},\quad a+b \le k.
\endalign
$$
Note that the case $a=0$ follows from (9.20) by adjoint
action of $x$.
By using in $a$ steps relations (8.6) of the form $C^h = 1,
\ C^x = 0$ we get
$$\gathered
x(j-2)^{k+1-a-b} h(j-1)^a r(\rho_1) \underset
\,{\text{sh}}\to\sim \\
\sum C_{defg} \,x(j-2)^{k+1-a-b} y(j-1)^d  y(j)^e h(j)^f
x(j)^g r(\rho),
\endgathered \tag{9.21}
$$
where $0 \le d \le a$. It will be enough to consider only
the terms with the 
leading term $\pi_1 \prec \pi$.

Consider three disjoint classes of $\pi_1 \prec \pi$ : (A)
$y(j-1)^d \subset \pi_1$
for some $d \ge 1$, (B) $\pi_1$ is not in A and 
$y(j-2)^p h(j-2)^q x(j-2)^{k+1-a-b-p-q} \subset \pi_1$ for
some
$p+q \ge 1$ and (C) $\pi_1$ is not in A and 
$ x(j-2)^{k+1-a-b} \subset \pi_1$.

Let $\pi_1$ be in the class A.

For two embeddings $\rho, \rho' \subset \pi_1$ (i.e.
 $\rho, \rho' \subset \pi_1/ y(j-1)^d$ by Lemma 9.5) by
induction hypothesis on
lenght terms $u(\rho \subset \pi_1)$ and  $u(\rho' \subset
\pi_1)$
are proportional modulo higher terms (with leading terms
$\succ \pi_1$). These
higher terms are again in class A.

Let $\pi_1$ be in the class B. 

For two embeddings  $\rho, \rho' \subset \pi_1$
(  $\rho, \rho' \subset \pi_1 / y(j-2)^p h(j-2)^q
x(j-2)^{k+1-a-b-p-q}$
because of Lemma 9.5), by induction hypothesis on lenght,
terms 
 $u(\rho \subset \pi_1)$ and  $u(\rho' \subset \pi_1)$ are
proportional
modulo higher terms (with leading terms $\succ \pi_1$).
These higher terms are
again in class B.

Now consider a term in (9.21) which has the leading term
$\pi_1$ in the class
C. Then $d=0$. For $\rho \subset \pi_1$ (i.e.
$\rho \subset \pi_1/y(j)^e h(j)^f x(j)^g$) the term $u(\rho
\subset \pi_1)$
is (by using relations of the form (10), (11) or (12))
proportional to 
$$
y(j)^{e'} h(j)^{f'} x(j)^{g'} r(x(j-2)^{k+1-a-b}
h(j-1)^{a'}x(j-1)^{b'}) 
$$
modulo terms with leading terms $\pi' \succ \pi_1$. 

Hence in a finite number of steps we can rewrite (first the
terms
in class C, and then the terms in class A and B)
the right hand side of (9.21), modulo $\Cal R_{(\pi)}$, as a
linear
combination of terms  $u(\rho \subset \pi_1)$, where for
each $\pi_1\prec\pi$ there
 is only one embedding. By applying Lemma 9.4  we get (9.1).

\bigskip

\noindent {\bf (13\,) and (14\,).} \  Let  $a\ge 0$, $b \ge
1$  and

$$\align
&\pi  = P y(j-1)^b h(j-1)^a y(j)^{k+1-a-b},\\
&\rho_1  = P y(j-1)^b,\quad \text{where}\\
&P  = h(j-2)^{k+1-b-c} x(j-2)^c,\quad a+b+c\le k, 
\quad \text{or}\\
&P  =  y(j-2)^{k+1-b-c} h(j-2)^c,\quad a+b \le k.
\endalign$$
Note that the case $a=0$ follows from (9.20) by adjoint
action of $x$.
Let us also remark that $\sh \pi_1 = \sh \pi$, $\pi_1 \prec
\pi$,
implies that for some $i \ge 0$
$$
\pi_1 \supset  y(j-1)^{b+i} h(j-1)^{a-i} y(j)^{k+1-a-b} =
\rho'.
$$
By using in $a$ steps the relation (8.6) of the form $C_h =
1$, $C_y = 0$ we get
$$\gathered
h(j-1)^a y(j)^{k+1-a-b} r(\rho_1)  \underset
\,{\text{sh}}\to\sim  \\
\sum C_{defg} \, y(j-2)^d h(j-2)^e x(j-2)^f x(j-1)^g
y(j)^{k+1-a-b} r(\rho),
\endgathered \tag{9.22}$$
where $0 \le g \le a$. For $g > 0$ the leading term $\pi_1$
is
greater than $\pi$. By the previous remark in the case $g=0$
we 
may use the induction hypothesis and express
terms in (9.22), modulo $\Cal R_{(\pi)}$, as a linear
combination of elements of the form
$$
 y(j-2)^{d'} h(j-2)^{e'} x(j-2)^{f'} r(\rho'),
$$
where $d'+e'+f'=k+1 - b$. Then each leading term $\pi_1
\prec \pi$ appears with
only one embedding and we may apply Lemma 9.4.

\heading 10. Linear independence of bases of standard
modules
\endheading

Let $\Lambda = k_0 \Lambda_0 + k_1 \Lambda_1$, $k_0, k_1 \in
\Bbb N$, and
$k=k_0 +k_1$.
Let $v_\Lambda$ be a highest weight vector of Verma module
$M(\Lambda)$. Let
 $M^1(\Lambda)$ be the maximal submodule of $M(\Lambda)$.

 Recall that we may choose a map 
$$
\lt(\bar{R}v_\Lambda) \rightarrow \bar R, \qquad
 \rho \mapsto r(\rho)
$$
such that $\lt(r(\rho)v_\Lambda) = \rho$. For
 $\rho \in \lt(\bar{R}v_\Lambda)$ and $\pi \in
(\lt(\bar{R}v_\Lambda))$,
$\rho \subset \pi$, we shall say that $\rho \subset \pi$ is
an {\it embedding}. For
an embedding $\rho \subset \pi$ we set 
$$
u(\rho \subset \pi) = u(\pi/\rho)r(\rho).
$$
Define a filtration ($M^1_{(\pi)} \mid \pi \in
(\lt(\bar{R}v_\Lambda)$) on
$M^1=M^1(\Lambda)=U(\tilde{\goth g})\bar{R}v_\Lambda=
U(\tilde{\goth n}_-)\bar{R}v_\Lambda$  by
$$
M^1_{[\pi]}=\lspan\{urv_\Lambda \mid u \in U(\tilde{\goth
n}_-), r \in \bar R , 
\lt(urv_\Lambda) \succcurlyeq \pi\},
$$
$$
M^1_{(\pi)} = \bigcup_{\pi'\succ \pi} M^1_{[\pi]}.
$$
(At this point this filtration should not be confused with
the filtration \newline
$(M^1 \cap M(\Lambda)_{[\pi]} \mid \pi \in (\lt(\bar{R}))$.)

The main result of this section is an analogue of Theorem
9.1:
\proclaim{Theorem 10.1} Let $\kappa \in
(\lt(\bar{R}v_\Lambda))$,
$\rho_1, \rho_2 \in \lt(\bar{R} v_\Lambda)$, $\rho_1, \rho_2
\subset
\kappa$. Then
$$
u(\rho_1 \subset \kappa)v_\Lambda \in \Bbb F^\times u(\rho_2
\subset \kappa)v_\Lambda
+ M^1 _{(\kappa)}. \tag{10.1}
$$
\endproclaim
\remark{Remark}
Note that Theorem 6.5.2 is a simple consequence of Theorem
10.1: Since \linebreak
$M^1(\Lambda)=U(\tilde{\goth n}_-)\bar{R}v_\Lambda $, 
the set of vectors
$$
u(\rho \subset \pi)v_\Lambda, \qquad  \rho \in \lt(\bar{R}
v_\Lambda),
 \pi \in (\lt(\bar{R} v_\Lambda)), \rho \subset \pi,
$$
is a spanning set of $M^1(\Lambda)$. Using Theorem 10.1 and
induction we see
 that (6.5.1) is a spanning set as well. Clearly this set is
linearly independent,
and Theorem 6.5.2 follows.
\endremark

The rest of this section is devoted to the proof of Theorem
10.1.

We may assume that the map $\rho \mapsto r(\rho)$ is such
that $r(\rho)$ is 
homogeneous of weight $\wt (\rho)$ and degree $|\rho |$. We
will construct
such a map $r \: \lt(\bar{R} v_\Lambda) \rightarrow \bar R$
by starting with a map
(denoted by the same $r$)
$$
r \: \lt(\bar{R}) \rightarrow \bar R, \qquad \rho \mapsto
r(\rho)
$$
defined by $r(\pi_{m\alpha}(n)) = r_{m\alpha}(n)$ for $n \in \Bbb Z$, $|m|
\le k+1$
(cf. Section 6.6).

\proclaim{Proposition 10.2} Let $ a,b \in [0,k+1]$,
 $a+b \le k+1$. Then 
$$\align
&r(x(-1)^a h^{k+1-a-b} x^b) v_\Lambda = 0 \quad \text{if}
\quad b > 0,
\tag{a} \\
&r(x(-1)^a y^b h^{k+1-a-b}) v_\Lambda = 0 \quad \text{if}
\quad
 a \le k_0 \ \text{and} \  b \le k_1,
\tag{b} \\
&\lt(r(x(-1)^a y^b h^{k+1-a-b}) v_\Lambda) =
x(-1)^a y^b  \quad \text{if} \quad  a > k_0 \ \text{or} \  b
> k_1,
\tag{c} \\
&\lt(r(h(-1)^a x(-1)^b y^{k+1-a-b}) v_\Lambda) =
h(-1)^a x(-1)^b y^{k+1-a-b}, \tag{d} \\
&\lt(r(y(-1)^a h(-1)^b y^{k+1-a-b}) v_\Lambda) =
y(-1)^a h(-1)^b y^{k+1-a-b}. 
\endalign$$
\endproclaim
\demo{Proof} First recall that
$$
x\cdot y^p = [x,y^p] = p y^{p-1} (h-(p-1)).
$$
For $r=r(y(-1)^a y^{k+1-a})$ we have
$$\gathered
(x^q \cdot r)v_\Lambda = x^q(r v_\Lambda) = \\
= \sum^q _{p=0} \binom{q}{p}(x^{q-p} \cdot y(-1)^a) (x^p
\cdot y^{k+1-a}) v_\Lambda \\
= \sum^q _{p=0} \binom{q}{p}(x^{q-p} \cdot y(-1)^a)\, p!\,
\binom{k+1-a}p \\
\cdot y^{k+1-a-p}(k_1 -(k-a))(k_1-(k-a)+1) \cdots \\
\cdots (k_1 - k-1+a+p) v_\Lambda,
\endgathered \tag{10.2}$$
where $y^{k+1-a-p} v_\Lambda =0$ for $k+1-a-p < 0$ and
$x^{q-p} \cdot y(-1)^a=0$
for $q-p > 2a$.

(a) Let $q=k+1+a+b=2a+k+1-a+b$, $b > 0$. Then $k+1-a > 0$
and
$x^p y^{k+1-a} v_\Lambda =0$ for $p \ge k+1-a+b$, and from
(10.2) 
we get $(x^q\cdot r)v_\Lambda = 0$.

(b) and (c). Let $q=k+1+a-b$. The leading term of $(x^q\cdot
r)v_\Lambda$ will be a
 nontrivial summand in (10.2) for $k+1-a-p$ maximal
possible, i.e. minimal
possible $p \ge q - 2a$. For $p=q - 2a = k+1-a-b$ the
interval in
$\Bbb Z$
$$
[k_1-(k-a), k_1-k-1+a+p] = [a - k_0, k_1 - b]
$$
does not contain zero for $a-k_0 > 0$ or $k_1-b < 0$, and
(c) follows.
Similarly, for $p \ge q - 2a$ we have
$k_1-k-1+a+p \ge k_1 - b$. Hence $(x^q\cdot r)v_\Lambda = 0$
if
$0 \in [a-k_0, k_1 - b]$ and (b) follows. 

(d) is obvious.\qed
\enddemo

Clearly $\lt(\bar R v_\Lambda)$ consists of elements
$\pi_{m\alpha}(n) = \lt (r_{m\alpha}(n))$ for 
$n \le -k-1$, $|m| \le k+1$, and the elements listed in (c)
and (d) of Proposition 10.2.
We define a map
$r \: \lt(\bar{R} v_\Lambda) \rightarrow \bar R$, $ \rho
\mapsto r(\rho)$, by
$$
r(x(-1)^a y^b) = r(x(-1)^a y^b h^{k+1-a-b}) \tag{10.3}
$$
for $a+b \le k$, $a > k_0$ or $b > k_1$, and as previously
defined $r$ on the set
$\lt(\bar{R} v_\Lambda) \cap \,\lt(\bar{R})$.
We shall say that the leading term $\rho = x(-1)^a y^b \in
\lt(\bar R v_\Lambda)$, 
$a+b \le k$, is {\it short}. We shall also say that
$\rho \in \lt(\bar R v_\Lambda)$, $\ell (\rho)=k+1$, is {\it
long}. 
Obviously the long $\rho \in \lt(\bar R v_\Lambda)$ are
characterized by
 equality $\lt(r(\rho)) = \lt(r(\rho)v_\Lambda)$. 

By Theorem 9.1 we have the following:

\proclaim{Lemma 10.3} For two embeddings $\rho_1, \rho_2
\subset \pi$, both
$\rho_1$ and $\rho_2$ long, (10.1) holds.
\endproclaim

Relations (8.4) and (8.5) constructed in Section 8 can be
written in the form
$$\gather
\sum_{j<0}((k\!+\!2)j\!-\!n)\{x(j)r_{(i-1)\alpha}(n\!-\!j)v_\Lambda -h(j)
r_{i\alpha}(n\!-\!j)v_\Lambda -
y(j) r_{(i+1)\alpha}(n\!-\!j)v_\Lambda\} \\
-n(k_0 +1-i)r_{i\alpha}(n)v_\Lambda + n y r_{(i+1)\alpha} (n)v_\Lambda = 0,
\endgather$$
$$\gather
\sum_{j<0}\{(k\!+\!2\!-\!i)x(j)r_{(i-1)\alpha}(n\!-\!j) +i
h(j)r_{i\alpha}(n\!-\!j) +
(k\!+\!2\!+\!i)y(j)r_{(i+1)\alpha}(n\!-\!j)\}v_\Lambda \\ 
+[(k+2)n+i(k_1 +1+i)]r_{i\alpha}(n)v_\Lambda + (k+2+i)y r_{(i+1)\alpha} (n)
v_\Lambda = 0.
\endgather$$
We shall also use linear combinations of (8.4) and (8.5) 
for $-k-1 \le n \le 0$ of the form
$$\gather
C_\Lambda r_{i\alpha}(n)v_\Lambda + C_y y r_{(i+1)\alpha} (n)v_\Lambda 
\tag{10.4} \\
+C^x x(-1)r_{(i-1)\alpha}(n+1)v_\Lambda + C^h h(-1)r_{i\alpha}(n+1)v_\Lambda
+
C^y y(-1)r_{(i+1)\alpha}(n+1)v_\Lambda  \\
+\sum_{j\le -2} a_j x(j)r_{(i-1)\alpha}(n-j)v_\Lambda +b_j h(j)r_{i\alpha}
(n-j)v_\Lambda+
c_j y(j)r_{(i+1)\alpha}(n-j)v_\Lambda = 0.
\endgather$$
\proclaim{Lemma 10.4} For two embeddings $\rho_1, \rho_2
\subset \pi$, both
$\rho_1$ and $\rho_2$ short, (10.1) holds.
\endproclaim

\demo{Proof}
Case 1. Let $\rho_1 = x(-1)^a y^b$, $a,b \ge 0$, $a + b \le
k$,
$\rho_2 = x(-1)\rho_1$.

Then $r(\rho_1) = r_{(i-1)\alpha} (n+1)$ for $n+1= -a$, 
$i-1 = a - b$. We can write (10.4) in the form
$$\gather
C r(x(-1)^{a+1} y^b)v_\Lambda + C_y y r(x(-1)^{a+1}
y^{b-1})v_\Lambda \\
+C^x x(-1)r(x(-1)^a y^b)v_\Lambda +C^h h(-1)r(x(-1)^a
y^{b-1})v_\Lambda  \\
+C^y y(-1)r(x(-1)^a y^{b-2})v_\Lambda + A = 0,
\endgather$$
where $A$ denotes a sum of terms with leading terms $\pi_1$,
$\ell(\pi_1) \le a+b$. Now we choose $C^x = 1$, $C_y = 0$
and we get
$$
x(-1)r(x(-1)^a y^b)v_\Lambda \in C r (x(-1)^{a+1}
y^b)v_\Lambda +
M^1_{(x(-1)^{a+1}y^b)}. \tag{10.5}
$$
Clearly (10.5) implies $C \ne 0$, so (10.1) holds.

Case 2. Let  $\rho_1 = x(-1)^a y^b$, $a,b \ge 0$, $a + b \le
k$,
$\rho_2 = y\rho_1$.

Then $r(\rho_1) = r_{(i+1)\alpha} (n)$ for $n= -a$, 
$i+1 = a - b$. We can write (10.4) in the form
$$\gather
C r(x(-1)^a y^{b+1})v_\Lambda + C_y y r(x(-1)^a
y^b)v_\Lambda \\
+C^x x(-1)r(x(-1)^{a-1} y^{b+1})v_\Lambda + C^h
h(-1)r(x(-1)^{a-1} y^b)v_\Lambda  \\
+C^y y(-1)r(x(-1)^{a-1} y^{b-1})v_\Lambda + A = 0,
\endgather$$
where $A$ denotes a sum of terms with leading terms $\pi_1$,
$\ell(\pi_1) \le a+b$. Now we choose $C_y=1$, $C^x=0$ and we
get
$$
y r(x(-1)^a y^b)v_\Lambda \in C r (x(-1)^a y^{b+1})v_\Lambda
+
M^1_{(x(-1)^a y^{b+1})}. \tag{10.6}
$$
Clearly (10.6) implies $C \ne 0$, so (10.1) holds.

In the general case when $\rho_2=x(-1)^{a'} y^{b'}$, $a',b'
\ge 0$,
$a'+b' \le k$ we can find a sequence $\rho_p=x(-1)^{a_p}
y^{b_p}$,
$p=1,\dots,s$, so that $(a,b)=(a_1,b_1)$,
$(a',b')=(a_s,b_s)$, and that
for $\rho_p, \rho_{p+1}$ we can apply either Case 1, Case 2
or Lemma 10.3.
\qed 
\enddemo

It remains to prove (10.1) for two embeddings $\rho_1,
\rho_2 \subset \pi$,
where one $\rho_i$ is short, the other long. First we have

\proclaim{Lemma 10.5} Let $\rho_1,\rho_2 \in \lt(\bar R
v_\Lambda)$,
$\pi = \rho_1 \rho_2$. Then
$$
u(\rho_1 \subset \pi)v_\Lambda \in \Bbb F^\times u(\rho_2
\subset \pi)v_\Lambda
+M^1 _{(\pi)}.
$$
\endproclaim

\demo{Proof} Because of Lemmas 10.3 and 10.4 it is enough to
consider the case when 
$\rho_1$ is long, $\rho_2$ is short. In this case the proof
is similar to the
proof of Lemma 9.3. \qed
\enddemo

We shall say that  $\rho_1,\rho_2 \in \lt(\bar R v_\Lambda)$
are
$\{y,x(-1)\}$-{\it linked} if there exists a sequence of
embeddings
$\sigma_1,\dots,\sigma_p \subset \rho_1 \cup \rho_2$ such
that 
$\sigma_1=\rho_1$, $\sigma_p=\rho_2$ and 
$$\gather
\text{either}\quad   \ell (\sigma_{i-1}\cup\sigma_i) < \ell
(\rho_1 \cup\rho_2)
\ \text{for each}\   i=2,\dots,p \tag{10.7}\\
\text{or}\quad  (\sigma_r \cup \sigma_{r+1})/\sigma_r
\subset \{y,x(-1)\},\  
(\sigma_r \cup \sigma_{r+1})/\sigma_{r+1} \subset
\{y,x(-1)\}. \tag{10.8}
\endgather$$
Note that the assumption (10.8) together with Lemma 10.3,
(10.5) or (10.6) imply for
$\pi = \sigma_r \cup \sigma_{r+1}$ that
$$
r(\sigma_r \subset \pi)v_\Lambda \in
 \Bbb F^\times r(\sigma_{r+1} \subset \pi)v_\Lambda+M^1
_{(\pi)}. \tag{10.9}
$$
By modifying the proof of Lemma 7.3.1 we get the following:
\proclaim{Lemma 10.6} Let $\rho_1,\rho_2 \in \lt(\bar R
v_\Lambda)$,
$\rho_1$ short, $\rho_2$ long, $\rho_1 \cap \rho_2 \ne
\varnothing$. Then
$\rho_1$ and $\rho_2$ are not $\{y,x(-1)\}$-linked only in
the following
cases
$$
\align
\rho_1 &=x(-1)^a y^b,\quad  
\rho_2 =h(-1)^{k+1-a-b} x(-1)^a y^b,\tag{4}\\
&\quad\qquad\qquad a+b \le k, \\
\rho_1 &=x(-1)^b y^a,\quad 
\rho_2 =y(-1)^{k+1-a} y^a, \tag{5}\\
&\quad\qquad\qquad a \ge 1,a+b \le k,\\ 
\rho_1 &=x(-1)^b y^a,\quad 
\rho_2 =y(-1)^{k+1-a-c} h(-1)^c y^a, \tag{7}\\
&\quad\qquad\qquad a\ge 1,  c \ge 1, a+b \le k,\\
\rho_1 &=x(-1)^b y^c,\quad 
\rho_2 =x(-2)^{k+1-a-b} h(-1)^a x(-1)^b,\tag{11}\\
&\quad\qquad\qquad b \ge 1,  a+b \le k, b+ c \le k.
\endalign$$
\endproclaim
We shall say that the cases listed in Lemma 10.6 are {\it
exceptional cases}.

It is clear that if (10.1) holds for $\rho_1,\rho_2 \subset
\pi$, then it
holds as well for $\rho_1,\rho_2 \subset \kappa\pi$ for any 
$\kappa \in \Cal P (\bar{B}_-)$. In particular it is enough
to prove (10.1)
for any two embeddings $\rho_1,\rho_2 \subset \rho_1 \cup
\rho_2$. By 
Lemmas 10.3 and 10.4 it is enough to show (10.1) for 
$\rho_1,\rho_2 \subset \rho_1 \cup \rho_2$, where one
$\rho_i$ is short, the other
one is long. The case $\rho_1 \cap \rho_2 = \varnothing$ is
proved by Lemma 10.5. From 
(10.9) we see that (10.1) holds in the case that a short
$\rho_1$ and a long $\rho_2$
satisfy (10.8). Moreover, if (10.1) holds for all embeddings
of length
less than $\ell (\rho_1\cup\rho_2)$, then (10.1) holds when
$\rho_1$ and $\rho_2$
 are $\{y, x(-1)\}$-linked. What is left is to prove (10.1)
in the exceptional cases. 
We shall prove this by double induction on degree $0, -1,
-2,\dots$ of
 $\rho_1 \cup \rho_2$ and, for a fixed degree, on length
$\ell (\rho_1\cup\rho_2)$.

We shall use the following modification of Lemma 9.4:

\proclaim{Lemma 10.7} Let
$$
u(\rho \subset \pi)v_\Lambda = \sum
C_{\rho',\pi'}\,u(\rho'\subset \pi')v_\Lambda,
$$
where for each $\pi' \prec \pi$ the term $u(\rho' \subset
\pi')$
appears for only one embedding  $\rho' \subset \pi'$. 
Then $C_{\pi'} = C_{\rho',\pi'}=0$
for  $\pi' \prec \pi$. In particular
$$
u(\rho \subset \pi)v_\Lambda \in \sum
C_{\rho',\pi}\,u(\rho'\subset \pi)v_\Lambda
+ M^1_{(\pi)}.
$$
\endproclaim

At some places in the proof it will be convenient to write 
$a \underset {-p}\to \sim \ b$ if parts of $\lt(a)$ and
$\lt(b)$
are of degree $0, -1,\dots, -p$, and $a - b$ is a linear
combination of
terms of the form $u(\rho \subset \pi)v_\Lambda$, where
$\pi$ has at least one
part of degree $\le - p - 1$.
\bigskip
\noindent{\bf Now we consider all exceptional cases:}
\bigskip
\noindent{\bf (4\,)} Let
$\pi = h(-1)^{k+1-a-b} x(-1)^a y^b h^{k+1-a-b}$, $a+b \le
k$,
and note that 
$$
r(\rho_1)=r(x(-1)^a y^b h^{k+1-a-b}).
$$
By using in $k+1-a-b$ steps the relations (10.4) of the form
$C^h=1$, 
$C^x=0$ we get
$$
h(-1)^{k+1-a-b}r(\rho_1)v_\Lambda \underset \,{-1}\to \sim \
\sum C_{cd}\, y(-1)^c y^d r(\rho)v_\Lambda.
$$
Consider two classes of terms with a leading term $\pi_1$ :
(A)
$y(-1)^c \subset \pi_1$ for some $c > 0$ and (B) terms of
the form
$y^d r(\rho)v_\Lambda$. Recall that by Proposition 10.2(a)
for the term 
$\dots r(\rho)v_{\Lambda}$ we may assume that $x$ is not a
part of $\rho$.
By Lemma 9.5(4) in the case (B) we have either
$$\align
&\rho=h(-1)^{k+1-b-p} x(-1)^p y^b, \quad p >
0,\quad\text{or}\\
&\rho=x(-1)^{k+1-b} y^b, \quad \pi_1=x(-1)^{k+1-b}
y^{k+1-a}.
\endalign$$
 In either case there is only one embedding.

In the case (A) we may apply the induction hypothesis on
$\pi_1/y(-1)^c$ and
rewrite the terms in the class (A) so that each $\pi_1$
appears with only
one embedding.

Finally we should consider the terms $u(\rho \subset
\pi_1)v_\Lambda$,
where $\pi_1$ contains a part of degree $\le -2$, say
$b(j)$. Since
$\vert \pi_1 \vert \ge -k -1 $, the part $b(j)$ cannot be a
part of any
embedding $\rho' \subset \pi_1$. Hence we may apply
induction hypothesis on 
$\pi_1/b(j)$ and rewrite a combination of such terms so that
each $\pi_1$
appears with only one embedding.

Now we can apply Lemma 10.7 and (10.1) follows.

\remark{Remark} Note that the proof for the case (4) is only
a slight
modification of the argument given in the exceptional case
(4) in Section 9.
With some other exceptional cases we have a similar
situation as well.
\endremark
\bigskip
\noindent{\bf (5\,)} Set $\pi =y(-1)^{k+1-a} x(-1)^b y^a
h^{k+1-a-b}$,
$a \ge 1$, $a+b \le k$, and note that
$$
r(\rho_1) = r(x(-1)^b y^a h^{k+1-a-b}).
$$

For $\sh \rho=(-1)^{k+1-p} 0^{\,p}$, $p \ge 1$, we shall
express the
 term $y(-1)r(\rho)$ by using the relation (10.4) of the
form $C^y=1$,
$C_y=0$. Note that
$$
y(-1)r(\rho)v_\Lambda \underset \,{-2}\to \sim \
C r(\rho')v_\Lambda + \sum C_{j,X}\,
X(j)r(\rho_{j,X})v_\Lambda,
$$
where $X \in\{x,h,y\}$, $j \in \{-1,-2\}$, and that
$\sh\rho'=
(-1)^{k+1-p+1} 0^{\,p-1}$, $\sh\rho_{j,X}=(-1)^{k+1-q} 0^q$,
$q \ge p$. We can use this type of relations to express
$y(-1)^{k+1-a} r(\rho_1)$ by applying the relations on terms
$u(\rho\subset\pi_1)$
with  $\sh\rho=(-1)^{k+1-p} 0^{\,p}$, $p \ge 1$, and no
parts of $\pi_1$ of
degree $\le -2$ (note that then $p \le k+1-b$):

$$\gather
y(-1)^{k+1-a}r(\rho_1)v_\Lambda  \\
\underset \,{-2}\to \sim \ 
 \sum C_{c\dots h}\, y(-2)^c h(-2)^d x(-2)^e y(-1)^f h(-1)^g
x(-1)^h r(\rho)v_\Lambda.
\endgather$$

For terms with $c+d+e > 0$ we may apply the induction
hypothesis and
rewrite their sum in such a way that for each leading term
$\pi'$ there is only
one embedding $\rho \subset \pi'$ which contributes with a
term 
$u(\rho \subset \pi')v_\Lambda$.

Now note that for the term $u v_\Lambda$ for which $\sh
\rho=
(-1)^{k+1}$ we have (by the way it was obtained)
$\ell(\lt(u v_\Lambda)) \le k+1-a+b < k+1+b$,
and hence $\lt(u v_\Lambda) \succ \lt(y(-1)^{k+1-a}
r(\rho_1)v_\Lambda)$.
Since we want to apply Lemma 10.7, we need not consider
further such terms. 
What is left to consider is terms with  $\sh
\rho=(-1)^{k+1-p} 0^{\,p}$,
$p \ge 1$. Then either $f=0$ or $f > 0$ and 
$\rho=y(-1)^{k+1-p} y^p$.

For a term $u v_\Lambda$, $u=h(-1)^g x(-1)^h r(\rho)$, $\wt
u=
\wt \pi= -(k+1)+b$, we must have (since $p \le k+1-b$)
$$
\rho = y(-1)^q h(-1)^{k+1-p-q} y^p \quad, \quad q \ge 0.
$$

Finally note that for terms omitted when writting
 $\underset \,{-2}\to \sim \ $, i.e. terms of the form 
$u(\rho \subset \pi')v_\Lambda$, where $\pi'$ has at least
one part of the form
$X(j)$, $j \le -3$, $X \in \{x,h.y\}$, we may apply the
induction
 hypothesis (since $X(j)$ cannot be a part of $\rho$) and
rewrite their sum in
such a way that for each leading term $\pi'$ there is only
one embedding 
$\rho \subset \pi'$ which contributes with a term
 $u(\rho \subset \pi')v_\Lambda$.

Hence we can apply Lemma 10.7.
$$
\align
&\rho_1 =x(-1)^b y^a, \quad a \ge 1,a+b \le k, \tag"{\bf
(7\,)}"\\
&\rho_2 = y(-1)^{k+1-a-c} h(-1)^c y^a, \quad c \ge 0.
\endalign
$$
Note that we have already considered the case $c=0$ in (5).
We keep the 
shape 
$$
(-1)^{k+1+b-a} 0^a
$$
fixed and prove (10.1) by induction on $c$, $0 \le c \le
k+1-a$. By the induction
hypothesis we have
$$
y(-1)^{k+2-a-c} h(-1)^{c-1} r(x(-1)^b y^a) v_\Lambda \in
$$
$$
\in \Bbb F^\times x(-1)^b r(y(-1)^{k+2-a-c} h(-1)^{c-1} y^a)
v_\Lambda
$$
$$
+ M^1_{(y(-1)^{k+2-a-c} h(-1)^{c-1} x(-1)^b y^a)}.
$$
We get a relation (10.1) for $\pi=\rho_1 \cup \rho_2$ by
adjoint action of $x$.

\bigskip

\noindent {\bf (11\,)}
\  Let $a\ge 0$, $b \ge 1$, $c\ge 0$, $a+b \le k$,  $b+c\le
k$, 
$$
\align
&\pi  = x(-2)^{k+1-a-b} h(-1)^a x(-1)^b y^c h^{k+1-b-c}, \\
&\rho_2  = x(-2)^{k+1-a-b} h(-1)^a x(-1)^b,
\endalign
$$
and note that $r(\rho_1) = r(x(-1)^b y^c h^{k+1-b-c})$.

First we consider the case $a=0$: the relation (9.20) gives
$$\gathered
y(-2)^{k+1-b} r(y(-1)^b y^{k+1-b}) v_\Lambda
 \underset \,{-2}\to\sim  y^{k+1-b} r(y(-2)^{k+1-b}
y(-1)^b)v_\Lambda\\
+\sum \Sb d+e+f+k+1=\ell (\pi)\\ d+e < k+1 \endSb  
C_{def}\,y(-2)^d y(-1)^e y^f r(y(-1)^p y^{k+1-p})v_\Lambda.
\endgathered\tag{10.10}
$$ 
By adjoint action of $x$ on (10.10) we get
$$\gathered
x(-2)^{k+1-b} r(h(-1)^p x(-1)^{b-p} y^{k+1-b}) v_\Lambda\\
 \underset \,{-2}\to\sim  y^{k+1-b} r(x(-2)^{k+1-b} h(-1)^p
x(-1)^{b-p})v_\Lambda\\
+\sum \Sb d+\dots +i< k+1 \endSb  
C_{d\dots j}\,y(-2)^d h(-2)^e x(-2)^f y(-1)^g h(-1)^h
x(-1)^i y^j r(\rho)v_\Lambda
\endgathered\tag{10.11}
$$
for $0 \le p \le b$, and 
$$\gathered
x(-2)^{k+1-b} r(x(-1)^b y^c) v_\Lambda
 \underset \,{-2}\to\sim  y^c r(x(-2)^{k+1-b}
x(-1)^b)v_\Lambda\\
+\sum \Sb d+\dots +i< k+1 \endSb  
C_{d\dots j}\,y(-2)^d h(-2)^e x(-2)^f y(-1)^g h(-1)^h
x(-1)^i y^j r(\rho)v_\Lambda,
\endgathered\tag{10.12}
$$ 
where for each summand $u(\rho\subset\pi_1)$ in (10.11) and
(10.12) $\rho$ has no 
parts of degree $-2$, so we may apply the induction
hypothesis on 
$\pi_1/y(-2)^d h(-2)^e x(-2)^f$ and then apply Lemma 10.7.
In particular (10.1) holds
in the case $a=0$.

Now we prove the case $a \ge 1$:
By using in $a$ steps relations (10.4) of the form $C^h = 1,
\ C^x = 0$ we get
$$
 h(-1)^a r(\rho_1) \underset \,{-1}\to\sim \sum C_{de} \,
y(-1)^d  y^e  r(\rho),
 \tag{10.13}
$$
where $0 \le d \le a$.
Multiply both sides of (10.13) with $x(-2)^{k+1-a-b}$; in
lack of better
notation we will write this as
$$
x(-2)^{k+1-a-b}\bigl(
 h(-1)^a r(\rho_1) \underset \,{-1}\to\sim \sum C_{de} \,
y(-1)^d  y^e  r(\rho)\bigr ).
 \tag{10.14}
$$

Consider four disjoint classes of $\pi_1 \prec
x(-2)^{k+1-a-b} h(-1)^a x(-1)^b y^c $; in 
the case $\sh \pi_1 \subset \pi$ three classes:
(A) $y(-1)^d \subset \pi_1$
for some $d \ge 1$, (B) $\pi_1$ is not in A and 
$y(-2)^p h(-2)^q x(-2)^{k+1-a-b-p-q} \subset \pi_1$ for some
$p+q \ge 1$ and (C) $\pi_1$ is not in A and 
$ x(-2)^{k+1-a-b} \subset \pi_1$, and the fourth class 
(D) $\sh \pi_1 \nsubseteq \pi$, i.e. $\pi_1$ has either a
part of degree $\le -3$
or more than $k+1-a-b$ parts of degree $-2$.

Let $\pi_1$ be in the class A.

For two embeddings $\rho, \rho' \subset \pi_1$ (i.e.
 $\rho, \rho' \subset \pi_1/ y(-1)^d$ by Lemma 9.5), by
induction hypothesis on
degree, terms $u(\rho \subset \pi_1)$ and  $u(\rho' \subset
\pi_1)$
are proportional modulo higher terms (with leading terms
$\succ \pi_1$). These
higher terms are in class A or D.

Let $\pi_1$ be in the class B. 

For two embeddings  $\rho, \rho' \subset \pi_1$
(  $\rho, \rho' \subset \pi_1 / y(-2)^p h(-2)^q
x(-2)^{k+1-a-b-p-q}$
because of Lemma 9.5), by induction hypothesis on degree,
terms 
 $u(\rho \subset \pi_1)$ and  $u(\rho' \subset \pi_1)$ are
proportional
modulo higher terms (with leading terms $\succ \pi_1$).
These higher terms are
in class B or D.

Let $\pi_1$ be in the class D.

Then for $\rho \subset \pi_1$ we have that $\sh \rho \subset
(-1)^{k+1}0^{k+1}$.
As above, for two embeddings  $\rho, \rho' \subset \pi_1$, 
by induction hypothesis on degree, terms 
 $u(\rho \subset \pi_1)$ and  $u(\rho' \subset \pi_1)$ are
proportional
modulo higher terms (with leading terms $\succ \pi_1$).
These higher terms are
again in class D.

Now consider a term in (10.14) which has the leading term
$\pi_1$ in the class
C. Then $d=0$. For $\rho \subset \pi_1$ (i.e.
$\rho \subset \pi_1/y^e$) the term $u(\rho \subset \pi_1)$
is (by using the relation of the form (10.11)) proportional
to 
$$
 y^{e'}r(x(-2)^{k+1-a-b} h(-1)^{a'}x(-1)^{b'}) 
$$
modulo terms with leading terms $\pi' \succ \pi_1$.

Hence in a finite number of steps we can rewrite (first the
terms
in class C, and then the terms in class A, B and D)
the right hand side of (10.14), modulo
$M^1_{(x(-2)^{k+1-a-b} h(-1)^a x(-1)^b y^c)}$, 
as a linear
combination of terms  $u(\rho \subset \pi_1)$, where for
each $\pi_1\prec\pi$ there
 is only one embedding. By applying Lemma 10.7  we get
(10.1).

\heading{11. Some combinatorial identities of
Rogers-Ramanujan type}
\endheading

\subhead{11.1. Lepowsky-Wakimoto product
formulas}\endsubhead

Let $A=(a_{ij})$ be an $(\ell + 1) \times (\ell + 1)$
generalized Cartan matrix
(GCM) so that the corresponding Kac-Moody Lie algebra $\goth
g(A)$ is an
affine Lie algebra. Use the usual notation for roots
 $\{\alpha_0, \alpha_1,\dots,\alpha_\ell \} \subset \goth
h^*$ and
coroots  $\{\alpha_0 ^\vee, \alpha_1 ^\vee,\dots,\alpha_\ell
^\vee\} \subset
 \goth h$, $a_{ij} = \langle \alpha_i ^\vee, \alpha_j
\rangle$. Let 
$\bold s=(s_0, s_1, \dots, s_\ell) \in (\Bbb
Z_{>0})^{\ell+1}$ and let 
$\rho_{\bold s} ^\vee \in \goth h$ be an element such that 
$\langle \rho_{\bold s} ^\vee,\alpha_i \rangle =
s_i$ for $i=0,1,\dots, \ell$. Then $\{\alpha_0 ^\vee,
\alpha_1 ^\vee,\dots,
\alpha_\ell ^\vee, \rho_{\bold s} ^\vee\}$ is a basis of
$\goth h$. Write $\rho^\vee$ 
instead of $\rho_{\bold s} ^\vee$ in the case $\bold s =
(1,1,\dots, 1)$.

Recall the Weyl-Kac character formula for the standard
$\goth g(A)$-module $L(\Lambda)$ with highest weight
$\Lambda$ (cf. \cite{K}):
$$
e^{-\Lambda} \ch\, L(\Lambda) = N(\Lambda + \rho\, ; A)/D(A),
\tag 11.1.1$$
where
$$\aligned
&N(\Lambda + \rho\, ; A) =
 \sum_{w \in W(A)} \varepsilon (w) e^{w(\Lambda +
\rho)-(\Lambda + \rho)},\\
&D(A)  = \sum_{w \in W(A)} \varepsilon (w) e^{w\rho- \rho}
= \prod_{\alpha > 0} (1 - e^{-\alpha})^{\text{dim}\goth
g(A)_\alpha}.
\endaligned$$
Then after applying the homomorphism
$e^{-\alpha} \mapsto q^{\langle \rho_{\bold s} ^\vee, \alpha
\rangle}$ 
(i.e. by $\rho_{\bold s} ^\vee$-specialization)
we get power series in $q$:
$$
 q^{\langle\Lambda,\rho_{\bold s}^\vee\rangle}
\ch^{\,\rho_{\bold s}^\vee} L(\Lambda)
= N^{\rho_{\bold s} ^\vee} (\Lambda + \rho\, ; A) /
D^{\rho_{\bold s} ^\vee} (A),
\tag 11.1.2$$
where
$$\aligned
&N^{\rho_{\bold s} ^\vee} (\Lambda + \rho\, ; A) = \sum_{w
\in W(A)}
\varepsilon(w) q^{- \langle\rho_{\bold s} ^\vee,w(\Lambda +
\rho)-(\Lambda + \rho)\rangle},\\
&D^{\rho_{\bold s} ^\vee} (A) 
 = \sum_{w \in W(A)}
\varepsilon(w) q^{- \langle\rho_{\bold s} ^\vee,w \rho -
\rho\rangle}
= \prod_{\alpha > 0}
(1 - q^{\langle\rho_{\bold s} ^\vee , \alpha
\rangle})^{\text{dim} \goth g (A)_\alpha}.
\endaligned$$

Note that $\langle\rho^\vee,w(\Lambda \!+\!
\rho)\!-\!(\Lambda \!+\! \rho)\rangle =
\langle\Lambda \!+\!\rho, w\rho^\vee \!-\!\rho^\vee\rangle$
implies Lepowsky's numerator
formula  $N^{\rho^\vee} (\Lambda + \rho\, ; A) = D^{\Lambda
+\rho} ({}^\tau \!A)$ 
(cf. \cite{L}). M.Wakimoto generalized this formula in the
following way (cf. \cite{W}):
Let $\tilde{A} = (\tilde{a}_{ij})$ be another GCM of affine
type with roots
$\{\tilde{\alpha}_0, \tilde{\alpha}_1,\dots ,
\tilde{\alpha}_\ell \} \subset
\tilde{\goth h}^*$ and coroots $\{\tilde{\alpha}_0 ^\vee,
\tilde{\alpha}_1 ^\vee,
\dots , \tilde{\alpha}_\ell ^\vee\} \subset \tilde{\goth
h}$,
$\tilde{a}_{ij} = \langle\tilde{\alpha}_i ^\vee,
\tilde{\alpha}_j \rangle$. Fix
$\tilde{\rho}^\vee, \langle \tilde{\rho}^\vee
,\tilde{\alpha}_i \rangle = 1$ for
$i=0,1,\dots, \ell$.
For fixed $\bold s=(s_0,s_1,\dots, s_\ell)$ define
isomorphism
$\Theta_{\bold s} : \goth h^* \rightarrow \tilde{\goth h}^*$
so that $\langle\tilde{\alpha}_i ^\vee, \Theta_{\bold s}
(\lambda) \rangle =
s_i \langle \alpha_i ^\vee, \lambda \rangle$, 
$\langle\tilde{\rho}^\vee, \Theta_{\bold s} (\lambda)
\rangle =
\langle \rho_{\bold s} ^\vee, \lambda \rangle$. Now assume
that $\tilde{a}_{ij} =
 s_i a_{ij} s_j ^{-1}$. Then it is easy to check that
$\tilde{r}_i =
\Theta_{\bold s} r_i \Theta_{\bold s} ^{-1}$, where
$\tilde{r}_i$ and $r_i$ are reflections
with respect to simple roots $\tilde{\alpha}_i$ and
$\alpha_i$ (respectively).
Hence $w \mapsto \tilde{w} = \Theta_{\bold s} w
\Theta_{\bold s} ^{-1}$ defines an isomorphism
between the Weyl groups $W(A)$ of $\goth g (A)$ and
$W(\tilde{A})$ of
 $\goth g (\tilde{A})$ (or $\goth g ({}^\tau\!\tilde{A})$).
We also have
$\langle \rho_{\bold s} ^\vee , w \lambda \rangle =
\langle\tilde{\rho}^\vee ,
 \tilde{w} (\Theta_{\bold s} \lambda) \rangle $.
This implies the first statement in the following theorem:
\proclaim{Theorem 11.1.1} \cite{W}. If \ 
$\tilde{a}_{ij} = s_i a_{ij} s_j ^{-1}$, \  $i,j =
0,1,\dots, \ell$, \ then
$$\align
&N^{\rho_{\bold s} ^\vee} (\Lambda + \rho\, ; A) =
D^{\Theta_{\bold s} (\Lambda+ \rho)} ({}^\tau\!\tilde{A}),
\tag"{(a)}"\\
&N^{\rho_{\bold s} ^\vee} (\Lambda + \rho\, ; A) =
N^{\Lambda + \rho}(\rho_{\bold s} ^\vee; {}^\tau\!{A}).
\tag"{(b)}"
\endalign$$
\endproclaim

Formula (a) implies that for certain choices of $\bold s$ 
the
$\rho_{\bold s} ^\vee$-specialization of characters of all
standard modules can be
written as infinite products.
The ``duality" formula (b) implies that for certain
(depending
on $\rho_{\bold s} ^\vee$) standard modules
all $(\Lambda + \rho)$-specializations of characters can be
written as infinite products.

Now let $\goth g (A)= \tilde{\goth g}= \goth {sl} (2,\Bbb
F)^{\!\sim}$, i.e.
$A = \left(\smallmatrix 2 & -2\\-2 & 2\endmatrix \right) $.
For $\bold s = (n,n)$,
$(n,2n)$ or $(2n, n)$ Theorem 11.1.1 may be applied. In the
case $\bold s=(n,n)$
we have $\tilde{A} = A$ and (a) is Lepowsky's numerator
formula. In the
case $(n, 2n)$ we have 
$\tilde{A} = \left(\smallmatrix 2 & -1\\-4 & 2\endmatrix
\right)$,
$\goth g({}^\tau\!\tilde{A})$ is of the type $A_2 ^{(2)}$,
$\tilde{\alpha}_1 ^\vee$ is a long root, $\tilde{\alpha}_0
^\vee$ is a short
root. The case $(2n, n)$ is related to $(n, 2n)$ by an outer
automorphism.
For $\bold s=(s_0,s_1)$, $A = \left(\smallmatrix 2 & -2\\-2
& 2\endmatrix \right) $,
and for $\mu \in \tilde{\goth h}^*$, $\mu(\tilde{\alpha}_0
^\vee)= m_0 \in \Bbb Z_{>0}$,
 $\mu(\tilde{\alpha}_1 ^\vee)= 2 m_1 \in \Bbb Z_{>0}$,  
$\tilde{A} = \left(\smallmatrix 2 & -1\\-4 & 2\endmatrix
\right)$, write
$$
P(s_0,s_1; q) = D^{\,\rho_{\bold s}^\vee}(A),\qquad 
Q( m_0, 2 m_1 ; q) = D^{\mu} ({}^\tau\!\tilde{A}).
$$
\proclaim{Lemma 11.1.2} Let  $s_0,s_1,m_0 \in \{ 1, 2, 3,
\dots \}$, 
$m_1 \in \{ \tfrac 12, 1, \tfrac 32, \dots \}$,
$s=s_0 + s_1$,  $m = m_0 + m_1$. Then
\bigskip
\noindent {\rm (a)}\quad
$P(s_0,s_1; q) = \prod_{r \equiv 0, s_0 \mod{s}} (1 - q^r) 
\prod_{r \equiv  s_1 \mod{s}} (1 - q^r)$,
\bigskip
\noindent {\rm (b)}\quad
$Q( m_0, 2 m_1 ; q) =  \prod\Sb r \equiv 0, \pm m_0
\mod{2m}\endSb (1 - q^r)
\prod\Sb r \equiv \pm 2 m_1 \mod{4m}\endSb(1 - q^r)$.
\bigskip
\noindent (Here is assumed $r \in \Bbb Z_{>0}$.)
\endproclaim

For $\bold s=(s_0, s_1)$ and $\Lambda = k_0 \Lambda_0 + k_1
\Lambda_1$ write
$d^{s_0, s_1} _{k_0, k_1} (q) = 
q^{\langle\Lambda,\rho_{\bold s}^\vee\rangle}
\ch^{\,\rho_{\bold s}^\vee} L(\Lambda)$.
Then Theorem 11.1.1 and Lemma 11.1.2 give Lepowsky-Wakimoto
product
formulas for specialized characters of standard modules for
 $\goth {sl}(2,\Bbb F)^\sim$:
$$\align
&d^{1,1} _{k_0,k_1} (q) = P(k_0\!+\!1, k_1\!+\!1 ; q)/
P(1,1; q),
\tag{11.1.3}\\
&d^{1,2} _{k_0,k_1} (q) = Q(k_0\!+\!1, 2(k_1\!+\!1) ; q)/
P(1,2; q),
\tag{11.1.4}\\
&d^{s_0,s_1} _{n-1,2n-1} (q) = Q(n s_0, 2 n s_1 ; q)/
P(s_0,s_1; q),
\tag{11.1.5}\\
&d^{s_0,s_1} _{2n-1,n-1} (q) = Q(n s_1, 2 n s_0 ; q)/
P(s_0,s_1; q),
 \tag{11.1.6}\\
&d^{s_0,s_1} _{n-1,n-1} (q) = P(ns_0,ns_1; q)/ P(s_0,s_1;
q).
  \tag{11.1.7}
\endalign$$
It is easy to see the following:
\proclaim{Proposition 11.1.3} Each of the products
(11.1.3)--(11.1.7) is of the form
$$
\prod_i \prod_{r\in B_i} (1+q^r) \prod_j \prod_{r\in C_j}
(1-q^r)^{-1} \tag 11.1.8
$$
for some congruence classes $B_1,\dots ,B_b \subset \Bbb
Z_{>0}$, $b \ge 0$,
$C_1,\dots ,C_c \subset \Bbb Z_{>0}$, $c \ge 1$. (These
classes need not be disjoint,
some may be empty.)
\endproclaim
\noindent As particular examples we have:
\bigskip
\noindent {\rm (11.1.9)}\quad 
$d^{1,1} _{n-\!1,n-\!1} (q)=
\prod_{r \text{ odd, } r\not\equiv b \mod{2b}} \,(1 -
q^r)^{-1}
 \prod_{r\not\equiv 0, n\! \mod{2n}}\, (1 - q^r)^{-1}$
\bigskip
      ${}\qquad \quad \cdot 
\prod_{1 \le j \le a} \prod_{r \equiv 2^{a-j}b
\mod{2^{a-j+1}b}}\,(1+q^r)$,\quad
for   $n = 2^a b$, $a \ge 0$, $b$ odd,
\bigskip
\noindent {\rm (11.1.10)}\quad 
$d^{1,1} _{n-1,2n-1} (q)=
 \prod_{r \ \text{odd}}\,(1-q^r)^{-1}
 \prod_{r\not\equiv 0, \pm n \mod{3n}}\, (1 - q^r)^{-1}$,
\bigskip
\noindent {\rm (11.1.11)}\quad 
$d^{1,1} _{k_0,k_1} (q)=
 \prod_{r \ \text{odd}}\,(1-q^r)^{-1}
 \prod_{r\not\equiv 0, \pm (k_0+1) \mod{(k+2)}}\, (1 -
q^r)^{-1}$
\bigskip
       ${}\qquad \ \quad {}$   for    $k = k_0 \!+ \!k_1 \ne
2 k_1$,
\bigskip
\noindent {\rm (11.1.12)}\quad
$d^{1,2} _{k_0,k_1} (q)=
 \prod\Sb r \not\equiv 0, \pm (k_0 + 1) \mod{2(k+2)}\\ 
r \not\equiv \pm 2(k_1 + 1) \mod{4(k+2)}\endSb
 \,(1 - q^r)^{-1}$
\bigskip
       ${}\qquad \ \quad {}$   for    $k = k_0 \!+ \!k_1 \ne
3 k_1 + 1$, 
\bigskip
\noindent {\rm (11.1.13)}\quad 
$d^{1,2} _{2n-1,n-1} (q)=
 \prod_{r \equiv \pm n \mod{6n}}\,(1+q^r)
 \prod_{r\not\equiv 0,\pm n, \pm 2n \mod{6n}} \,(1 -
q^r)^{-1}$,
\bigskip
\noindent {\rm (11.1.14)}\quad
$d^{1,2} _{n-1,n-1} (q)=
 \prod_ {r \not\equiv 0 \mod{n}} \,(1 - q^r)^{-1}$.
\bigskip
For nonempty subsets $A_1, A_2, \dots ,A_s \subset \Bbb
Z_{>0}$, $s \ge 2$, set
$
 A = A_1 \bigsqcup A_2 \bigsqcup \dots \bigsqcup A_s,
$
a disjoint union of sets. We may call elements of $\Cal
P(A)$ {\it colored partitions}
with parts in $A$, where for $i \in \{1, \dots ,s\}$ and $j
\in A_i$ we say that
$j$ is of {\it color} $i$ and of {\it weight} $|j| = j \in \Bbb Z_{>0}$.
We define the {\it degree} 
$|\pi|$ of $\pi \in \Cal P(A)$ as before: $|\pi| = \sum_{a
\in A}\,\pi (a)|a|$.
It is clear that infinite products of the form (11.1.8) may
be
interpreted as generating functions of partition functions
for colored
partitions defined by congruence conditions.

\subhead{11.2. Difference conditions and partition
ideals}\endsubhead

Let $A$ be a nonempty set and $\Cal P(A)$ the set of all
partitions 
with parts in $A$. By following \cite{A1} and \cite{A2}
we shall say that $\Cal C \subset \Cal P(A)$ is a
{\it partition ideal} in $\Cal P(A)$ if $\kappa \in \Cal C$ and
$\pi \in \Cal P(A)$ 
implies $\kappa \cap \pi \in \Cal C$, or equivalently, if
$\mu \in \Cal C$, 
$\nu \in \Cal P(A)$, $\nu \subset \mu$, implies $\nu \in
\Cal C$.
If $\Cal C, \Cal I\subset \Cal P(A)$, $\Cal P(A) = \Cal C
\bigcup \Cal I$,
$\Cal C \bigcap \Cal I = \emptyset$, then it is easy to see
that $\Cal C$
is a partition ideal in $\Cal P(A)$ if and only if $\Cal I$
is an ideal in the monoid $\Cal P(A)$. In particular,
$$
\Cal C_\Lambda = \Cal P(\bar B_-)\backslash (\ell t(\bar R
v_\Lambda))
$$
is a partition ideal in $\Cal P(\bar B_-)$. By Theorem 6.5.5
the partition
ideal $\Cal C_\Lambda$ determines the character of standard 
$\goth g(A)$-module $L(\Lambda)$ with highest weight
$\Lambda$:
$$
e^{-\Lambda} \ch\, L(\Lambda) = \sum_{\pi \in \Cal
C_\Lambda}
e^{\wt \pi +|\pi|\delta}.\tag 11.2.1
$$

In Section 6.6 we have described the set $\ell t(\bar R
v_\Lambda)$, and as
a consequence we may describe the partition ideal $\Cal
C_\Lambda$,
$\Lambda = k_0 \Lambda_0 + k_1 \Lambda_1$, $k=k_0 + k_1$, as
the
set of all colored partitions $\pi$ in $\Cal P(\bar B_-)$
satisfying
{\it difference conditions}:
$$
\gathered
\pi (y(j-1)) + \pi (h(j-1)) + \pi (y(j))  \le k,\\
\pi (h(j-1)) + \pi (x(j-1)) + \pi (y(j))  \le k,\\
\pi (x(j-1)) + \pi (y(j)) + \pi (h(j)) \le k,\\
\pi (x(j-1)) + \pi (h(j)) + \pi (x(j))  \le k,
\endgathered \tag{11.2.2}
$$
and {\it initial conditions}:
$$
\pi (x(-1)) \le k_0, \quad \pi (y(0))  \le k_1.\tag 11.2.3
$$

For $\Bbb B \subset \Bbb Z_{> 0}$ we shall write an ordinary
 partition $f \: \Bbb B \rightarrow \Bbb N$
(with parts in $\Bbb B$) by its values $(f_r \mid r \in \Bbb
B)$.
For $\bold s = (s_0, s_1)$, $s_0, s_1 \ge 1$, $s_0 \neq
s_1$, $s=s_0+s_1$, 
we have a bijection $\bar{B}_- \rightarrow \Bbb B^{s_0,
s_1}$, where
$$
\Bbb B^{s_0,s_1}=
 \{r \in \Bbb Z_{ > 0} \mid r \equiv 0, \pm s_1
\mod{s}\},\tag 11.2.4
$$
defined by the specialization
$$
\alignedat2
&y(-j) \mapsto \deg y(-j)  = s_1 + js, && \quad
 j \in \Bbb Z_{ \ge 0}, \\
&h(-j) \mapsto \deg h(-j)  = js, && \quad
 j \in \Bbb Z_{> 0},  \\
&x(-j) \mapsto \deg x(-j)  = -s_1 + js, && \quad
 j \in \Bbb Z_{> 0}.  
\endalignedat \tag{11.2.5}
$$
This bijection extends to an isomorphism of monoids
$\Cal P (\bar{B}_- )\rightarrow \Cal P (\Bbb B^{s_0,s_1})$.
Hence we have a bijection of partition ideals
$$
\Cal C_\Lambda \rightarrow \Cal C^{s_0,s_1} _{k_0,k_1}\, ,
$$
where $\Cal C^{s_0,s_1} _{k_0,k_1} \subset \Cal P (\Bbb
B^{s_0,s_1} )$ is the set of all
partitions $f \: \Bbb B^{s_0,s_1} \rightarrow \Bbb N$ which
satisfy {\it difference
conditions}:
$$
\gathered
f_{js+s_1} + f_{js} + f_{js-s_0} \le k,\\
f_{js} + f_{js-s_1} + f_{js-s_0}  \le k,\\
f_{js+s_0} + f_{js+s_1} + f_{js} \le k,\\
f_{js+s_0} + f_{js} + f_{js-s_1}  \le k,
\endgathered \tag{11.2.6}
$$
and {\it initial conditions}:
$$
f_{s_0} \le k_0, \quad f_{s_1} \le k_1.\tag 11.2.7
$$
As a particular example in the case $\bold s=(1,2)$ we have 
$\Bbb B^{1,2} = \Bbb Z_{> 0}$ and the partition ideal $\Cal
C^{1,2} _{k_0,k_1}$
in $\Cal P (\Bbb Z_{>0})$ is defined by
$$
\gathered
f_{3j+2} + f_{3j} + f_{3j-1}  \le k,\\
f_{3j} + f_{3j-1} + f_{3j-2}  \le k,\\
f_{3j+2} + f_{3j+1} + f_{3j}  \le k,\\
f_{3j+1} + f_{3j} + f_{3j-2}  \le k,\\
f_1 \le k_0, \quad f_2 \le k_1 .
\endgathered \tag{11.2.8}
$$

For a nonempty set $A$ we shall write a partition $\pi \: A
\rightarrow \Bbb N$
(with parts in $A$) by its values $(\pi_a \mid a \in A)$.
Now consider the case of
principal specialization, i.e. $\bold s = (1,1)$. Set 
$$
A_1 = \{ 1,2,3, \dots \} \cong \Bbb Z_{>0},\quad
A_2 = \{ \underline 1,\underline 2,\underline 3, \dots \} \cong \Bbb
Z_{>0},\quad
A = A_1 \bigsqcup A_2.
$$
We shall say that $j \in A_1$ is of degree $|j| = j \in \Bbb
Z_{>0}$ and of plain color, and
that $\underline j \in A_2$ is of degree $|\underline j| = j \in \Bbb
Z_{>0}$ and of underline color. Set
$
\underline B = A_1 \bigsqcup \,\{ \underline j \in A_2 \mid j \
\text{odd} \}.
$
Then we have a bijection $\bar{B}_- \rightarrow \underline B$
defined by the specialization
$$
\alignedat2
&y(-i) \mapsto 2i+1 \in A_1, && \quad
 i \in \Bbb Z_{ \ge 0}, \\
&h(-i) \mapsto 2i \in A_1, && \quad
 i \in \Bbb Z_{> 0},  \\
&x(-i) \mapsto \underline{2i -1} \in A_2, && \quad
 i \in \Bbb Z_{> 0}.  
\endalignedat \tag{11.2.9}
$$
This bijection extends to an isomorphism of monoids
$\Cal P (\bar{B}_- )\rightarrow \Cal P (\underline B)$.
Hence we have a bijection of partition ideals
$$
\Cal C_\Lambda \rightarrow \Cal C^{1,1} _{k_0,k_1}\, ,
$$
where $\Cal C^{1,1} _{k_0,k_1} \subset \Cal P (\underline B)$ is
the set of all colored
partitions $\pi \: \underline B \rightarrow \Bbb N$ which satisfy
{\it difference
conditions}:
$$
\gathered
\pi_{2i+1} + \pi_{2i} + \pi_{2i-1} \le k,\\
\pi_{2i} + \pi_{\underline {2i-1}} + \pi_{2i-1}  \le k,\\
\pi_{\underline {2i+1}} + \pi_{2i+1} + \pi_{2i} \le k,\\
\pi_{\underline {2i+1}} + \pi_{2i} + \pi_{\underline {2i-1}}  \le k,
\endgathered \tag{11.2.10}
$$
and {\it initial conditions}:
$$
\pi_{\underline 1} \le k_0, \quad \pi_{1} \le k_1.\tag 11.2.11
$$

\subhead{11.3. Combinatorial identities}\endsubhead

Formulas (11.1.1) and (11.2.1) express the character
$\ch\,L(\Lambda)$ 
in two different ways which may have different combinatorial
interpretations.
In particular, if we take a specialization (11.2.5) for $s_0
\ne s_1$, then the
specialized character $d^{s_0, s_1} _{k_0, k_1} (q)$ is a
generating function
for $\Cal C^{s_0,s_1} _{k_0,k_1}$-partition function, where
the  
partition ideal $\Cal C^{s_0,s_1} _{k_0,k_1}$ is defined by
difference
conditions (11.2.6), initial conditions (11.2.7) and
congruence conditions
(11.2.4). On the other side, for certain choices of $\bold
s$ and/or $\Lambda$
we get Lepowsky-Wakimoto product formulas (11.1.4)--(11.1.7)
for
specialized characters $d^{s_0, s_1} _{k_0, k_1} (q)$, and
these  have
a combinatorial interpretation in terms of congruence
conditions. As a consequence
we get a series of Rogers-Ramanujan type combinatorial
identities. In the
case of principal specialization (11.2.9) (i.e. for $\bold s
= (1,1)$) the
product formula (11.1.3) on one side and the conditions
(11.2.10)--(11.2.11) on
the other side give a series of Rogers-Ramanujan type
combinatorial identities
for colored partitions.

 For example, in the case $\bold s=(1,2)$ and $k_0 = k_1 =
n-1$ we have a product
formula (11.1.14) for $d^{1,2} _{n-1,n-1} (q)$, so (11.2.8)
for the partition ideal 
$\Cal C^{1,2} _{n-1,n-1}$ in $\Cal P (\Bbb Z_{>0})$  gives
that for every 
$m \in \Bbb N$ the number of partitions $f$ such that $ |f|
= m$,
$$
f_j > 0 \quad \text{implies} \quad j \not\equiv 0 \mod{n} ,
$$
equals the number of partitions $f$ such that $ |f| = m$,
$$
\gathered
f_{3j+2} + f_{3j+1} + f_{3j}  \le 2n -2,\\
f_{3j+2} + f_{3j} + f_{3j-1}  \le 2n - 2,\\
f_{3j+1} + f_{3j} + f_{3j-2}  \le 2n -2,\\
f_{3j} + f_{3j-1} + f_{3j-2}  \le 2n - 2,\\
f_1 \le n - 1, \quad f_2 \le n - 1 .
\endgathered 
$$
Of course, this also equals the number of partitions $f$ of
$m$ such that each part 
of $f$ appears
at most $n-1$ times, and the example given in the
introduction is for the $(1,2)$-specialization
of the level 2 standard module $L(\Lambda_0+\Lambda_1)$.

By some sort of coincidence, the combinatorial identities for
the $(1,2)$-speciali\-za\-tion of the
fundamental $\goth {sl}(2,\C)\sptilde$-modules are identical with
the combinatorial identities
for the level 3 modules in the principal picture for the type
$A^{(2)}_2$ affine Lie
algebra (cf. \cite{C1}, \cite{C2}, \cite{A3}). Note that in
this case (i.e. the level 1
$\goth {sl}(2,\C)\sptilde$-modules) Lemma 9.2 is
(essentially) all that is needed for
the proof of linear independence.

Another example that we may take is  the case $\bold
s=(1,1)$ and $k_0 = 1$, $ k_1 = 2$. 
Then we have a product
formula (11.1.11) for $d^{1,1} _{1,2} (q)$, so
(11.2.10)--(11.2.11) for the partition ideal 
$\Cal C^{1,1} _{1,2}$ in $\Cal P (\underline B)$  gives that for
every 
$m \in \Bbb N$ the number of colored partitions $\pi \in
\Cal P (\underline B)$ 
such that $|\pi| = m$, 
$$
\pi_a > 0 \quad \text{implies} \quad a \in \{\underline i \mid i
\  \text{odd}\}
\bigsqcup \,\{ i \mid i \equiv \pm 1 \mod{5}\} ,
$$
equals the number of colored partitions $\pi \in \Cal P
(\underline B)$ such that $ |\pi| = m$,
$$
\gathered
\pi_{\underline {2i+1}} + \pi_{2i+1} + \pi_{2i} \le 3,\\
\pi_{2i+1} + \pi_{2i} + \pi_{2i-1} \le 3,\\
\pi_{\underline {2i+1}} + \pi_{2i} + \pi_{\underline {2i-1}}  \le 3,\\
\pi_{2i} + \pi_{\underline {2i-1}} + \pi_{2i-1}  \le 3,\\
\pi_{\underline 1} \le 1, \quad \pi_{1} \le 2.
\endgathered $$
In a similar way we get that for every 
$m \in \Bbb N$  the number of colored partitions $\pi \in
\Cal P (\underline B)$ 
such that $|\pi| = m$, 
$$
\pi_a > 0 \quad \text{implies} \quad a \in \{\underline i \mid i
\  \text{odd}\}
\bigsqcup \,\{ i \mid i \equiv \pm 2 \mod{5}\} ,
$$
equals the number of colored partitions $\pi \in \Cal P
(\underline B)$ such that $ |\pi| = m$,
$$
\gathered
\pi_{\underline {2i+1}} + \pi_{2i+1} + \pi_{2i} \le 3,\\
\pi_{2i+1} + \pi_{2i} + \pi_{2i-1} \le 3,\\
\pi_{\underline {2i+1}} + \pi_{2i} + \pi_{\underline {2i-1}}  \le 3,\\
\pi_{2i} + \pi_{\underline {2i-1}} + \pi_{2i-1}  \le 3,\\
\pi_{\underline 1} \le 0, \quad \pi_{1} \le 3.
\endgathered $$

\enddocument